\documentclass[lettersize,journal]{IEEEtran}
\usepackage{amsmath,amsfonts}
\usepackage{algorithmic}
\usepackage{algorithm}
\usepackage{array}
\usepackage[caption=false,font=normalsize,labelfont=sf,textfont=sf]{subfig}
\usepackage{textcomp}
\usepackage{stfloats}
\usepackage{url}
\usepackage{verbatim}
\usepackage{graphicx}
\usepackage{cite}
\usepackage{mathrsfs}
\usepackage{float}
\usepackage{enumitem}
\makeatletter

\newcommand{\Rmnum}[1]{\expandafter\@slowromancap\romannumeral #1@}
\makeatother
\newtheorem{theorem}{Theorem}
\newtheorem{assumption}{Assumption}
\newtheorem{lemma}{Lemma}
\newtheorem{example}{Example}
\newtheorem{remark}{Remark}
\newtheorem{proposition}{Proposition}
\newtheorem{corollary}{Corollary}
\begin{document}
	
	\title{Switching Event-Triggered Control of Nonlinear Parabolic PDE Systems via Galerkin/Neural-Network-Based Modeling Approach}
	
	\author{Xiao-Yu Sun, Chuan Zhang, Huai-Ning Wu, Xian-Fu Zhang  
		\thanks{This work was supported by the National Natural Science Foundation of China (Grant Nos. 62003189 and 61973189); the China Postdoctoral Science Foundation (Grant No. 2020M672024); the Natural Science Foundation of Shandong Province (Grant No. ZR2024QF012). $\textit{(Corresponding author: Chuan Zhang.)}$
			
			Xiao-Yu Sun is with the School of Mathematical Sciences Qufu Normal University, Qufu 273165, China (e-mail: sunxiaoyu156520@163.com).
			
			Chuan Zhang is with the School of Mathematical Sciences Qufu Normal University, Qufu 273165, China (e-mail: zhangchuan@qfnu.edu.cn).
			
			Huai-Ning Wu is with the Science and Technology on Aircraft Control Laboratory, School of Automation Science and Electrical Engineering, Beihang University, Beijing 100191, China, and also with Hangzhou International Innovation Institute, Beihang University, Hangzhou 311115,
			China (e-mail: whn@buaa.edu.cn).
			
			Xian-Fu Zhang is with the School of Control Science and Engineering, Shandong University, Jinan 250061, China (e-mail: zhangxianfu18@126.com).}}

\markboth{Journal of \LaTeX\ Class Files,~Vol.~14, No.~8, August~2021}%
{Shell \MakeLowercase{\textit{et al.}}: A Sample Article Using IEEEtran.cls for IEEE Journals}


\maketitle

\begin{abstract}
	This paper focuses on switching event-triggered output feedback control for a class of parabolic partial differential equation (PDE) systems subject to unknown nonlinearities and external bounded disturbance. Initially, the PDE systems is properly separated into a finite-dimensional ordinary differential equation (ODE) slow system and an infinite-dimensional ODE fast system based on Galerkin technique, especially the slow system can characterize the dominated dynamics. Then, a three-layer neural network is employed to approximate the unknown nonlinearities, and Levenberg-Marquardt algorithm is adopted to get a relative accurate slow system. Subsequentaly, a switching event-triggered control scheme is developed, and a waiting time subject to the triggered condition is implemented to avoid the Zeno behavior and convert the slow system into a switching system. In the following, the stability and $\textbf{\textit{H}}_\infty$ performance issues of the closed-loop system are discussed, and the controllers are displayed in terms of bilinear matrix inequalities (BMIs). Novel algorithms are proposed to convert the BMIs into linear matrix inequalities (LMIs). Additionally, a sub-optimal switching event-triggered controller is obtained using an iterative optimization approach based on LMIs. Finally, simulation on the catalytic rod reaction model and traffic flow model demonstrate the effectiveness of switching event-triggered control strategy.
\end{abstract}

\begin{IEEEkeywords}
	Event-triggered control, parabolic partial differential equation systems, multilayer neural network, unknown nonlinearities, linear matrix inequality.
\end{IEEEkeywords}

\section{Introduction}
\IEEEPARstart{D}{URING} recent years, partial differential equation (PDE) systems have been widely used in many technical domains. The reason for this is that the dynamics of many engineering systems and real-world physical processes exhibit complex intertwinements between the temporal and spatial dimensions, such as heat conduction \cite{4982632,9619762}, fluid flow \cite{6573344,9234521}, elastic vibration \cite{8352583}, and diffusion reaction \cite{6243175}, etc. These complex processes cannot be sufficiently represented by a small set of attributes. However, the theory of PDE provides a more realistic foundation for the development of control approaches and enables accurate modeling and analysis of these processes. It provides strong tools for handling complex system control issues, allowing for improved control \cite{5406054}, reliable control \cite{7161345}, and optimal control \cite{10759746}, and so on. Existing researches have already recognized the critical role of PDE systems in industrial production.

Beginning with the idea of approximation inertial manifolds \cite {1997Finite}, the Galerkin method is fast becoming a key instrument in the research of parabolic PDE systems. The eigenspectrum of the spatial differential operator in parabolic PDE systems can be distinctly characterized by a finite-dimensional, slowly evolving component, complemented by an infinite-dimensional, stable, and rapidly decaying segment. This property makes Galerkin method capable of lowering the complexity of parabolic PDE systems by reducing them to ordinary differential equation (ODE) systems \cite{9873976,6979213}. Inspired by previous works \cite{9269405,7885587}, this paper employs the Galerkin method to conduct a finite-dimensional modeling approach for PDE systems, thereby the control design can be relatively finely tuned to target the primary dynamic characteristics of the system.

Over the past few years, event-triggered control has attracted a lot of attention as a way to reduce computational and communication resources. Traditional time-triggered control has achieved many fruitful results during the past decades, for example, active fault-tolerant control of sampled-data nonlinear PDE systems \cite{https://doi.org/10.1002/rnc.1799}, Lyapunov-based sampled-data control \cite{LIU2012102}, and fault-tolerant stochastic sampled-data fuzzy control \cite{10018178}, etc. But it requires systems to gather information, communicate, and update the control at predefined intervals, resulting in a huge waste of unnecessary resources. In contrast, event-triggered control significantly lowers the amount of data that is collected, transmitted, and processed by only starting pertinent operations when specific requirements are met \cite{7563411,6068223,6257423}. So far, valuable achievements have been published these years, such as Selivanov et al. explored a distributed event-triggered controller of diffusion semilinear PDEs \cite{SELIVANOV2016344}, Chen et al. investigated backstepping control of a class of linear parabolic PDEs \cite{doi:10.1137/23M1556836}, and so forth.
 
Nevertheless, the Zeno behavior can potentially arise in using a static continuous event-triggered controller \cite{4303247,10579772,9551964}, compromising its suitability for real-world system applications. This detrimental behavior can be circumvented by employing periodic event-triggered control \cite{6310015}. But it should be noted that not all available information is fully exploited in this process. To address these issues, Selivanov et al. proposed a switching approach \cite{7355303} and the switching event-triggered mechanism is explained as follows. For discrete time series $t_0<t_1<t_2<\cdots<t_k<\cdots$, the event-triggered instant is:
\begin{equation*}
\begin{aligned}
	t_{k+1}&= min\{
	t \geq  t_k + h  |			
	\boldsymbol{e}(t)^T \Omega\boldsymbol{e}(t)\geq \varepsilon \boldsymbol{y}^T(t)\Omega\boldsymbol{y}(t)\nonumber
	\},
\end{aligned}
\end{equation*}
where $h$ represents waiting time, and it is selected as the maximum sampling value preserving the stability to avoid Zeno behavior. In contrast to the static event-triggered mechanism \cite{4303247}, switching event-triggered method can further reduce computational and communication resources. And it also permits longer sampling periods for equivalent values of the event-trigger parameter, compared to periodic event-trigger mechanisms \cite{6310015}. This method can essentially lower the number of measurements sent over the communication network in addition to avoiding the Zeno behavior \cite{7355303,ZHANG2023119599}.
Furthermore, in real life, nonlinear phenomena are not uncommon. There have been many notable efforts in the research of nonlinear PDE systems, such as uncertain nonlinearities was considered in \cite{9269405,7885587}. Nonetheless, there is still a striking lack of research on PDE systems with unknown nonlinearities. To tackle this challenge, an advantageous approach that this article introduces to study unknown nonlinearities is the multilayer neural networks (MNNs) \cite{4469946}. The learning capabilities of neural networks (NNs) have advanced significantly during decades of research \cite{8183423,809035,8706985,8939337}, providing unmatched benefits in fitting nonlinear functions. Additionally, external disturbances also present significant challenges to this work. 

However, there are relatively limited study findings on switching event-triggered control issues for PDE systems, which greatly motivate our investigation. In this paper, we will focus on switching event-triggered control for  a class of parabolic PDE systems with unknown nonlinearities. The primary contributions and innovations are as follows:
\begin{enumerate}[label=\textbullet]
	\item This is the first attempt to employ switching event-triggered strategy for the control issues of the PDE system. Several earlier works have been focused on sampled-data control or static event-triggered control while examining parabolic PDE systems utilizing Galerkin approaches in \cite{7885587,9551964}, comparatively speaking, the switching event-triggered control used in our work performs significantly better in terms of resource conservation and control efficiency.  
	\item Unknown nonlinearity is handled by MNN technique in our work. Takagi-Sugeno fuzzy technology has been utilized to copy with uncertain nonlinearity in \cite{4469891}, while it encounters difficulties when faced with unknown nonlinearity. Differently, in light of the model analyzer \cite{1981control} and Back-Propagation (BP) NNs with Levenberg-Marquardt (LM) algorithm, we set MNN to reconstruct unknown nonlinearities with a relative accurate manner.
	\item Semi-globally uniformly ultimately bounded (SGUUB) and $H_\infty$ performance stable of the closed-loop system are discussed via Lyapunov functional method. The controllers are given in terms of bilinear matrix inequalities (BMIs), and novel algorithms are proposed to solve BMIs and optimize the attenuation level. Particularly, the exponentially stable is also derived at the absence of external disturbance.
\end{enumerate}

The residual section is structured as outlined hereinafter. Section \ref{s2} displays lemmas that will be required for the subsequent analysis. Section \ref{s3} provides the problem description and control objectives. Section \ref{s4} primarily overcomes the design challenges of the switching event-triggered controller and introduces stability analysis methods. And the effectiveness of these advanced designs is tested in Section \ref{s5}. Ultimately, a brief conclusion is ultimately reached in Section \ref{s6}.

$\textit{Notations:}$ $\mathbb{R}$ indicates the real number set, $\mathbb{R}^{n}$, $\mathbb{R}^{n\times m}$ represent, respectively, the $n$ dimensional Euclidean space and the collection of all real $n\times m$ matrices, $\left\| \cdot \right\|$ represents a matrix's determinant, a vector's Euclidean norm or a matrix's spectral norm, $\mathscr{L}^{2}(\Omega;\mathbb{R}^{n})$ represents a collection $\{\xi:\Omega\to\mathbb{R} ^{n}$ and $\left\| \xi \right\|_{2,\Omega }< \infty \}$, and $\left\| \xi \right\|_{2,\Omega }=\sqrt{\int_{\Omega}\xi^{T}(w)\xi(w) dw}$, $\langle  \sigma(w),\varphi(w)\rangle=\int_{\Omega} \sigma(w)\varphi(w)dw$. $P^{T}$ indicates the conversion of matrix $P$, and
\begin{equation*}
	\begin{bmatrix}
		M+M^{T}&N\\N^{T}&Q
	\end{bmatrix}=
	\begin{bmatrix}
		M+\ast&N\\\ast&Q
	\end{bmatrix}.
\end{equation*}
\section{Preliminaries}\label{s2}
\begin{lemma}[see \cite{809035,8620233}]\label{l1} For any bounded initial condition, if a continuously differentiable and positive definite function $V(\psi,t)$ is defined, s.t. the inequality is met:
	\begin{equation*}
		\gamma_{1}(\mid\psi\mid)\leq V(\psi,t)\leq\gamma_{2}(\mid\psi\mid),\ 	\dot{V}\leq-CV+D,
	\end{equation*}
	where $C,D>0$ are positive constants, then it can be obtained
	\begin{equation*}
		0\leq V(t)\leq V(0)\mathrm{e}^{-Ct}+\frac{ D}{C}(1-e^{-Ct}),
	\end{equation*}
	then the solution $\psi(t)$ is SGUUB.
\end{lemma}
\begin{lemma}[see \cite{SEURET20132860}]\label{l2} (Jensen Inequality) For any symmetric positive definite matrix $X\in\mathbb{R}^{n\times n}$, scalars $\gamma_{1}$ and $\gamma_{2}$ satisfy $\gamma_ {1}<\gamma_{2}$, vector function $\omega: [\gamma_{1},\gamma_{2}]\longrightarrow\mathbb{R}^{n}$, such that the following integral can be defined, then the inequality is satisfied:
	\begin{equation*}\begin{aligned}
			\int_{\gamma_{1}}^{\gamma_{2}}\omega^{T}(s)X\omega(s)ds&\geq\frac{1}{\gamma_{2}-\gamma_{1}}\bigg[\int_{\gamma_{1}}^{\gamma_{2}}\omega(s)ds\bigg]^{T}\\
			&\ \times X\bigg[\int_{\gamma_{1}}^{\gamma_{2}}\omega(s)ds\bigg].
	\end{aligned}\end{equation*}
\end{lemma}
\section{Problem Formulation}\label{s3}
\subsection{Description of PDE System}
Consider a class of parabolic PDE systems with unknown nonlinearities:
\begin{equation}
	\left\{\begin{aligned}
		\frac{\partial\bar{\xi}(p,t)}{\partial t}&=z_1(p)\frac{\partial\bar{\xi}(p,t)}{\partial p}+\frac{\partial(z_2(p)\frac{\partial\bar{\xi}(p,t)}{\partial p})}{\partial p}+f(\bar{\xi}(p,t))\\
		&\ +b_{2}(p)u(t)+b_{1}(p)d(t),\\
		y(t)&=\int_{\Omega}\bar{c}(p)\bar{\xi}(p,t)dp,
	\end{aligned}\right.
\end{equation}
with the boundary conditions:
\begin{equation}
	h_{1,a}\bar{\xi}(\alpha_{a},t)+h_{2,a}\frac{\partial\bar{\xi}(\alpha_{a},t)}{\partial p}=0,\ a=1,2,
\end{equation}
and the initial condition:
\begin{equation}
	\bar{\xi}(p,0)=\bar{\xi}_{0}(p),
\end{equation}
where $\Omega=[\alpha_{1},\alpha_{2}]\subset \mathbb{R}$ is the process defining domain, $p\in\Omega$ represents spatial coordinates, $t\in[0,\infty]$ represents time, $\bar{\xi}(p,t)\in \mathbb{R}$ signifies state, $u(t)\in \mathbb{R}^{n_{u}}$ denotes the control input, $d(t)\in \mathbb{R}^{n_{d}}$ stand for disturbances, $y(t)\in \mathbb{R}^{n_{y}}$ stands for output. $z_i(p)>0,\ i=1,2$ stand for real continuous functions about $p$, and $f(\bar{\xi}(p,t))$ stands for an unknown nonlinear function satisfying $f(0)=0$ and locally Lipschitz continuous, $b_{2}(p)\in \mathscr{L}^{2}(\Omega;\mathbb{R}^{n_{u}})$, $b_{1}(p)\in \mathscr{L}^{2}(\Omega;\mathbb{R}^{n_{d}})$, $c(p)\in \mathscr{L}^{2}(\Omega;\mathbb{R}^{n_{y}})$. $h_{a}$, $a=1,2$ is constant. $\bar{\xi}_{0}(p)$ stands for a smooth enough function about $p$.

\begin{assumption}\label{a1}
	$d(t)$ is a constrained external disturbance generated with the input, satisfies:
	\begin{equation*}\mid d(t)\mid\leq D_{1}(D_{1}>0).\end{equation*}
\end{assumption}

For the PDE system (1)–(3), we consider the following $H_\infty$ control performance:
\begin{equation*} \int_{0}^{T}{y^{T}( \varrho)}y( \varrho)d \varrho\leq\gamma^{2}\int_{0}^{T}{d^{T}( \varrho)}d( \varrho)d \varrho,
\end{equation*}
where $T $ is the final time of control.

The design objectives of this paper are as follows:
\begin{enumerate}[label=\textbullet]
	\item In the presence of disturbances, the controller is designed to attain SGUUB of PDE systems (1)-(3).
	\item In the absence of disturbances, the controller is designed to achieve exponential stability of PDE systems (1)-(3).
	\item To meet the performance objectives, the controller is designed to achieve $H_\infty$ performance stable of PDE systems (1)-(3).
\end{enumerate}
%
\subsection{Fast and Slow Separation of PDE System}
In this subsection, semigroup theory is used to discuss the existence of solutions for PDE systems (1)–(3). The Galerkin approach is then used to develop the finite dimensional approximation of the PDE system. Assuming $\mathscr{H}=\mathscr{L}^{2}(\Omega;\mathbb{R})$ is the Hilbert space defined on $\Omega$, defining the spatial differential operator $\mathscr{A}$ of PDE systems (1)-(3):
\begin{equation*}
	\mathscr{A}\bar{\xi}=z_1(p)\frac{\partial\bar{\xi}(p,t)}{\partial p}+\frac{\partial (z_2(p)\frac{\partial\bar{\xi}(p,t)}{\partial p})}{\partial p},
\end{equation*}
and the domain is 
$\textit{D}(\mathscr{A})=\{\bar{\xi}\in \mathscr{H}\mid\bar{\xi}$ and $\frac{\partial \bar{\xi}}{\partial p}$ is absolutely continuous,\ $\frac{\partial\bar{\xi}(p,t)}{\partial p},\frac{\partial^2 \bar{\xi}}{\partial p^2}\in \mathscr{H}$, boundary conditions $h_{1,a}\bar{\xi}(\alpha_{a},t)+h_{2,a}\frac{\partial\bar{\xi}(\alpha_{a},t)}{\partial p}=0,\ a=1,2\}$.

After that, the system can be rewritten:
\begin{equation}\label{e4}
	\dot{\bar{\xi}}(t)=\mathscr{A}\bar{\xi}(t)+\bar{f}(\bar{\xi}(t))+\mathscr{B}_{2}u(t)+\mathscr{B}_{1}d(t)
	,\end{equation}
where $\bar{f}(\bar{\xi}(t))=f(\bar{\xi}(\cdot,t))$, $\mathscr{B}_{2}=b_{2}$, $\mathscr{B}_{1}=b_{1}$, $\bar{\xi}(0)=\bar{\xi}_{0}(\cdot)$.

Similar to the process of 2.1 in the book \cite{curtain2012introduction}, and referring to the Theorem 1.4 in the literature \cite{pazy2012semigroups}, there exists a unique solution to equation (\ref{e4}) by means of the separation of variables:
\begin{equation*}
	\begin{aligned}
		\bar{\xi}(t)&=F(t)\bar{\xi}_{0}+\int_{0}^{t}F(t-s)(f(\bar{\xi}(s))\\
		&\ +\mathscr{B}_{2}u(s)+\mathscr{B}_{1}d(s))ds,\ t\in[0,T],
	\end{aligned}
\end{equation*}
where $u(t)\in \mathscr{L}^2([0,T];\mathbb{R}^{n_{u}})$, $d(t)\in \mathscr{L}^2([0,T];\mathbb{R}^{n_{d}})$, initial conditions $\bar{\xi}(0)\in \textit{D}(\mathscr{A})$.

Define the eigenvalue problem for the operator $\mathscr{A}$: $\mathscr{A}\phi_{j}(p)=\lambda_{j}\phi_{j}(p)$, $j=1,2$,$\ldots$,$\infty$, where $\lambda_ {j}$ and $\phi_{j}(p)$ denote the $j$th eigenvalue and the corresponding standard orthogonal eigenfunction, respectively, and the operator $\mathscr{A}$ is a self-concomitant operator.

\begin{assumption}\label{a3}
	All the eigenvalues of $\mathscr{A}$ are arranged, i.e. $\lambda_{k}\geq\lambda_{k+1}$, and there exists a finite positive integer $m$, s.t. $\lambda_{m+1}<0$ and $\frac{|\lambda_m|}{|\lambda_{m+1}|}=O(\varepsilon)$, where $\varepsilon=\frac{|\lambda_L|}{|\lambda_{m+1}|}<1$ is a small positive
	number and $\lambda_L$ is the largest nonzero eigenvalue.
\end{assumption}

For each PDE state vector piecewise smooth function can be presented as follow:
\begin{equation}\label{e5}
	\bar{\xi}(p,t)=\sum_{i=1}^{\infty}\xi_{i}(t)\phi_{i}(p),
\end{equation}
where $\xi_{i}(t)$ is a time-varying coefficient.

\begin{remark}\label{r2}
The physical importance of each modal component resulting from modal decomposition is unique. In the context of vibration systems, for example, specific modes individually correlate to various vibration frequencies and mode shapes. In the context of structural dynamics, modal decomposition offers a clear illustration of the characteristics linked to each vibration mode when it comes to the vibration analysis of structural elements like beams and plates.
\end{remark}

Based on Assumption 2, applying Galerkin method, the infinite-dimensional ODE system that follows can be produced by simultaneously differentiating both sides of equation (\ref{e5}) with regard to $t$:
\begin{equation}\label{e6}
	\left\{\begin{aligned}
		\dot{\xi}_{s}(t)&=A_{s}\xi_{s}(t)+f_{s}(\xi_s(t),\xi_f(t))+B_{2,s}u(t)+B_{1,s}d(t),\\
		\dot{\xi}_{f}(t)&=A_{f}\xi_{f}(t)+f_{f}(\xi_s(t),\xi_f(t))+B_{2,f}u(t)+B_{1,f}d(t),\\
		y(t)&=C_{s}\xi_{s}(t)+C_{f}\xi_{f}(t),\\
		\xi_{s}(0)&=\xi_{s,0},\ \xi_{f}(0)=\xi_{f,0},
	\end{aligned}\right.
\end{equation}
where
\begin{equation*}\begin{aligned}
		\xi_{s}(t)&=[\xi_{1}(t)\ldots \xi_{m}(t)]^{T}\in \mathbb{R}^{m},\\
		\xi_{f}(t)&=[\xi_{m+1}(t)\ldots \xi_{\infty}(t)]^{T}\in l^2,\\
		\xi(t)&=[{\xi_{s}}^T(t),\ {\xi_{f}}^T(t)]^T\in l^2,\\
		A_{s}&=\{\lambda_{1}\ldots\lambda_{m}\},\ A_{f}=\{\lambda_{m+1}\ldots\lambda_{\infty}\},\\
		\varphi_{s}(p)&=[\phi_{1}(p)\ldots\phi_{m}(p)]^T,\\
		\varphi_{f}(p)&=[\phi_{m+1}(p)\ldots\phi_{\infty}(p)]^T,\\
		B_{h,\beta}&=\langle\varphi_{\beta}(p),b_{h}(p)\rangle,\ h=1,2,\\
		C_{\beta}&=\int_{\Omega}\bar{c}(p){\varphi_{\beta}}^T(p)dp,\ \xi_{\beta,0}=\langle\varphi_{\beta}(p),\ \bar{\xi}_{0}(p)\rangle,\\
		f_{\beta}(\xi_s(t),\xi_f(t))&=\langle\varphi_{\beta}(p),f(\bar{\xi}(p,t))\rangle,\ \beta=s,f. 
	\end{aligned}
\end{equation*}

Then the $H_\infty$ control performance problem can be similarly rewritten as
\begin{equation*} \int_{0}^{T}{y_s^{T}(t)}y_s(t)dt\leq\gamma_s^{2}\int_{0}^{T}{d^{T}(t)}d(t)dt,
\end{equation*}
where $y_{s}(t)=C_{s}\xi_{s}(t)$.

The singular perturbation model that follows is obtained by increasing the $\xi_{f}$-subsystem in equation (\ref{e6}) by $\varepsilon$:
\begin{equation}\label{e7}
	\left\{\begin{aligned}
		\dot{\xi}_{s}(t)&=A_{s}\xi_{s}(t)+f_{s}(\xi_s(t),\xi_f(t))+B_{2,s}u(t)+B_{1,s}d(t),\\
		\varepsilon\dot{\xi}_{f}(t)&=A_{f\varepsilon}\xi_{f}(t)+\varepsilon f_{f}(\xi_s(t),\xi_f(t))\\
		&\ +\varepsilon B_{2,f}u(t)+\varepsilon B_{1,f}d(t),\\
		y(t)&=C_{s}\xi_{s}(t)+C_{f}\xi_{f}(t),
	\end{aligned}\right.
\end{equation}
where $A_{f\varepsilon}=\varepsilon A_{f}$.

In accordance with the singular perturbation theory \cite{2002Nonlinear}, based on $\lambda_ {m+1}<0$, the fast subsystem has been found to be globally exponentially stable. Nonlinear finite-dimensional slow system can be gotten:
\begin{equation}\label{e8}
	\left\{\begin{aligned}
		\dot{\xi}_{s}(t)&=A_{s}\xi_{s}(t)+f_{s}(\xi_{s}(t),0)+B_{2,s}u(t)+B_{1,s}d(t),\\
		y(t)&=C_{s}\xi_{s}(t)+\varphi(t),
	\end{aligned}\right.
\end{equation}
where $\varphi(t)=y(t)-y_{s}(t)$.

\begin{remark}
	To ensure clarity, we initially address the challenge of designing a finite-dimensional switching event-triggered controller for slow system (\ref{e8}), utilizing the output $y_s(t)$ of the slow system rather than the actual measurement output $y(t)$.
\end{remark}
\subsection{NNs Approximation of Finite-Dimensional ODE System}
\begin{figure}[!t]
	\centering
	\includegraphics[width=0.4\textwidth]{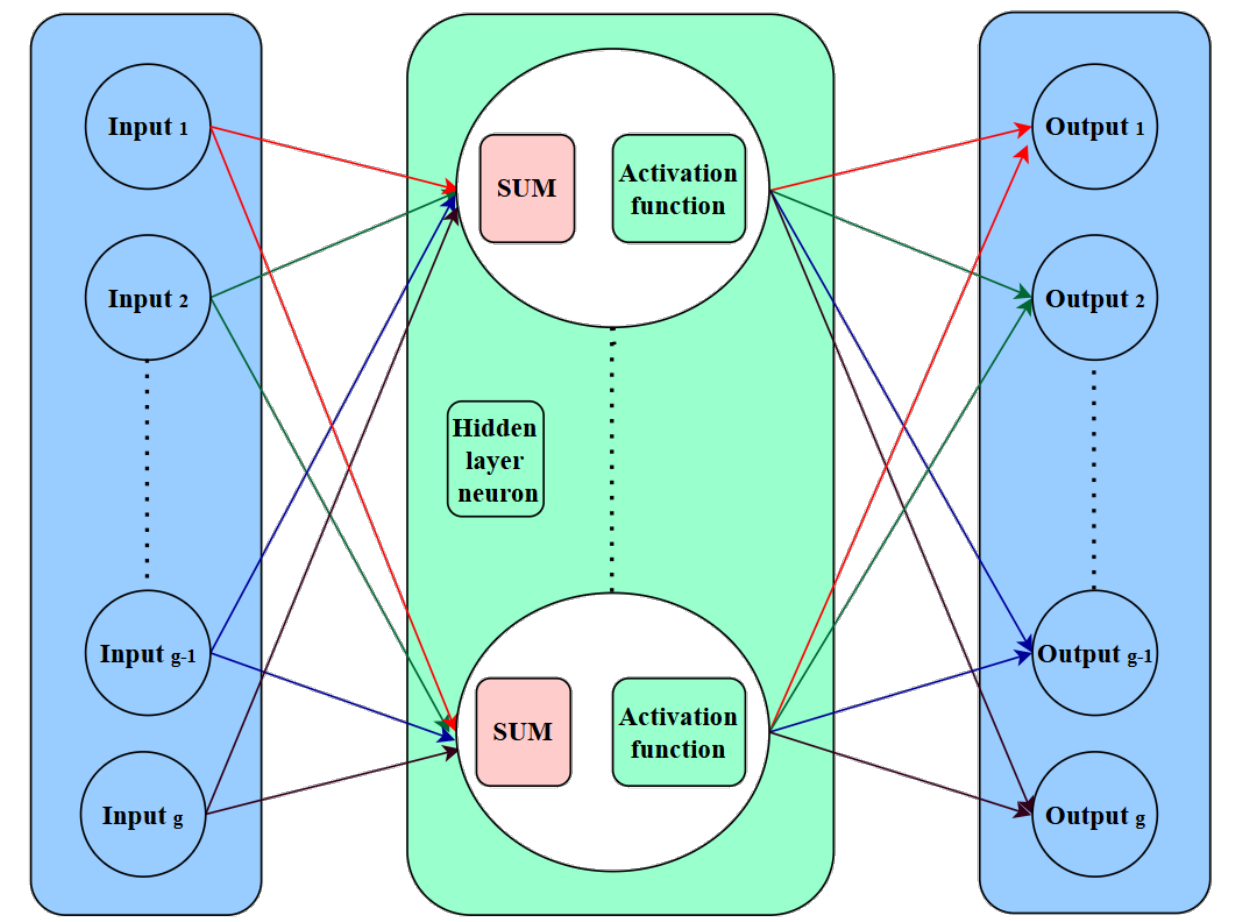}
	\caption{Three-Layer Neural Network Architecture.}
	\label{Fig. 1}
\end{figure}

The Galerkin approach implies that the slow subsystem accurately captures the primary dynamic features of the original system. But due to the unknown nonlinear function $f_{s}(\xi_{s}(t),0)$, the finite-dimensional slow system (\ref{e8}) is still inapplicable. As shown in Fig. \ref{Fig. 1}, we set up a three-layer NN to approximate complex nonlinear functions.

The matrix-vector notation describes a MNN featuring a single hidden layer and the absence of bias terms, expressed elegantly as follow:
\begin{equation}\label{e9}
	f_{nn}(\xi_s(t),\boldsymbol{V},\boldsymbol{W})=\boldsymbol{W}\boldsymbol{\mu}(\boldsymbol{V}\xi_s(t)),
	\end{equation}
where $f_{nn}(\xi_s(t),\boldsymbol{V},\boldsymbol{W})\in\mathbb{R}^{m}$ stands for the output of the network,
$\boldsymbol{V}\in\mathbb{R}^{n_h\times m}$ represents the first-to-second layer interconnection weight
matrix, $\boldsymbol{W}\in\mathbb{R}^{m\times n_{h}}$ represents the second-to-third layer interconnection weight matrix, and $n_h$ is the number of hidden neurons.
$\boldsymbol{\mu}(\cdot):\mathbb{R}^{n_h} \mapsto \mathbb{R}^{n_h}$ denotes the activation function vector of
the network, which is defined as
\begin{equation*}\boldsymbol{\mu}(\boldsymbol{\varsigma})=[\mu_1(\varsigma_1),\mu_2(\varsigma_2),\ldots,\mu_{n_h}(\varsigma_{n_h})],\end{equation*}
where $\boldsymbol{\varsigma} = [ \varsigma_{1},\varsigma_{2},\ldots, \varsigma_{n_{h}}]^{T} \in \mathbb{R}^{n_{h}},\mu_{i}( \cdot ) (i\in \mathbb{H}\overset{\underset{\mathrm{def}}{}}{=} \{1, 2,\ldots,n_{h}\})$ is the activation function of the $i$th hidden neuron. $\mu_{i}(\cdot)$ has the following assumption in this paper.
\begin{assumption}
	$\mu_{i}(\cdot)$ is differentiable. $\mu_{i}(\varsigma_{i})\in[-q_{i},q_{i}],q_{i}>0$, $\forall\varsigma_{i}\in\mathbb{R}$ and $\mu_{i}(0)=0$.
\end{assumption}

A common activation function that adheres to Assumption 3 is the bipolar sigmoid function, which exhibits symmetry with respect to the origin
\begin{equation}\label{e10}
	\mu_i(\varsigma_i)=q_i\left(\frac{2}{1+\exp(-\varsigma_i/r_i)}-1\right),
	\end{equation}
where $q_i$ and $r_i$ are parameters of the function.

\begin{remark}
	If the activation function is discontinuous, the gradient at the points of discontinuity becomes undefined, which impedes the proper functioning of the BP algorithm. Consequently, it becomes challenging to effectively train NNs to fit nonlinear functions. Therefore, when selecting and utilizing activation functions for hidden layers, continuity is a crucial factor to consider. Moreover, the continuity of the activation function facilitates stable learning and convergence of NNs in the process of fitting nonlinear functions.
\end{remark}

Subsequently, we can establish the learning rules for the MNN based on (\ref{e9}) and (\ref{e10}). The method employed in this paper is LM algorithm. Compared with the BP algorithm, this algorithm has a better global search ability in fitting nonlinear processes. 

It is necessary to get the state vector $\bar{\xi}(p,t)$ of the PDE systems (1)-(3).  But it is impossible that the state $\bar{\xi}(p,t)$ is completely available. Thus we apply the following method to measure $\bar{\xi}(p_i,t),i=1\ldots m,$ at spatial positions. From (\ref{e5}), it can be devised:
\begin{equation*}
	\bar{\xi}(p_i,t)\approx\varphi_s^T(p_i)\xi_{s}(t).
\end{equation*}

Define $\bar{\xi}(p,t)=[\bar{\xi}(p_1,t),\ldots,	\bar{\xi}(p_m,t)]^T$, $\bar{\varphi}_s(p)=[\varphi_s(p_1),\ldots, \varphi_s(p_m)]^T$, then one holds
\begin{equation*}
	\bar{\xi}(p,t)\approx\bar{\varphi}_s(p)\xi_{s}(t),
\end{equation*}
then, if $\bar{\varphi}_s(p)$ is not singular, it can be obtained:
\begin{equation*}
	\xi_{s}(t)\approx  \bar{\varphi}_s^{-1}(p)\bar{\xi}(p,t),
\end{equation*}
\begin{remark}
	It is assumed that the sampling locations are well chosen in the above analysis, so that the $\varphi_s(p)$ is not singular. When measurements from more than $m$ spatial locations are attainable, a least squares approach can be employed to estimate $\xi_s(t)$ with precision
	\begin{equation*}
		\xi_{s}(t)\approx  (\bar{\varphi}_s^T\bar{\varphi}_s)^{-1}(p)\bar{\varphi}_s^T(p)\bar{\xi}(p,t).
	\end{equation*}
	For more detail regarding this issue, one can
	refer to \cite{1981control}.
\end{remark} 

According to the definition of the derivative, when a sufficiently small measurement interval $\Delta T_s$ exists, the following equation will be satisfied:
\begin{equation*}
	\dot{\xi}_s(t)=\frac{\xi_s(t+\Delta T_s)-\xi_s(t)}{\Delta T_s},
\end{equation*}
so the MNN target outputs can be deduced under $u(t)=0$:
\begin{equation}\label{e11}
	\begin{aligned}
		f_s(\xi_s(t),0)& =\dot{\xi}_{s}(t)-A_{s}\xi_{s}(t) \\
		&\approx\frac{\xi_{s}(t+\Delta T_{s})-\xi_{s}(t)}{\Delta T_{s}}-A_{s}\xi_{s}(t).
	\end{aligned}
\end{equation}
\begin{remark}
In reality, numerous scenarios emerge wherein the disturbance $d(t)$ coincides with $u(t)$. For instance, consider the catalytic rod reaction; if the feeding of reactants is entirely halted, thereby terminating all inputs, the chemical reaction within the reactor will inevitably cease due to the depletion of reactants. Consequently, the disturbance arising from feed factors will gradually dissipate, ultimately converging to a disturbance value of zero. 
\end{remark}

Based on the approximation theorem \cite{Stable2001} for
MNNs, there exists a positive scalar $\delta>0$, such that the optimized weight matrices $\boldsymbol{W}^*$ and $\boldsymbol{V}^*$ can be found as:
\begin{equation*}
	(\boldsymbol{W}^*,\boldsymbol{V}^*)=\arg\min_{(\boldsymbol{W},\boldsymbol{V})}\left\{\max_{\xi_s\in\mathcal{D}}\|f_s(\xi_s,0)
	-f_{nn}(\xi_s,\boldsymbol{V},\boldsymbol{W})\|\right\}\end{equation*}
where $\mathcal{D}$ is a compact subset of $\mathbb{R}^{m}$, such that
\begin{equation*}\max_{\xi_s\in\mathcal{D}}\|f_s(\xi_s,0)-f_{nn}(\xi_s,\boldsymbol{V}^*,\boldsymbol{W}^*)\|\leq\delta\|\xi_s\|\\, \end{equation*}
The error expression can be derived:
\begin{equation}\label{e12}
	\delta=\sup_{\xi_s(t)\neq0}\frac{\|f_s(\xi_s(t),0)-f_{nn}(\xi_s(t),\boldsymbol{V}^*,\boldsymbol{W}^*)\|}{\|\xi_s(t)\|}.\end{equation}
	
For clarity, the approximation process of nonlinear functions of system (\ref{e8}) can be concisely summarized as follows:
\begin{enumerate}[label=\textbullet]
	\item Step 1: Set up the training set. A polyhedral region $\mathcal{D}=\{\xi_{s}\in\mathbb{R}^{m},-\rho_{j1}\leq \xi_{j}(t)\leq\rho_{j2},\ j=1,2,\dots,m\}$, where $\rho_{j1}$ and $\rho_{j2}$ are non-negative. As the MNN inputs, choose a set of initial state conditions in $\mathcal{D}$. Next, apply each state with $u(t)=0,\ d(t)=0$, measure the next state after $\Delta T_{s}$ seconds. The NN target outputs can therefore be obtained as (\ref{e11}). 

	\item Step 2: Network initialization. Initialize the weight matrices $\boldsymbol{V}$ and $\boldsymbol{W}$ of NN randomly. Build a hidden layer based on (\ref{e9}) in which all activation functions are selected as (\ref{e10}) and the number of hidden neurons $n_h$ is predetermined. 

	\item Step 3: Network training by Algorithm \ref{alg1}.

	\item Step 4: Estimate the bound $\delta$. According to (\ref{e12}), we can estimate $\delta$.
	\end{enumerate}

Define $\mu_{i}^{\prime}(\varsigma_{i})\stackrel{\Delta}{=}(\mathrm{d}\mu_{i}(\varsigma_{i}))/(\mathrm{d}\varsigma_{i})$, and respectively, $g_{\mathrm{i,min}}\overset{\underset{\mathrm{def}}{}}{=}\operatorname*{min}_{\varsigma_{i}}\mu_{i}^{\prime}(\varsigma_{i}),\ g_{i,\mathrm{max}}\overset{\underset{\mathrm{def}}{}}{=}\operatorname*{max}_{\varsigma_{i}}\mu_{i}^{\prime}(\varsigma_{i}).$

Then one has:
\begin{equation}\label{e13}
	\mu_i(\varsigma_i)/\varsigma_i\in[g_{i,\min},g_{i,\max}],\quad i\in\mathbb{H},
\end{equation}
and from (11), we can obtain $g_{i,\mathrm{min}} = 0$ and $g_{i,\mathrm{max}} = \frac{q_i}{2r_i}$, $i\in\mathbb{H}$.

\begin{algorithm}[H]
	\caption{Network training with LM algorithm.}\label{alg:alg1}
	\begin{algorithmic}
		\STATE 
		\STATE {\textsc{Step 3.1}}
		\STATE \hspace{0.5cm}$ \textbf{LM algorithm initialization. } $
		\STATE \hspace{0.5cm} Initialize iteration step size $=\mu_0$. Let accuracy of convergence $= \epsilon_c$. Input $X_i=[x_{i,1},\ldots,x_{i,N_1}],i=m$. Set $k=0$, $k_{\max}=N$.
		\STATE {\textsc{Step 3.2}}
		\STATE \hspace{0.5cm}$ \textbf{Calculate NN output and error.} $
		\STATE \hspace{0.5cm} Output $f_{i,n,n}=\boldsymbol{W_i}\boldsymbol{\mu}({\boldsymbol{V_i}}X_{i})$. 
		
		\ \ \ \ \ Error $E_i=\frac{1}{2N_{1}}\sum_{j=1}^{N_{1}}\left\|f_{i,n,n,j}-f_{i,j}(\xi_s,0)\right\|^2$.
		\STATE {\textsc{Step 3.3}}
		\STATE \hspace{0.5cm}$ \textbf{Gradient computation.} $
		\STATE \hspace{0.5cm}  
		$\begin{aligned}
			\frac{\partial E_i}{\partial \boldsymbol{W_i}}&=\frac{1}{N_1}\sum_{j=1}^{N_{1}}(f_{i,n,n,j}-f_{s,j}(\xi_s,0))\odot \boldsymbol{\mu}, \\
			\frac{\partial E_i}{\partial \boldsymbol{V_i}}&=\frac{1}{N_1}\sum_{j=1}^{N_{1}}(f_{i,n,n,j}-f_{s,j}(\xi_s,0))\boldsymbol{W_i}^T\odot \boldsymbol{\mu}'X_{i}^T. 		\end{aligned}$
		Construct Jacobian matrix $J_{ij}=\frac{\partial E_i}{\partial {x_{i,j}}}$.
		\STATE {\textsc{Step 3.4}}
		\STATE \hspace{0.5cm}$ \textbf{Parameter update.} $
		\STATE \hspace{0.5cm}  Let $h=(J^T J+\tilde{\epsilon}I)+J^T e$, where $\tilde{\epsilon}$ is tuning parameter, $e=f_{i,n,n,j}-f_{i,j}(\xi_s,0)$. $\boldsymbol{W_i}=\boldsymbol{W_i}-h_{\boldsymbol{W_i}}$, $\boldsymbol{V_i}=\boldsymbol{V_i}-h_{\boldsymbol{V_i}}$. 
		Calculate $E_{i,new}$ .
		\STATE {\textsc{Step 3.5}}
		\STATE \hspace{0.5cm}$ \textbf{Convergence judgment.} $
		\STATE \hspace{0.5cm} 
		If $|E_{i}-E_{i,new}|<\epsilon_c$, proceed to the next step;
		
		\ \ \ \ Else if $|E_{i}-E_{i,new}|\geq\epsilon_c$, update $\tilde{\epsilon}$, set $k=k+1$, $E_{i}=E_{i,new}$, return to Step 3.4;
		
		\ \ \ \ \ Else if $k=N$, proceed to the next step.
		\STATE {\textsc{Step 3.6}}
		\STATE \hspace{0.5cm}$ \textbf{Output $\boldsymbol{W}^*$ and $\boldsymbol{V}^*$.} $
	\end{algorithmic}
	\label{alg1}
\end{algorithm}
Let $\boldsymbol{\varsigma}(t)=\boldsymbol{V}^*\xi_s(t)$. During the aforementioned process, the MNN outputs can be gotten $f_{nn}(\xi_s,\boldsymbol{V}^*,\boldsymbol{W}^*)=\boldsymbol{W}^*\boldsymbol{\mu}(\boldsymbol{\varsigma}(t))$. Thus the system (\ref{e8}) can be rewritten as:
\begin{equation}\label{e14}
	\left\{\begin{aligned}
		\dot{\xi}_{s}(t)&=A_{s}\xi_{s}(t)+	\boldsymbol{W}^{*}\boldsymbol{\mu}(\boldsymbol{\varsigma}(t))+\Delta\boldsymbol{A}(\xi_{s}(t))\\
		&\ +B_{2,s}u(t)+B_{1,s}d(t),\\
		y_s(t)&=C_{s}\xi_{s}(t),
	\end{aligned}\right.
\end{equation}
where $\Delta \boldsymbol{A}(\xi_s(t))=f_s(\xi_s(t),0)-\boldsymbol{W}^*\boldsymbol{\mu}(\boldsymbol{\varsigma}(t))$.

\section{Switching Event-Triggered Control Design}\label{s4}
\begin{figure}[!t]
	\centering
	\includegraphics[width=0.4\textwidth]{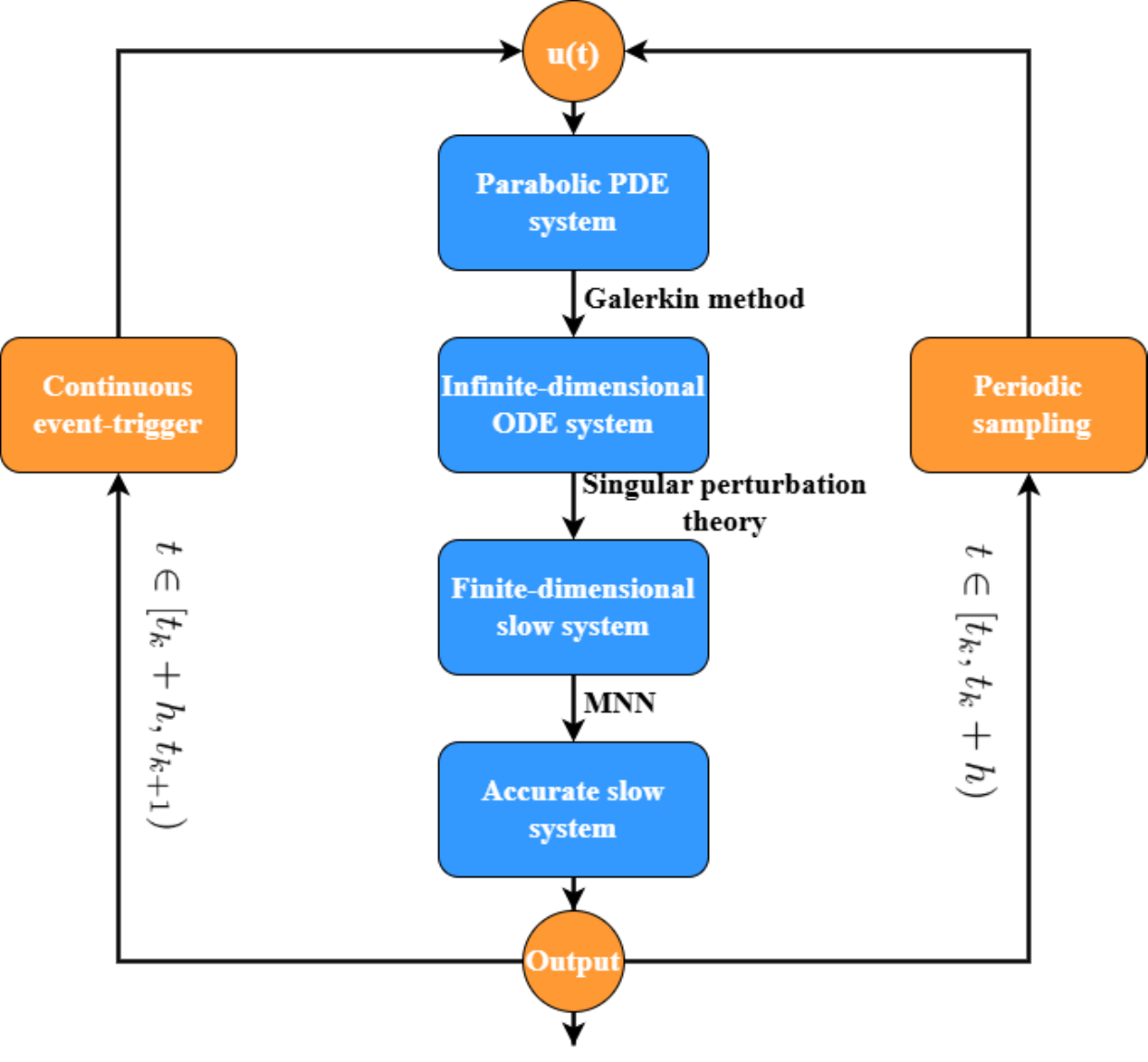}
	\caption{Switching event-triggered controller design  for	parabolic PDE systems.}
	\label{Fig.2}
\end{figure}

To accomplish stability of the relative accurate system (\ref{e14}), the switching event-triggered output feedback controllers and $H_\infty$ control design are proposed in this section. The controller design is given as Fig. \ref{Fig.2}.
\\
\subsection{Design of Switching Event-Triggered Controller}
To reduce communication resources, a switching event-triggered controller is designed in this subsection. On one hand, this approach fundamentally avoids the Zeno behavior. On the other hand, the proposed switching control strategy guarantees stability of the closed-loop system.

The following event-triggered condition is presented as:
\begin{equation}\label{e15}
	\begin{aligned}t_{k+1}&=\min\{t\geq t_{k}+h\mid(y_s(t)-y_s(t_{k}))^{T}\Lambda (y_s(t)-y_s(t_{k}))\\&\geq\epsilon y_s^{T}(t)\Lambda  y_s(t)\},\end{aligned}\end{equation}
where the matrix $\Lambda>0$, scalars $\epsilon\geq 0$, $h> 0$. 

The switching event-triggered output feedback controller is designed as:
%
\begin{equation}\label{e16}
	u(t)=\left\{	\begin{aligned}
		&KC_s(\int_{t-\tau(t)}^{t}\dot{\xi}_s(s) ds-\xi_s(t)),\ t\in[t_k,t_k+h],\\
		&-KC_s(\frac{e(t)}{C_s}+\xi_s(t)),\ t\in[t_k+h,t_{k+1}],\\
	\end{aligned}\right.
\end{equation}
where \begin{equation*}\begin{aligned}
		\tau(t)&=t-t_{k}\leq h,\quad t\in[t_{k},t_{k}+h),\\e(t)&=y_s(t_{k})-y_s(t),\quad t\in[t_{k}+h,t_{k+1}).\end{aligned}\end{equation*}

By incorporating the controller (\ref{e16}) into the system (\ref{e14}), we can obtain the following switching closed-loop system:
\begin{equation}\label{e17}
	\begin{aligned}
		\dot{\xi}_{s}(t)&=(A_{s}-B_{2,s}KC_s)\xi_{s}(t)+	\boldsymbol{W}^{*}\boldsymbol{\mu}(\boldsymbol{\varsigma}(t))+\Delta\boldsymbol{A}(\xi_{s}(t))\\
		&\ +\chi(t) B_{2,s}KC_s\int_{t-\tau(t)}^{t}\dot{\xi}_s(s) ds\\
		&\ -(1-\chi(t)) B_{2,s}Ke(t)+B_{1,s}d(t),\\
	\end{aligned}
\end{equation}
where \begin{equation*}
	\chi(t)=\left\{\begin{aligned}
		&1,\ t\in[t_k,t_k+h],\\
		&0,\ t\in[t_k+h,t_{k+1}].\\
	\end{aligned}\right.
\end{equation*}

For closed-loop system (\ref{e17}), the Lyapunov functional is designed as follows:
\begin{equation}\label{e18}
	V(t)=\sum_{j=1}^{4}V_{j}(t),
\end{equation}
where
\begin{equation*}\begin{aligned}
		V_{1}(t)&={\varrho}_{1}^T(t) P{\varrho}_{1}(t),\\
		V_{2}(t)&=(h-\tau(t))\int_{t-\tau(t)}^{t}e^{2\alpha (s-t)}\dot{\xi}_{s}^{T}(s)U\dot{\xi}_{s}(s)ds,\\
		V_{3}(t)&=(h-\tau(t))[\xi_s^T(t)(\frac{Q_1}{2}+\ast)\xi_s(t)+\xi_s^T(t-\tau(t))\\
		&\ \times(\frac{-2Q_2+Q_1}{2}+\ast)\xi_s(t-\tau(t))
		+2\xi_s^T(t)\\
		&\ \times(-Q_1+Q_2)\xi_s(t-\tau(t))],\\
		V_{4}(t)&=\sum_{i=1}^{n_h}2\varpi_i\int_{0}^{\varsigma_i(t)}\mu_i(\rho ) (g_{i,max}-g_{i,min})d\rho ,\\
		\varrho_1(t)&=\begin{bmatrix}\xi_s(t)\\ \xi_s(t-\tau)
		\end{bmatrix},\
		P=\begin{bmatrix}
			P_{11}&P_{12}\\
			P_{21}&P_{22}\\
		\end{bmatrix},
\end{aligned}\end{equation*}
and positive-definite matrices $P\in\mathbb{R}^{2m\times 2m}>0$, $U\in\mathbb{R}^{m\times m}>0$, with suitable matrices $P_{11}$,  $P_{12}$,  $P_{21}$,  $P_{22}$, $Q_1$, $Q_2$$\in\mathbb{R}^{m\times m}$, scalars $\alpha$, $\varpi_i>0,\ i\in\mathbb{H}$.

When $\chi(t)=1$, $V_{\chi(t)=1}(t)=\sum_{j=1}^{4}V_{j}(t)$, we proceed to conduct the following analysis.

Define $\nu_1(t)=\frac{1}{\tau(t)}\int_{t-\tau(t)}^{t}\dot{\xi}_{s}(s)ds$, it is satisfied that
\begin{equation*}
	\tau(t)\nu_1(t)=\xi_{s}(t)-\xi_{s}(t-\tau(t)),
\end{equation*} 
and $\displaystyle \lim_{\tau(t) \to 0} \nu_1(t)=\dot{\xi}_s(t)$. From Lemma 2, one has:
\begin{equation}\label{e19}
	\int_{t-\tau(t)}^t\dot{\xi}^T(s)U\dot{\xi}(s)\mathrm{d}s\geq\tau(t)\nu_1^TU\nu_1,\end{equation}

Derive $V_{1}(t)$ about $t$ and substitute (\ref{e17}) into it to obtain:
\begin{equation}\label{e20}
	\begin{aligned}
		\dot{V}_{1}(t)&={\varrho}_{1}^T(t) P\dot{\varrho}_{1}(t)\\
		&=\xi_{s}^{T}(t)[P_{11}(A_{s}-B_{2,s}KC_s)+\ast]\xi_{s}(t)\\
		&\ +2\xi_{s}^{T}(t)P_{11}B_{2,s}KC_s\tau(t)\nu_1(t)+2\xi_{s}^{T}(t)P\boldsymbol{W}^{*}\boldsymbol{\mu}(\boldsymbol{\varsigma}(t))\\
		&\ +2\xi_{s}^{T}(t)P_{11}\Delta\boldsymbol{A}({\xi}_{s}(t))+2\xi_{s}^{T}(t)P_{11}B_{1,s}d(t)\\
		&\ +2\xi_s^T(t-\tau)P_{21}\dot{\xi}_s(t).
	\end{aligned}
\end{equation}

Deriving $V_{2}(t)$ derivative with regard to $t$, we can attain:
\begin{equation}
	\begin{aligned}\label{e21}
		\dot{V}_{2}(t)&=-\int_{t-\tau(t)}^t\mathrm{e}^{2\alpha(s-t)}\dot{\xi}^T(s)U\dot{\xi}(s)\mathrm{d}s+(h-\tau(t))\dot{\xi}^T(t)\\
		&\ \times U\dot{\xi}(t)-2\alpha(h-\tau(t))\int_{t-\tau(t)}^{t}e^{2\alpha (s-t)}\dot{\xi}_{s}^{T}(s)U\dot{\xi}_{s}(s)ds
		.	\end{aligned}\end{equation}

By finding the derivative of $V_{3}(t)$ in relation to $t$, we can derive:
\begin{equation}\label{e22}
	\begin{aligned}
		\dot{V}_{3}(t)&=-\begin{bmatrix}\xi_s(t)\\\xi_s(t-\tau(t))\end{bmatrix}^T\begin{bmatrix}\frac{Q_1}{2}+\ast&-Q_1+Q_2\\\\\ast&\frac{-2Q_2+Q_1}{2}+\ast\end{bmatrix}\\
		&\ \times \begin{bmatrix}\xi_s(t)\\\xi_s(t-\tau(t))\end{bmatrix} +(h-\tau(t))[\dot{\xi}_s^{T}(t)(Q_1+\ast)\xi_s(t) \\
		&\ +2\dot{\xi}_s^T(t)(-Q_1+Q_2)\xi_s(t-\tau(t))].
	\end{aligned}
\end{equation}

Calculating the $V_{4}(t)$ derivative with respect to $t$ to get:
\begin{equation}\label{e23}
	\begin{aligned}
		\dot{V}_{4}(t)&=2\boldsymbol{\mu}^T(\boldsymbol{\varsigma})\boldsymbol{\Omega}(G_\mathrm{max}-G_\mathrm{min})\dot{\boldsymbol{\varsigma}}\\
		&=2\boldsymbol{\mu}^T(\boldsymbol{\varsigma})\boldsymbol{\Omega}G\boldsymbol{V}^*\dot{\xi}_s(t).
\end{aligned}\end{equation}
where \begin{equation*}\begin{aligned}
		\boldsymbol{\Omega}& =\mathrm{diag}\{\varpi_{1},\varpi_{2},\ldots,\varpi_{n_{h}}\}, \\
		G_{\mathrm{min}}& =\mathrm{diag}\{g_{1,\mathrm{min}},g_{2,\mathrm{min}},\ldots,g_{n_{h},\mathrm{min}}\},\\
		G_{\mathrm{max}}& =\mathrm{diag}\{g_{1,\mathrm{max}},g_{2,\mathrm{max}},\ldots,g_{n_{h},\mathrm{max}}\},\\
		G&=G_{\mathrm{max}}-G_{\mathrm{min}}. 
\end{aligned}\end{equation*}
From (\ref{e12}), for any $\beta_1>0$, it can be obtained:
\begin{equation}\label{e24}
	\beta _1[\delta^2 \xi_s^T(t)\xi_s(t)-\Delta\boldsymbol{A}^T(\xi_s)\Delta\boldsymbol{A}(\xi_s)]\geq 0,
\end{equation}
From (\ref{e13}), for any $l_i>0,i\in\mathbb{H}$, we can get:
\begin{equation*}
	\sum_{i=1}^{n_h}l_i[\mu_i(\varsigma_i)-g_{i,min}\varsigma_i][g_{i,max}\varsigma_i-\mu_i(\varsigma_i)]\geq 0,\ i\in\mathbb{H},
\end{equation*}
i.e.
\begin{equation}\label{e25}
	[\boldsymbol{\mu}(\boldsymbol{\varsigma} )-G_{min}\boldsymbol{\varsigma} ]^TL[G_{max}\boldsymbol{\varsigma} ]-\boldsymbol{\mu}(\boldsymbol{\varsigma })]\geq 0.
\end{equation}
where $L=diag\{l_1,\ldots,l_{n_h}\}$.

According to Young's inequality, for $t\in[t_k,t_k+h]$, the inequalities that follows is presented:
\begin{equation}\label{e26}
	\begin{aligned}
		2\dot{\xi}_{s}^T(t)NB_{1,s}d(t)&\leq\frac{1}{\chi^{2}}\dot{\xi}_{s}^T(t)NB_{1,s}(NB_{1,s})^{T}\dot{\xi}_{s}(t)+{\chi}^{2}D_{1}^{2},\\
		2{\xi}_{s}^T(t)P_{11}B_{1,s}d(t)&\leq\frac{1}{\chi^{2}}{\xi}_{s}^T(t)P_{11}B_{1,s}(P_{11}B_{1,s})^{T}\xi_{s}(t)+\chi^{2}D_{1}^{2}.
\end{aligned}\end{equation}
When $\chi(t)=0$, $V_{\chi(t)=0}(t)=\sum_{j=1}^{4}V_{j}(t)$, we can deduce:
\begin{equation}\label{e27}
	\begin{aligned}
		\dot{V}_{1}(t)&=\xi_{s}^{T}(t)[P_{11}(A_{s}-B_{2,s}KC_s)+\ast]\xi_{s}(t)\\
		&\ -2\xi_{s}^{T}(t)P_{11}B_{2,s}Ke(t)+2\xi_{s}^{T}(t)P_{11}\boldsymbol{W}^{*}\boldsymbol{\mu}(\boldsymbol{\varsigma}(t))\\
		&\ +2\xi_{s}^{T}(t)P_{11}\Delta\boldsymbol{A}(\xi_{s}(t))+2\xi_{s}^{T}(t)P_{11}B_{1,s}d(t)\\
		&\ +2\xi_s^T(t-\tau)P_{21}\dot{\xi}_s(t).
	\end{aligned}
\end{equation}

Similar to the above analysis, for $t\in[t_k+h,t_{k+1}]$, the inequalities that follows can be gotten:
\begin{equation}\label{e28}
	\begin{aligned}
		2\dot{\xi}_{s}^T(t)NB_{1,s}d(t)&\leq\frac{1}{(1-\chi)^{2}}\dot{\xi}_{s}^T(t)NB_{1,s}(N B_{1,s})^{T}\dot{\xi}_{s}(t)\\
		&\ +(1-\chi)^{2}D_{1}^{2},\\
		2{\xi}_{s}^T(t)P_{11}B_{1,s}d(t)&\leq\frac{1}{(1-\chi)^{2}}{\xi}_{s}^T(t)P_{11}B_{1,s}(P_{11}B_{1,s})^{T}\xi_{s}(t)\\
		&\ +(1-\chi)^{2}D_{1}^{2}.
\end{aligned}\end{equation}

\begin{theorem}\label{thm1} Given scalars $\alpha$, $\epsilon$, $\beta_2$, matrix $\Lambda$, and the weights $\boldsymbol{V}^{*},\boldsymbol{W}^{*}$ obtained through MNN. For the closed-loop system (\ref{e17}), if there exist positive-definite matrices $P$, $U$, and matrices $Q_1,Q_2$, scalars $\varpi_i>0$, appropriate matrices $M_1$, $M_2$, $M_3$, $N$, and $L>0$, scalar $\beta_1>0$, with suitable switching controller gain matrix $K$, such that the following matrix inequalities are feasible:
		\begin{equation}\label{e31}
		\tilde{\Xi}>0,
	\end{equation}
		\begin{equation}\label{e29}
		\boldsymbol{\Phi}_{1}<0,\ \boldsymbol{\Phi}_{2}<0,\
		\boldsymbol{\Phi}_{3}<0,\ \boldsymbol{\Phi}_{4}<0,
	\end{equation}
	where
	\begin{equation*}\begin{aligned}
			\boldsymbol{\Phi}_{1}&=
			\begin{bmatrix}
				\eta_ {1,11}&\eta_ {1,12}&\eta_ {1,13}&0&\eta_ {1,15}&P_{11}\\
				\ast&\eta_ {1,22}&\eta_ {1,23}&0&\eta_ {1,25}&N\\
				\ast&\ast&\eta_ {1,33}&0&0&0\\
				\ast&\ast&\ast&-\Lambda&0&0\\
				\ast&\ast&\ast&\ast&-L&0\\
				\ast&\ast&\ast&\ast&\ast&-\beta_1\\
			\end{bmatrix},\\
			\boldsymbol{\Phi}_{2}&=
			\begin{bmatrix}
				\eta_ {1,11}&\eta_ {2,12}&\eta_ {1,13}&0&\eta_ {1,15}&P_{11}&\eta_ {2,17}\\
				\ast&\eta_ {2,22}&\eta_ {2,23}&0&\eta_ {1,25}&N&\eta_ {2,27}\\
				\ast&\ast&\eta_ {1,33}&0&0&0&\eta_ {2,37}\\
				\ast&\ast&\ast&-\Lambda&0&0&0\\
				\ast&\ast&\ast&\ast&-L&0&0\\
				\ast&\ast&\ast&\ast&\ast&-\beta_1&0\\
				\ast&\ast&\ast&\ast&\ast&\ast&\eta_ {2,77}\\
			\end{bmatrix},\\
			\boldsymbol{\Phi}_{3}&=
			\begin{bmatrix}
				\eta_ {1,11}&\eta_ {1,12}&\eta_ {1,13}&\eta_ {3,14}&\eta_ {1,15}&P_{11}\\
				\ast&\eta_ {1,22}&\eta_ {1,23}&\eta_ {3,24}&\eta_ {1,25}&N\\
				\ast&\ast&\eta_ {1,33}&0&0&0\\
				\ast&\ast&\ast&-\Lambda&0&0\\
				\ast&\ast&\ast&\ast&-L&0\\
				\ast&\ast&\ast&\ast&\ast&-\beta_1\\
			\end{bmatrix},\\
			\boldsymbol{\Phi}_{4}&=
			\begin{bmatrix}
				\eta_ {1,11}&\eta_ {2,12}&\eta_ {1,13}&\eta_ {3,14}&\eta_ {1,15}&P_{11}&\eta_ {4,17}\\
				\ast&\eta_ {2,22}&\eta_ {2,23}&\eta_ {3,24}&\eta_ {1,25}&N&\eta_ {2,27}\\
				\ast&\ast&\eta_ {1,33}&0&0&0&\eta_ {2,37}\\
				\ast&\ast&\ast&-\Lambda&0&0&0\\
				\ast&\ast&\ast&\ast&-L&0&0\\
				\ast&\ast&\ast&\ast&\ast&-\beta_1&0\\
				\ast&\ast&\ast&\ast&\ast&\ast&\eta_ {2,77}\\
			\end{bmatrix},\\
			\tilde{\Xi}&=
			\begin{bmatrix}
				\frac{hQ_1+\ast}{2}+P_{11}&h(-Q_1+Q_2)+P_{12}\\
				h(-Q_1^T+Q_2^T)+P_{21}&h(\frac{Q_1-2Q_2}{2}+\ast)+P_{22}\\
			\end{bmatrix},\\
			\eta_ {1}&=[P_{11}(A_{s}-B_{2,s}KC_s)+\ast]-(\frac{Q_1}{2}+\ast)\\
			&\ +\beta_1\delta^2- \boldsymbol{\zeta}_1^TL \boldsymbol{\zeta}_2+2M_1+\epsilon C_s^{T}(t)\Lambda  C_s(t),\\
			\eta_{1,11}&=\eta_{1}+P_{11}B_{1,s}(PB_{1,s})^T,\\
			\eta_{1,12}&=\frac{1}{2}h(Q_1+\ast)+(NA_s)^T-(NB_{2,s}KC_s)^T+M_2^T,\\
			\eta_{1,13}&=Q_1-Q_2-M_1+M_3^T,\\			
			\eta_ {1,15}&=P_{11}\boldsymbol{W}^*+\frac{1}{2}\boldsymbol{\zeta}_3^T L^T,\\
			\eta_{1,22}&=hU-2N+NB_{1,s}(NB_{1,s})^T,\\
			\eta_{1,23}&=h(-Q_1+Q_2)-M_2+P^T_{21},	\ \eta_ {1,25}=\boldsymbol{\zeta}_4^T\boldsymbol{\Omega}^T+N\boldsymbol{W}^*,\\
			\eta_ {1,33}&=(Q_2-\frac{Q_1}{2}+\ast)-2M_3,\\
			\eta_{2,12}&=(NA_s)^T-(NB_{2,s}KC_s)^T+M_2^T,\\
			\eta_{2,17}&=h(P+N)B_{2,s}KC_s-hM_1,\\
			\eta_{2,22}&=-2N+NB_{1,s}(NB_{1,s})^T,\\
			\eta_{2,23}&=-M_2+P^T_{21},\ \eta_{2,27}=-hM_2,\\
			\eta_{2,37}&=-hM_3,\ 	\eta_{2,77}=-he^{-2\alpha h}U,\\
			\eta_ {3,14}&=-P_{11}B_{2,s}K,\ \eta_{3,24}=-NB_{2,s}K,\\
			\eta_{4,17}&=-hM_1,\ \boldsymbol{\zeta}_1=G_{min}\boldsymbol{V}^*,\  \boldsymbol{\zeta}_2=G_{max}\boldsymbol{V}^*,\\ \boldsymbol{\zeta}_3&=(G_{min}+G_{max})\boldsymbol{V}^*,\ \boldsymbol{\zeta}_4=G\boldsymbol{V}^*.
	\end{aligned}\end{equation*}
	
	Then the closed-loop system (\ref{e17}) is SGUUB.
\end{theorem}
\begin{IEEEproof} Firstly, condition (\ref{e31}) implies that for any $\chi(t)$, with positive matrices $P,U>0$, and scalar $\varpi_i>0$, there exists inequalities that holds true
	\begin{equation*}
		\bar{\zeta}_1\|\xi_s(t)\|^2\leq V(t)\leq \bar{\zeta}_2\|\xi_s(t)\|^2,
	\end{equation*}
	where $\bar{\zeta}_1>0$, $\bar{\zeta}_2>0$. 
	For switching controller design (\ref{e16}), we make the following discussion.
	
	From (\ref{e15}), it is obviously satisfied that
	\begin{equation*}
		\epsilon y_s^{T}(t)\Lambda  y_s(t)-e^T(t)\Lambda e(t)\geq0.
	\end{equation*}
	
	The inequality below is valid for any matrix $M_1,M_2,M_3\in\mathbb{R}^{m\times m}$:
	\begin{equation}\label{e32}
		\begin{aligned}
			0&=2[\xi_s^T(t)M_1+\dot{\xi}_s^T(t)M_2+\xi_s^T(t-\tau)M_3]\\
			\ &\times [\xi_s(t)-\xi_s(t-\tau)-\tau(t)\nu_1(t)]
			.\end{aligned}
	\end{equation}
	
	If $\chi(t)=1$, the inequality below is valid for any matrix $N\in\mathbb{R}^{m\times m}$:
	\begin{equation}\label{e33}
		\begin{aligned}
			0&=2\dot{\xi}_s^T(t)N[(A_{s}-B_{2,s}KC_s)\xi_{s}(t)+	\boldsymbol{W}^{*}\boldsymbol{\mu}(\boldsymbol{\varsigma}(t))\\&\ +\Delta\boldsymbol{A}(\boldsymbol{\xi}_{s}(t))
			+ B_{2,s}KC_s\tau(t)\nu_1(t)+B_{1,s}d(t)-\dot{\xi}_s(t)].\\
		\end{aligned}
	\end{equation}
	
	According to (\ref{e18})-(\ref{e25}) and (\ref{e32})-(\ref{e33}), the inequality mentioned below stands true:
	\begin{equation*}
		\begin{aligned}
			&\ \dot{V}(t)\leq 2\xi_s^T(t-\tau)P_{21}\dot{\xi}_s(t)\\
			&\ +\xi_{s}^{T}(t)[P_{11}(A_{s}-B_{2,s}KC_s)+\ast]\xi_{s}(t)\\
			&\ +2\xi_{s}^{T}(t)P_{11}B_{2,s}KC_s\tau(t)\nu_1(t)+2\xi_{s}^{T}(t)P_{11}\boldsymbol{W}^{*}\boldsymbol{\mu}(\boldsymbol{\varsigma}(t))\\
			&\ +2\xi_{s}^{T}(t)P_{11}\Delta\boldsymbol{A}(\xi_{s}(t))+2\xi_{s}^{T}(t)P_{11}B_{1,s}d(t) +(h-\tau(t))\\
			&\ \times \dot{\xi}^T(t)U\dot{\xi}(t)-[2\alpha(h-\tau(t))+1]e^{-2\alpha h}\tau(t)\nu_{1}^{T}(t)U\nu_{1}(t)\\
													\end{aligned}
		\end{equation*}
	\begin{equation}\label{e34}
	\begin{aligned}
			&\ -\begin{bmatrix}\xi_s(t)\\\xi_s(t-\tau(t))\end{bmatrix}^T\begin{bmatrix}\frac{Q_1}{2}+\ast&-Q_1+Q_2\\\\\ast&\frac{-2Q_2+Q_1}{2}+\ast\end{bmatrix}\\
			&\ \times \begin{bmatrix}\xi_s(t)\\\xi_s(t-\tau(t))\end{bmatrix} +(h-\tau(t))[\dot{\xi}_s^{T}(t)(Q_1+\ast)\xi_s(t) \\
			&\ +2\dot{\xi}_s^T(t)(-Q_1+Q_2)\xi_s(t-\tau(t))] +2\boldsymbol{\mu}^T(\boldsymbol{\varsigma})\boldsymbol{\Omega}G\boldsymbol{V}^*\dot{\xi}_s(t)\\
			&\ +\beta _1[\delta^2 \xi_s^T(t)\xi_s(t)-\Delta\boldsymbol{A}^T(\xi_s)\Delta\boldsymbol{A}(\xi_s)]\\
			&\ 	+[\boldsymbol{\mu}(\boldsymbol{\varsigma} )-G_{min}\boldsymbol{\varsigma} ]^TL[G_{max}\boldsymbol{\varsigma}-\boldsymbol{\mu}(\boldsymbol{\varsigma} )]\\
			&\ +2\dot{\xi}_s^T(t)N[(A_{s}-B_{2,s}KC_s)\xi_{s}(t)+	\boldsymbol{W}^{*}\boldsymbol{\mu}(\boldsymbol{\varsigma}(t))\\&\ +\Delta\boldsymbol{A}(\xi_{s}(t))
			+B_{2,s}KC_s\tau(t)\nu_1(t)-\dot{\xi}_s(t)]\\
			&\ +2[\xi_s^T(t)M_1+\dot{\xi}_s^T(t)M_2+\xi_s^T(t-\tau)M_3]\\ &\ \times [\xi_s(t)-\xi_s(t-\tau)-\tau(t)\nu_1(t)]\\
			&\ +\epsilon \xi_s^{T}(t)C_s^T\Lambda  C_s\xi_s(t)-e^T(t)\Lambda e(t).
		\end{aligned}
	\end{equation}
	
	Then, it can be derived:
	\begin{equation}\label{e35}
		\begin{aligned}
			\dot{V}(t)&\leq \bar{\xi}_{1}^T\Phi_1\bar{\xi}_{1} +2\dot{\xi}_{s}^T(t)NB_{1,s}d(t)\\
			&\ +2{\xi}_{s}^T(t)P_{11}B_{1,s}d(t)
			,\end{aligned}\end{equation}
	where
	\begin{equation*}\begin{aligned}
			\bar\xi_{1}^T=&[\xi_{s}^T(t),\dot{\xi}_{s}^T(t),\xi_{s}^T(t-\tau),e^T(t),\\
			&\boldsymbol{\mu}^T(t),\Delta\boldsymbol{A}^T(\xi_{s}(t)),\nu_1^T(t)].
	\end{aligned}\end{equation*}
	
	\begin{equation*}
		\Phi_{1}=
		\begin{bmatrix}
			\eta_ {1}&\eta_ {2}&\eta_ {1,13}&0&\eta_ {1,15}&P&\eta_ {3}\\
			\ast&\eta_ {4}&\eta_ {5}&0&\eta_ {1,25}&N&\eta_ {6}\\
			\ast&\ast&\eta_ {1,33}&0&0&0&\eta_ {7}\\
			\ast&\ast&\ast&-\Lambda&0&0&0\\
			\ast&\ast&\ast&\ast&-L&0&0\\
			\ast&\ast&\ast&\ast&\ast&-\beta_1&0\\
			\ast&\ast&\ast&\ast&\ast&\ast&\eta_ {8}\\
		\end{bmatrix},
	\end{equation*}
	\begin{equation*}\begin{aligned}
			\eta_ {2}&=\frac{1}{2}(h-\tau(t))(Q_1+\ast)+(NA_s)^T\\
			&\ -(NB_{2,s}KC_s)^T+M_2^T,\\
			\eta_ {3}&=(P_{11}+N)B_{2,s}KC_s\tau(t)-M_1\tau(t),\\
			\eta_ {4}&=(h-\tau(t))U-2N,\\ 
			\eta_ {5}&=(h-\tau(t))(-Q_1+Q_2)-M_2+P^T_{21},\\
			\eta_ {6}&=-M_2\tau(t),\ \eta_ {7}=-M_3\tau(t),\\
			\eta_ {8}&=-[2\alpha(h-\tau(t))+1]e^{-2\alpha h}\tau(t)U.\\
	\end{aligned}\end{equation*}
	
	According to (\ref{e26}) and (\ref{e35}), it can be satisfied:
	\begin{equation}\label{e36}
		\begin{aligned}
			\dot{V}(t)\leq \bar{\xi}_{1}^T\Phi_2\bar{\xi}_{1}+2D_{1}^{2}
			,\end{aligned}\end{equation}
	where \begin{equation*}
		\Phi_{2}=
		\begin{bmatrix}
			\eta_ {1,11}&\eta_ {2}&\eta_ {1,13}&0&\eta_ {1,15}&P_{11}&\eta_ {3}\\
			\ast&\eta_ {9}&\eta_ {5}&0&\eta_ {1,25}&N&\eta_ {6}\\
			\ast&\ast&\eta_ {2,33}&0&0&0&\eta_ {7}\\
			\ast&\ast&\ast&-\Lambda&0&0&0\\
			\ast&\ast&\ast&\ast&-L&0&0\\
			\ast&\ast&\ast&\ast&\ast&-\beta_1&0\\
			\ast&\ast&\ast&\ast&\ast&\ast&\eta_ {8}\\
		\end{bmatrix},\end{equation*}
	where $\eta_{9}=\eta_{4}+NB_{1,s}(NB_{1,s})^T$.
	
	If condition (\ref{e29}) is satisfied, this implies that for matrix $\Phi_{2}$, both cases where $\tau(t)\to 0$ and $\tau(t)\to h$ fulfill the negative-definite condition, leading to the conclusion that:
	\begin{equation*}
		\Phi_{2}<0,
	\end{equation*}
	thus $\lambda(\phi_2)<0$, then we can extract from (\ref{e36}):
	\begin{equation*}
		\dot{V}(t)\leq-2\alpha_1 \left\| \bar{\xi}_1\right\|^2+2D_{1}^{2}
		,\end{equation*}
	where $\alpha_1=\frac{1}{2}|\lambda_{\max}(\phi_2)|>0$.
	
	And then it can be gotten:
	\begin{equation*}
		\dot{V}(t)\leq-2\alpha_1 \left\| \xi_s\right\|^2+2D_{1}^{2}
		,\end{equation*}		
	it can be satisfied:
	\begin{equation*}
		\dot{V}(t)\leq-2\alpha_1 (\alpha_2+\overline{\alpha}_3)^{-1} V(t)+2D_{1}^{2}
		,\end{equation*}		
	where \begin{equation*}\begin{aligned}
			\alpha_2&=\max\{\lambda_{\max}P,\ \lambda_{\max}Q,\ \lambda_{\max}U\},\\ \overline{\alpha}_3&=\sup_{\xi_s(t)\neq0}\frac{\sum_{i=1}^{n_h}2\omega_i\int_0^{ \varsigma_i(t)}\mu_i(\rho)(g_{i,\max}-g_{i,\min})d\rho}{\|\xi_s(t)\|^2},\\
			Q&=
			\begin{bmatrix}
				\frac{Q_1+\ast}{2}&(-Q_1+Q_2)\\
				(-Q_1^T+Q_2^T)&(\frac{Q_1-2Q_2}{2}+\ast)\\
			\end{bmatrix}.\\
	\end{aligned}	\end{equation*}
	
	From Lemma 1, one has:
	\begin{equation}\label{e37}
		V(t)\leq V(0)e^{-2\bar{\alpha}_1 t}+\frac{D_{1}^{2}}{\bar{\alpha}_1}(1-e^{-2\bar{\alpha}_1 t}),
	\end{equation}
	where $\bar{\alpha}_1=\alpha_1 (\alpha_2+\overline{\alpha}_3)^{-1}$.
	
	Else if $\chi(t)=0$, the inequality below is valid for any matrix $N\in\mathbb{R}^{m\times m}$:
	\begin{equation*}\begin{aligned}
			0&=2\dot{\xi}_s^T(t)N[(A_{s}-B_{2,s}KC_s)\xi_{s}(t)+	\boldsymbol{W}^{*}\boldsymbol{\mu}(\boldsymbol{\varsigma}(t))\\&\ +\Delta\boldsymbol{A}(\xi_{s}(t))
			-B_{2,s}Ke(t)+B_{1,s}d(t)-\dot{\xi}_s(t)].\\
		\end{aligned}
	\end{equation*}
	
	Similar to the above analysis, it can be obtained
	\begin{equation}\label{e38}
		\begin{aligned}
			\dot{V}(t)&\leq \bar{\xi}_{1}^T\Phi_3\bar{\xi}_{1} +2\dot{\xi}_{s}^T(t)NB_{1,s}d(t)\\
			&\ +2{\xi}_{s}^T(t)PB_{1,s}d(t)
			,\end{aligned}\end{equation}
	where 	\begin{equation*}\begin{aligned}
			\Phi_{3}&=
			\begin{bmatrix}
				\eta_ {1}&\eta_ {2}&\eta_ {1,13}&\eta_ {3,14}&\eta_ {1,15}&P_{11}&\eta_ {10}\\
				\ast&\eta_ {4}&\eta_ {5}&\eta_ {33}&\eta_ {1,25}&N&\eta_ {6}\\
				\ast&\ast&\eta_ {1,33}&0&0&0&\eta_ {7}\\
				\ast&\ast&\ast&-\Lambda&0&0&0\\
				\ast&\ast&\ast&\ast&-L&0&0\\
				\ast&\ast&\ast&\ast&\ast&-\beta_1&0\\
				\ast&\ast&\ast&\ast&\ast&\ast&\eta_ {8}\\
			\end{bmatrix},\\
			\eta_ {10}&=-M_1 \tau(t).\\
		\end{aligned}
	\end{equation*}

	From inequality (\ref{e28}), it is concluded that:
	\begin{equation*}\begin{aligned}
			\dot{V}(t)\leq \bar{\xi}_{1}^T\Phi_4\bar{\xi}_{1}+2D_{1}^{2}
			,\end{aligned}\end{equation*}
	where \begin{equation*}
		\Phi_{4}=
		\begin{bmatrix}
			\eta_ {1,11}&\eta_ {2}&\eta_ {1,13}&\eta_ {3,14}&\eta_ {1,15}&P_{11}&\eta_ {10}\\
			\ast&\eta_ {9}&\eta_ {5}&\eta_ {3,24}&\eta_ {1,25}&N&\eta_ {6}\\
			\ast&\ast&\eta_ {1,33}&0&0&0&\eta_ {7}\\
			\ast&\ast&\ast&-\Lambda&0&0&0\\
			\ast&\ast&\ast&\ast&-L&0&0\\
			\ast&\ast&\ast&\ast&\ast&-\beta_1&0\\
			\ast&\ast&\ast&\ast&\ast&\ast&\eta_ {8}\\
		\end{bmatrix}.
	\end{equation*}

	If condition (\ref{e29}) is satisfied, this implies that for matrix $\Phi_{4}$, both cases where $\tau(t)\to 0$ and $\tau(t)\to h$ fulfill the negative-definite condition, leading to the conclusion that:
	\begin{equation*}
		\Phi_{4}<0.
	\end{equation*}
	
	Similarly, we can derive that
	\begin{equation}\label{e39}
		\dot{V}(t)\leq-2\alpha_4 \left\| \bar{\xi}_1\right\|^2+2D_{1}^{2}
		,\end{equation}
	where $\alpha_4=\frac{1}{2}|\lambda_{\max}(\phi_4)|>0$.
	
	According to (\ref{e39}), it is similar to obtain the inequality:
	\begin{equation}\label{e40}
		V(t)\leq V(0)e^{-2\bar{\alpha}_2 t}+\frac{D_{1}^{2}}{\bar{\alpha}_2}(1-e^{-2\bar{\alpha}_2 t}),
	\end{equation}
	where $\bar{\alpha}_2=\alpha_4 (\alpha_2+\overline{\alpha}_3)^{-1}$.
	
	In summary, from (\ref{e37}) and (\ref{e40}), the closed-loop system (\ref{e17}) is SGUUB. The proof is complete.
\end{IEEEproof}

By analyzing the Lyapunov functional, we obtain the stability criterion of the system (\ref{e17}). The switching event-triggered control problem is transformed into a BMI problem, but its feasibility cannot be resolved immediately using the MATLAB toolbox. Therefore, we will present a novel algorithm to convert the BMI problem in Theorem 1 into a LMI problem.

Define $\Xi_{1}=NB_{2,s}K$, $\Xi_{2}=P_{11}B_{2,s}K$. From $(\boldsymbol{\beta}_{2}^{(1/2)}\boldsymbol{\zeta}_1+\beta_{2}^{-(1/2)}\boldsymbol{\zeta}_2)^T L(\beta_{2}^{(1/2)}\boldsymbol{\zeta}_1+\beta_{2}^{-(1/2)}\boldsymbol{\zeta}_2)~\ge~0$, it can be derived:
\begin{equation*}-\boldsymbol{\zeta}_1^TL\boldsymbol{\zeta}_2-\boldsymbol{\zeta}_2^TL\boldsymbol{\zeta}_1\leq\beta_2\boldsymbol{\zeta}_1^TL\boldsymbol{\zeta}_1+\beta_2^{-1}\boldsymbol{\zeta}_2^TL\boldsymbol{\zeta}_2.\end{equation*}

From the Schur complement, the following transformation can be carried out through the method of variable substitution. Convert the BMIs (\ref{e29}) in Theorem 1 to LMIs:

\begin{equation}\label{e41}
	\begin{aligned}
	\tilde{\Phi}_{1}&=
	\begin{bmatrix}
		\Psi_{1,11}&\Psi_{1,12}\\
		\ast&\Psi_{1,22}
	\end{bmatrix}<0,\ 	\tilde{\Phi}_{2}=
	\begin{bmatrix}
		\Psi_{2,11}&\Psi_{1,12}\\
		\ast&\Psi_{1,22}
	\end{bmatrix}<0,\\
%
	\tilde{\Phi}_{3}&=
	\begin{bmatrix}
		\Psi_{3,11}&\Psi_{1,12}\\
		\ast&\Psi_{1,22}
	\end{bmatrix}<0,\ 	\tilde{\Phi}_{4}=
	\begin{bmatrix}
		\Psi_{4,11}&\Psi_{1,12}\\
		\ast&\Psi_{1,22}
	\end{bmatrix}<0,
	\end{aligned}
\end{equation}
where
\begin{equation*}\begin{aligned}
		\Psi_{1,11}=&
		\begin{bmatrix}
			\bar{\Psi}_{1,11}&\eta_ {1,12}&\eta_ {1,13}&0&\eta_ {1,15}&P_{11}\\
			\ast&\bar{\Psi}_{1,22}&\eta_ {1,23}&0&\eta_ {1,25}&N\\
			\ast&\ast&\eta_ {1,33}&0&0&0\\
			\ast&\ast&\ast&-\Lambda&0&0\\
			\ast&\ast&\ast&\ast&-L&0\\
			\ast&\ast&\ast&\ast&\ast&-\beta_1\\
		\end{bmatrix},\\
		\Psi_{2,11}=&
		\begin{bmatrix}
			\bar{\Psi}_{1,11}&\eta_ {2,12}&\eta_ {1,13}&0&\eta_ {1,15}&P_{11}&\eta_ {2,17}\\
			\ast&-2N&\eta_ {2,23}&0&\eta_ {1,25}&N&\eta_ {2,27}\\
			\ast&\ast&\eta_ {1,33}&0&0&0&\eta_ {2,37}\\
			\ast&\ast&\ast&-\Lambda&0&0&0\\
			\ast&\ast&\ast&\ast&-L&0&0\\
			\ast&\ast&\ast&\ast&\ast&-\beta_1&0\\
			\ast&\ast&\ast&\ast&\ast&\ast&\eta_ {2,77}\\
		\end{bmatrix},\\
		\Psi_{3,11}=&
		\begin{bmatrix}
			\bar{\Psi}_{1,11}&\eta_ {1,12}&\eta_ {1,13}&\eta_ {3,14}&\eta_ {1,15}&P_{11}\\
			\ast&\bar{\Psi}_{1,22}&\eta_ {1,23}&\eta_ {3,24}&\eta_ {1,25}&N\\
			\ast&\ast&\eta_ {1,33}&0&0&0\\
			\ast&\ast&\ast&-\Lambda&0&0\\
			\ast&\ast&\ast&\ast&-L&0\\
			\ast&\ast&\ast&\ast&\ast&-\beta_1\\
		\end{bmatrix},\\
		\Psi_{4,11}=&
		\begin{bmatrix}
			\bar{\Psi}_{1,11}&\eta_ {2,12}&\eta_ {1,13}&\eta_ {3,14}&\eta_ {1,15}&P_{11}&\eta_ {4,17}\\
			\ast&-2N&\eta_ {2,23}&\eta_ {3,24}&\eta_ {1,25}&N&\eta_ {2,27}\\
			\ast&\ast&\eta_ {1,33}&0&0&0&\eta_ {2,37}\\
			\ast&\ast&\ast&-\Lambda&0&0&0\\
			\ast&\ast&\ast&\ast&-L&0&0\\
			\ast&\ast&\ast&\ast&\ast&-\beta_1&0\\
			\ast&\ast&\ast&\ast&\ast&\ast&\eta_ {2,77}\\
		\end{bmatrix},\\
			\bar{\Psi}_{1,12}=&
		\begin{bmatrix}
			\delta I_m&\boldsymbol{\zeta}_1^T&\boldsymbol{\zeta}_2^T\\
			0_{m\times m}&0_{m\times n_h}&0_{m\times n_h}\\
			0_{4m\times m}&0_{4m\times n_h}&0_{4m\times n_h}\\
		\end{bmatrix},\\
		\bar{\Psi}_{1,13}=&
		\begin{bmatrix}
			C_s^T&P_{11}B_{1,s}&0_{m\times 1}\\
			0_{m\times 1}&0_{m\times 1}&NB_{1,s}\\
			0_{4m\times 1}
			&0_{4m\times 1}&0_{4m\times 1}\\
		\end{bmatrix},\\
		\Psi_{1,22}=&
		\begin{bmatrix}
			-\bar{\beta}_1&0&0&0&0&0\\
			\ast&-\bar{\beta}_2&0&0&0&0\\
			\ast&\ast&-\bar{\beta}_2&0&0&0\\
			\ast&\ast&\ast&-\Lambda_1&0&0\\
			\ast&\ast&\ast&\ast&-I_m&0\\
			\ast&	\ast&\ast&\ast&\ast&-I_m\\
		\end{bmatrix},\\
			\Psi_{1,12}=&
	\begin{bmatrix}
		\bar{\Psi}_{1,12}&\bar{\Psi}_{1,13}\\
	\end{bmatrix},\	\bar{\Psi}_{1,22}=hU-2N,\\
		\bar{\Psi}_{1,11}=&[P_{11}A_{s}-\Xi_2 C_s+\ast]-(\frac{Q_1}{2}+\ast)+2M_1,\\
		\bar{\beta}_1=&\beta_{1}^{-1}I_m,\ \bar{\beta}_2=\beta_{2}^{-1}L^{-1},\ \Lambda_1=( \epsilon \Lambda)^{-1}I_{n_h}.
\end{aligned}\end{equation*}
\begin{remark}
	It is evident that the solutions to the LMIs (\ref{e31}) and (\ref{e41}) do not necessarily align with the solutions to the BMIs (\ref{e29}) and (\ref{e31}). However, the intricacy of the cross-terms involving $\Xi_1$ and $\Xi_2$ poses a significant challenge in directly deriving a corresponding LMI. Consequently, we turn to Algorithm \ref{alg.1} to obtain an indirect LMI problem, ensuring that the solutions to this indirect LMI fall within the feasible region defined by the solutions outlined in Theorem 1.
\end{remark}
\begin{algorithm}[H]
	\caption{Solving controller gain $K$ of LMIs (\ref{e31})(\ref{e41}).}\label{alg:alg2}
	\begin{algorithmic}
		\STATE 
		\STATE {\textsc{Step 1}}
		\STATE \hspace{0.5cm}$ \textbf{Set up the relationship. } $
		\STATE \hspace{0.5cm} Letting $N=r_1P_{11}$, where $r_1\in\mathbb{R}^{m\times m}.$ 
		\STATE {\textsc{Step 2}}
		\STATE \hspace{0.5cm}$ \textbf{Searching for an Appropriate $r_1$. } $
		\STATE \hspace{0.5cm} Solving the LMIs (\ref{e31}) (\ref{e41}) using $N$ directly. Given that $P \neq N$, we can derive:
		%
		%
		
		\   \  \ \ \ If feasible, $r_1=NP_{11}^{-1}$, proceed to the next step;
		
			\   \  \ \ \ \  Else if infeasible, adjust the parameters and solve again.
		\STATE {\textsc{Step 3}}
		\STATE \hspace{0.5cm}$ \textbf{Refine the LMIs. } $
		\STATE \hspace{0.5cm} Substituting the $r_1$ obtained in the previous step, we obtain $N=r_1P_{11}$, $\Xi_1=r_1\Xi_2$. By incorporating this into the LMIs (\ref{e31}) (\ref{e41}), we can derive LMIs whose feasible solutions lie within the feasible region of the BMIs (\ref{e29}) and LMI (\ref{e31}), and then solve LMIs (\ref{e31}) (\ref{e41}).
%
		
		\ \	\ \ \ \ \  If feasible, $	K=B_{2,s}^{-1}P_{11}^{-1}\Xi_2$, proceed to the next step;
		
		\   \  \ \ Else if infeasible, return to Step 2 and readjust the parameters.
		\STATE {\textsc{Step 4}}
		\STATE \hspace{0.5cm}$ \textbf{Outout $K$.} $
		
	\end{algorithmic}
	\label{alg.1}
\end{algorithm}

\begin{theorem}\label{thm2} If the LMIs (\ref{e31}) (\ref{e41}) have feasible solution by Algorithm \ref{alg.1}, then the appropriate controller gain can be solved 
	\begin{equation*}
		K=B_{2,s}^{-1}P_{11}^{-1}\Xi_2,
	\end{equation*}
	to ensure the closed-loop system (\ref{e17}) be SGUUB by switching event-triggered control.
\end{theorem}

\begin{proposition}
	According to Theorem 1, there exists an appropriate controller gain matrix $K$ that enables the system to achieve SGUUB, and satisfies:
	\begin{equation}\label{e43}
		\begin{aligned}
			-\sqrt{\frac{2}{\bar{\alpha}_6}}D_1\leq\parallel \xi_s(t)\parallel\leq\sqrt{\frac{2}{\bar{\alpha}_6}}D_1,
	\end{aligned}\end{equation}
	where\begin{equation*} \begin{aligned} \bar{\alpha}_6&=\bar{\alpha}_5(\alpha_5+\underline{\alpha}_3),\
			\bar{\alpha}_4=\min\{\alpha_{1},\ \alpha_4\},\\ \bar{\alpha}_5&=\bar{\alpha}_4 (\alpha_2+\overline{\alpha}_3)^{-1},\ 
			\alpha_5=\min\{\lambda_{\min} P,\lambda_{\min} Q,\lambda_{\min} U\},\\
			\underline{\alpha}_3&=\inf_{\xi_s(t)\neq0}\frac{\sum_{i=1}^{n_h}2\omega_i\int_0^{ \varsigma_i(t)}\mu_i(\rho)(g_{i,\max}-g_{i,\min})d\rho}{\|\xi_s(t)\|^2}.\\
	\end{aligned}\end{equation*}	
\end{proposition}
\begin{IEEEproof}
	For any $\chi(t)$, one has
	\begin{equation*}
		\dot{V}(t)\leq-2\bar{\alpha}_4 \left\| \bar{\xi}_1\right\|^2+2D_{1}^{2}
		,\end{equation*}
	then it can be obtained		
	\begin{equation*}
		\dot{V}(t)\leq-2\bar{\alpha}_4 (\alpha_2+\overline{\alpha}_3)^{-1} V(t)+2D_{1}^{2}
		.\end{equation*}
	From Lemma 1, we can get:
	\begin{equation}\label{e44}
		V(t)\leq V(0)e^{-2\bar{\alpha}_5 t}+\frac{2D_{1}^{2}}{\bar{\alpha}_5}(1-e^{-2\bar{\alpha}_5 t}),
	\end{equation} 
	Consider (\ref{e44}), it can be derived:
	\begin{equation*}\begin{aligned}
			(\alpha_5+\underline{\alpha}_3)\parallel \xi_s(t)\parallel^2&\leq(\alpha_2+\overline{\alpha}_3)\parallel \xi_s(0)\parallel^2e^{-2\bar{\alpha}_5 t}\\
			&\ +\frac{2D_{1}^{2}}{\bar{\alpha}_5}(1-e^{-2\bar{\alpha}_5 t}),
		\end{aligned}
	\end{equation*}
	Then it can be gotten:
	\begin{equation*}\begin{aligned}
			\parallel \xi_s(t)\parallel&\leq\sqrt{(\alpha_6\parallel \xi_s(0)\parallel^2-\frac{2D_{1}^{2}}{\bar{\alpha}_6})e^{-2\bar{\alpha}_5 t}+\frac{2D_{1}^{2}}{\bar{\alpha}_6}},\\
				\parallel \xi_s(t)\parallel&\geq-\sqrt{(\alpha_6\parallel \xi_s(0)\parallel^2-\frac{2D_{1}^{2}}{\bar{\alpha}_6})e^{-2\bar{\alpha}_5 t}+\frac{2D_{1}^{2}}{\bar{\alpha}_6}},
		\end{aligned}
	\end{equation*}
	where $\alpha_6=\frac{(\alpha_2+\overline{\alpha}_3)}{	(\alpha_5+\underline{\alpha}_3)}$. From this, the boundary at $t=T\to\infty$ can be derived:
	\begin{equation*}\begin{aligned}
				-\sqrt{\frac{2}{\bar{\alpha}_6}}D_1\leq\parallel \xi_s(t)\parallel\leq\sqrt{\frac{2}{\bar{\alpha}_6}}D_1.
	\end{aligned}\end{equation*}
Then the system state boundaries of closed-loop system (\ref{e17}) under switching event triggered controller are obtained. The proof is complete.
\end{IEEEproof}

 Define $\pounds_1=[I_{2},0_{2\times 2},0_{2\times 2},0_{2\times 1},0_{2\times 15},0_{2\times 1},0_{2\times 2}]$, $\pounds_2=[0_{2\times 2},I_{2},0_{2\times 2},0_{2\times 1},0_{2\times 15},0_{2\times 1},0_{2\times 2}]$. Let us now analyze the effects of using switching event-triggered controller in absence of disturbances.
\begin{corollary}
	Consider $d(t)\equiv 0$. Given scalars $\alpha$, $\epsilon$, $\beta_2$, matrix $\Lambda$, and the weights $\boldsymbol{V}^{*},\boldsymbol{W}^{*}$ obtained through MNN. For closed-loop system (\ref{e17}), if there exist positive-definite matrices $P$, $U>0$, scalar $\varpi_i>0$,  appropriate matrices $M_1$, $M_2$, $M_3$, $N$, and $L>0$, scalar $\beta_1>0$, with suitable switching controller gain matrix $K$, such that the following inequalities are feasible: 
	\begin{equation}\label{eq1}
		\hat{\boldsymbol{\Phi}}_1<0,\ \hat{\boldsymbol{\Phi}}_2<0,\ \hat{\boldsymbol{\Phi}}_3<0,\ \hat{\boldsymbol{\Phi}}_4<0,\ 
	\end{equation}
	where \begin{equation*}
		\begin{aligned}
			\hat{\boldsymbol{\Phi}}_1&=\boldsymbol{\Phi}_1-\pounds_1^TP_{11}B_{1,s}(P_{11}B_{1,s})^T\pounds_1,\\
	\hat{\boldsymbol{\Phi}}_2&=\boldsymbol{\Phi}_2-\pounds_2^TNB_{1,s}(NB_{1,s})^T\pounds_2,\\
		\hat{\boldsymbol{\Phi}}_3&=\boldsymbol{\Phi}_3-\pounds_1^TP_{11}B_{1,s}(P_{11}B_{1,s})^T\pounds_1,\\
	\hat{\boldsymbol{\Phi}}_4&=\boldsymbol{\Phi}_4-\pounds_2^TNB_{1,s}(NB_{1,s})^T\pounds_2,\\
		\end{aligned}
	\end{equation*}
	with LMI (\ref{e31}), then the closed-loop system (\ref{e17}) is exponentially stable.
\end{corollary}
\begin{IEEEproof}
	From (\ref{e35}), when $d(t)\equiv0$, the following inequality can obviously be derived if $\chi(t)=1$:
	\begin{equation*}
		\dot{V}(t)\leq \bar{\xi}_{1}^T\Phi_1\bar{\xi}_{1}.
	\end{equation*} 
	If (\ref{eq1}) is satisfied, it can be obtained:
	\begin{equation*}
		\dot{V}(t)\leq-2\tilde{\alpha}_1 \left\| \bar{\xi}_1\right\|^2,
	\end{equation*}
	where $\tilde{\alpha}_1=\frac{1}{2}|\lambda_{\max}(\Phi_1)|$.
	Then we can get:
	\begin{equation*}
		\dot{V}(t)\leq-2\tilde{\alpha}_1 \left\| \xi_s(t)\right\|^2.
	\end{equation*}
	And then it can be derived:
	\begin{equation}\label{eq2}
		\dot{V}(t)\leq-2\tilde{\alpha}_1 (\alpha_2+\overline{\alpha}_3)^{-1} V(t).
	\end{equation}
	From (\ref{eq2}),we can obtain: 
	\begin{equation*}
		V(t)\leq V(0)e^{-2\tilde{\alpha}_1 (\alpha_2+\overline{\alpha}_3)^{-1}t}.
	\end{equation*}	
	Thus it can be gotten:
	\begin{equation}\label{eq3}
		\left\|\xi_s(t)\right\|\leq	\alpha_6^{\frac{1}{2}}\left\|\xi_s(0)\right\|e^{-\tilde{\alpha}_1 (\alpha_2+\overline{\alpha}_3)^{-1}t}.
	\end{equation}	
	From (\ref{e38}), if $\chi(t)=0$, one has :
	\begin{equation*}
		\dot{V}(t)\leq \bar{\xi}_{1}^T\Phi_3\bar{\xi}_{1}.
	\end{equation*} 
	If (\ref{eq1}) is satisfied, it can be obtained:
	\begin{equation*}
		\dot{V}(t)\leq-2\tilde{\alpha}_2 \left\| \bar{\xi}_1\right\|^2,
	\end{equation*}
	where $\tilde{\alpha}_2=\frac{1}{2}|\lambda_{\max}(\Phi_3)|$. Similar to the above analysis, the conclusion can obviously be deduced:
	\begin{equation}\label{eq4}
		\left\|\xi_s(t)\right\|\leq	\alpha_6^{\frac{1}{2}}\left\|\xi_s(0)\right\|e^{-\tilde{\alpha}_2 (\alpha_2+\overline{\alpha}_3)^{-1}t}.
	\end{equation}	
	
	In summary, combining (\ref{eq3}) and (\ref{eq4}), the closed-loop system (\ref{e17}) is exponentially stable. The proof is complete.
\end{IEEEproof}

The BMIs problem in Corollary 1 can be solved similar to Theorems 1-2, $N=r_2P_{11}$, then we can obtain:
\begin{equation*}
		K=B_{2,s}^{-1}P_{11}^{-1}\Xi_2.
\end{equation*}

\subsection{$H_\infty$ Performance Control Design}
The $H_\infty$ performance control issue of the slow system is discussed in this subsection. By designing appropriate controller gain, the impact of external disturbances is suppressed tremendously.

\begin{theorem}\label{thm3} Given scalars $\alpha$, $\epsilon$, $\beta_2$, matrix $\Lambda$, and the weights $\boldsymbol{V}^{*},\boldsymbol{W}^{*}$ obtained through MNN. For the closed-loop system (\ref{e17}), if there exist positive-definite matrices $P$, $U>0$, and matrices $Q_1,Q_2$, scalars $\varpi_i>0$, appropriate matrices $M_1$, $M_2$, $M_3$, $N$, and $L>0$, scalar $\beta_1>0$, with suitable switching controller gain matrix $K$, with a prescribed attenuation level $\gamma_s$, s.t. the following inequalities are feasible with LMI (\ref{e31}) (\ref{e41}):
		\begin{equation}\label{e45}
		\bar{\Phi}_{5,1}<0,\ \bar{\Phi}_{5,2}<0,\
		\bar{\Phi}_{6,1}<0,\ \bar{\Phi}_{6,2}<0,
	\end{equation}
	where
	\begin{equation*}
		\begin{aligned}
			\bar{\Phi}_{5,1}&=
			\begin{bmatrix}
				\Upsilon_1&\Upsilon _2\\
				\ast&\Upsilon _3\\ 
			\end{bmatrix},\
			\bar{\Phi}_{6,1}=
			\begin{bmatrix}
				\Upsilon_1&\Upsilon_6\\
				\ast&\Upsilon_3\\
			\end{bmatrix},\\
			\bar{\Phi}_{5,2}&=
			\begin{bmatrix}
				\Upsilon_4&\Upsilon_2&\Upsilon_5\\
				\ast&\Upsilon_3&0\\
				\ast&\ast&\eta_ {2,77}
			\end{bmatrix},\			
			\bar{\Phi}_{6,2}=
			\begin{bmatrix}
				\Upsilon_4&\Upsilon_6&\Upsilon_7\\
				\ast&\Upsilon_3&0\\
				\ast&\ast&\eta_ {2,77}
			\end{bmatrix},\\
			\Upsilon_1&=
			\begin{bmatrix}
				\eta_ {5,11}&\eta_ {1,12}&\eta_ {1,13}\\
				\ast&\eta_ {5,22}&\eta_ {1,23}\\
				\ast&\ast&\eta_ {1,33}\\
			\end{bmatrix},\\
			\Upsilon_2&=
			\begin{bmatrix}
				0&\eta_ {1,15}&P_{11}&\eta_ {5,17}\\
				0&\eta_ {1,25}&N&\eta_ {5,27}\\
				0&0&0&0\\
			\end{bmatrix},\\
			\Upsilon_3&=
			\begin{bmatrix}
				-\Lambda&0&0&0\\
				\ast&-L&0&0\\
				\ast&\ast&-\beta_1&0\\
				\ast&\ast&\ast&-\gamma_s^2\\
			\end{bmatrix},\\
			\Upsilon_4&=
			\begin{bmatrix}
				\eta_ {5,11}&\eta_ {2,12}&\eta_ {1,13}\\
				\ast&-2N&\eta_ {2,23}\\
				\ast&\ast&\eta_ {1,33}\\
			\end{bmatrix},\
			\Upsilon_5=
			\begin{bmatrix}
				\eta_ {2,17}\\ \eta_ {2,27}\\ \eta_ {2,37}
			\end{bmatrix},\\
			\Upsilon_6&=
			\begin{bmatrix}
				\eta_ {3,14}&\eta_ {1,15}&P_{11}&\eta_ {5,17}\\
				\eta_ {3,24}&\eta_ {1,25}&N&\eta_ {5,27}\\
				0&0&0&0\\
			\end{bmatrix},\
			\Upsilon_7=
			\begin{bmatrix}
				\eta_ {4,17}\\ \eta_ {2,27}\\ \eta_ {2,37}
			\end{bmatrix},\\
											\end{aligned}
		\end{equation*}
	\begin{equation*}
	\begin{aligned}
				\eta_ {5,11}&=[P_{11}(A_{s}-B_{2,s}KC_s)+\ast]-(\frac{Q_1}{2}+\ast)\\
			&\ +\beta_1\delta^2- \boldsymbol{\zeta}_1^TL \boldsymbol{\zeta}_2+2M_1+ C_s^{T}(t)(\epsilon\Lambda+1)  C_s(t),\\
			\eta_ {5,17}&=PB_{1,s},\ \eta_ {5,22}=hU-2N,\ 	\eta_{5,27}=NB_{1,s}.\\
		\end{aligned}
	\end{equation*}
	Then the closed-loop system (\ref{e17}) satisfies the $H_\infty$ performance stable with a prescribed attenuation level $\gamma_s$.
\end{theorem}

\begin{IEEEproof} 
	Based on Lyapunov functional, if $\chi(t)=1$, the inequality that follows is true:
	\begin{equation}\label{e47}
		\begin{aligned}
			\dot{V}(t)+ y_s^{T}(t)y_s(t)-\gamma_s^{2}{d(t)}^{T}d(t)&\leq \bar{\xi}_{2}^T\Phi_5\bar{\xi}_{2}
			,\end{aligned}\end{equation}
	where
	\begin{equation*}\begin{aligned}
			\bar{\xi}_{2}&=[\xi_{s}^T(t),\dot{\xi}_{s}^T(t),\xi_{s}^T(t-\tau),e(t),\boldsymbol{\mu}^T(t),\\
			&\ \Delta\boldsymbol{A}^T(\xi_{s}(t)),d(t),\nu_1^T(t)],\\
	\end{aligned}\end{equation*}
	\begin{equation*}\begin{aligned}
			\Phi_{5}&=
			\begin{bmatrix}
				\eta_ {5,11}&\eta_ {2}&\eta_ {1,13}&0&\eta_ {1,15}&P_{11}&\eta_ {5,17}&\eta_ {3}\\
				\ast&\eta_ {4}&\eta_ {5}&0&\eta_ {1,25}&N&\eta_ {5,27}&\eta_ {6}\\
				\ast&\ast&\eta_ {1,33}&0&0&0&0&\eta_ {7}\\
				\ast&\ast&\ast&-\Lambda&0&0&0&0\\
				\ast&\ast&\ast&\ast&-L&0&0&0\\
				\ast&\ast&\ast&\ast&\ast&-\beta_1&0&0\\
				\ast&\ast&\ast&\ast&\ast&\ast&-\gamma_s^2&0\\
				\ast&\ast&\ast&\ast&\ast&\ast&\ast&\eta_ {8}\\
			\end{bmatrix}.\\
	\end{aligned}\end{equation*}
	
	If (\ref{e45}) is satisfied, it is implied that
	\begin{equation*}
		\Phi_{5}<0,
	\end{equation*}
	then we can derive the following inequality by (\ref{e47})
	\begin{equation*}\begin{aligned}
			\dot{V}(t)+ y_s^{T}(t)y_s(t)-\gamma_s^{2}{d(t)}^{T}d(t)&\leq 0.
	\end{aligned}\end{equation*}
	
	By integrating both sides of the above equation simultaneously, we can obtain
	\begin{equation*}
		\int_{0}^{T}(\dot{V}(s)+ y_s^{T}(s)y_s(s)-\gamma_s^{2}{d(s)}^{T}d(s))ds\leq 0,
	\end{equation*}
	then it can be satisfied that
	\begin{equation*}
		\int_{0}^{T} y_s^{T}(s)y_s(s)ds\leq -\int_{0}^{T}\dot{V}(s)ds+	\gamma_s^{2}\int_{0}^{T}{d(s)}^{T}d(s))ds.
	\end{equation*}
	
	Thus, we have
	\begin{equation}\label{e48}
		\int_{0}^{T} y_s^{T}(s)y_s(s)ds\leq 	\gamma_s^{2}\int_{0}^{T}{d(s)}^{T}d(s))ds.
	\end{equation}
	
	Else if $\chi(t)=0$, similar to the above derivation process, we can obtain
	\begin{equation}\label{e49}
		\begin{aligned}
			\dot{V}(t)+ y_s^{T}(t)y_s(t)-\gamma_s^{2}{d(t)}^{T}d(t)&\leq \bar{\xi}_{2}^T\Phi_6\bar{\xi}_{2}
			,\end{aligned}\end{equation}
	where 
	\begin{equation*}\begin{aligned}
			\Phi_{6}&=
			\begin{bmatrix}
				\eta_ {5,11}&\eta_ {2}&\eta_ {1,13}&\eta_ {3,14}&\eta_ {1,15}&P&\eta_ {5,17}&\eta_ {10}\\
				\ast&\eta_ {4}&\eta_ {5}&\eta_ {3,24}&\eta_ {1,25}&N&\eta_ {5,27}&\eta_ {6}\\
				\ast&\ast&\eta_ {1,33}&0&0&0&0&\eta_ {7}\\
				\ast&\ast&\ast&-\Lambda&0&0&0&0\\
				\ast&\ast&\ast&\ast&-L&0&0&0\\
				\ast&\ast&\ast&\ast&\ast&-\beta_1&0&0\\
				\ast&\ast&\ast&\ast&\ast&\ast&-\gamma_s^2&0\\
				\ast&\ast&\ast&\ast&\ast&\ast&\ast&\eta_ {8}\\
			\end{bmatrix}.\\
	\end{aligned}\end{equation*}
	
	If (\ref{e45}) can be satisfied, the following inequality can be obtained
	\begin{equation*}
		\Phi_{6}<0,
	\end{equation*}
	then we can derive similar to be situation $\chi(t)=1$ by (\ref{e49}):
	\begin{equation}\label{e50}
		\int_{0}^{T} y_s^{T}(s)y_s(s)ds\leq 	\gamma_s^{2}\int_{0}^{T}{d(s)}^{T}d(s))ds.
	\end{equation}
	
	To sum up, from (\ref{e48}) (\ref{e50}) and Theorem 1, the closed system (\ref{e17}) can satisfy the $H_\infty$ performance stable with a prescribed attenuation level $\gamma_s$.
	The proof is complete.
\end{IEEEproof}

Converting the BMIs problem (\ref{e45}) in Theorem 3 to a LMI problem similar to the process of Theorem 2:

\begin{equation}\label{e51}
	\begin{aligned}
	\tilde{\Phi}_{5}&=
	\begin{bmatrix}
		\tilde{\Psi}_{5,11}&	\tilde{\Psi}_{5,12}\\
		\ast&	\tilde{\Psi}_{5,22}
	\end{bmatrix}<0,\ 	
	\tilde{\Phi}_{6}=
	\begin{bmatrix}
		\tilde{\Psi}_{6,11}&	\tilde{\Psi}_{6,12}\\
		\ast&	\tilde{\Psi}_{6,22}
	\end{bmatrix}<0,\\
%
	\tilde{\Phi}_{7}&=
	\begin{bmatrix}
		\tilde{\Psi}_{7,11}&	\tilde{\Psi}_{5,12}\\
		\ast&	\tilde{\Psi}_{5,22}
	\end{bmatrix}<0,\ 	
	\tilde{\Phi}_{8}=
	\begin{bmatrix}
		\tilde{\Psi}_{8,11}&	\tilde{\Psi}_{8,12}\\
		\ast&	\tilde{\Psi}_{6,22}
	\end{bmatrix}<0,
	\end{aligned}
\end{equation}
where\begin{equation*}
	\begin{aligned}
		\tilde{\Psi}_{5,11}&=
		\begin{bmatrix}
			\bar{\Theta}_1&\bar{\Theta}_2\\
			\ast&\bar{\Theta}_3\\
		\end{bmatrix},\	\tilde{\Psi}_{5,22}=
		\begin{bmatrix}
			-\gamma_s^2&0\\
			\ast&\bar{\Theta}_4
		\end{bmatrix},\\
		\tilde{\Psi}_{5,12}&=
		\begin{bmatrix}
			\Theta_1& \Theta _2
		\end{bmatrix},\ \ \tilde{\Psi}_{6,12}=
		\begin{bmatrix}
			\Theta_3& \Theta _2
		\end{bmatrix},\\											
		\tilde{\Psi}_{6,11}&=
		\begin{bmatrix}
			\bar{\Theta}_5&\bar{\Theta}_6\\
			\ast&\bar{\Theta}_3\\
		\end{bmatrix},\ \tilde{\Psi}_{6,22}=
		\begin{bmatrix}
			\bar{\Theta}_7&0\\
			\ast&\bar{\Theta}_4\\
		\end{bmatrix},\\
		\tilde{\Psi}_{7,11}&=
		\begin{bmatrix}
			\bar{\Theta}_1&\bar{\Theta}_8\\
			\ast&\bar{\Theta}_3\\
		\end{bmatrix},\ \tilde{\Psi}_{8,11}=
		\begin{bmatrix}
			\bar{\Theta}_5&\bar{\Theta}_8\\
			\ast&\bar{\Theta}_3\\
		\end{bmatrix},\\
		\tilde{\Psi}_{8,12}&=
		\begin{bmatrix}
			\Theta_4& \Theta _2
		\end{bmatrix},\
		\bar{\Theta}_1=
		\begin{bmatrix}
			\bar{\Psi}_{1,11}&\eta_ {1,12}\\
			\ast&\eta_ {5,22}\\
		\end{bmatrix},\\
		\bar{\Theta}_2&=
		\begin{bmatrix}
			\eta_ {1,13}&0&\eta_ {1,15}&P_{11}\\
			\eta_ {1,23}&0&\eta_ {1,25}&N\\
		\end{bmatrix},\\
		\bar{\Theta}_3&=
		\begin{bmatrix}
			\eta_ {1,33}&0&0&0\\
			\ast&-\Lambda&0&0\\
			\ast&\ast&-L&0\\
			\ast&\ast&\ast&-\beta_1\\
		\end{bmatrix},\\
		\bar{\Theta}_4&=	
		\begin{bmatrix}
			-\bar{\beta}_1&0&0&0\\
			\ast&-\bar{\beta}_2&0&0\\
			\ast&\ast&-\bar{\beta}_2&0\\
			\ast&\ast&\ast&-1-\Lambda_1\\
		\end{bmatrix},\\
		\bar{\Theta}_5&=
		\begin{bmatrix}
			\bar{\Psi}_{1,11}&\eta_ {2,12}\\
			\ast&-2N\\
		\end{bmatrix},\ 	\Theta_1=
		\begin{bmatrix}
		\eta_ {5,17}&	\delta I_m\\
		\eta_ {5,27}&0_{m\times m}\\
		0_{4m\times 1}&0_{4m\times m}\\
		\end{bmatrix},\\
		\bar{\Theta}_6&=
		\begin{bmatrix}
			\eta_ {1,13}&0&\eta_ {1,15}&P_{11}\\
			\eta_ {2,23}&0&\eta_ {1,25}&N\\
		\end{bmatrix},\ \bar{\Theta}_7=
		\begin{bmatrix}
		-\gamma_s^2&0\\
		\ast&\eta_{2,77}\\
		\end{bmatrix},\\
		\bar{\Theta}_8&=
		\begin{bmatrix}
			\eta_ {1,13}&\eta_ {3,14}&\eta_ {1,15}&P_{11}\\
			\eta_ {1,23}&\eta_ {3,24}&\eta_ {1,25}&N\\
		\end{bmatrix},\\
		\Theta_2&=
		\begin{bmatrix}
			\zeta_1^T&\zeta_2^T&C_s^T\\
			0_{m\times n_h}&0_{m\times n_h}&0_{m\times 1}\\
			0_{4m\times n_h}&0_{4m\times n_h}&0_{4m\times 1}\\
		\end{bmatrix},\\
		\Theta_3&=
		\begin{bmatrix}
			\eta_ {5,17}&\eta_ {2,17}&\delta I_m\\
			\eta_ {5,27}&\eta_ {2,27}&0_{m\times m}\\
			0_{m\times 1}&\eta_ {2,37}&0_{m\times m}\\
			0_{3m\times 1}&0_{3m\times m}&0_{3m\times m}\\
		\end{bmatrix},\\
		\Theta_4&=
		\begin{bmatrix}
			\eta_ {5,17}&\eta_ {4,17}&\delta I_m\\
			\eta_ {5,27}&\eta_ {2,27}&0_{m\times m}\\
			0_{m\times 1}&\eta_ {2,37}&0_{m\times m}\\
			0_{3m\times 1}&0_{3m\times m}&0_{3m\times m}\\
		\end{bmatrix}.\\
	\end{aligned}
\end{equation*}

\begin{theorem}\label{thm4} If LMI (\ref{e31}) (\ref{e41}) (\ref{e51}) has feasible solutions similar to Algorithm \ref{alg.1}, $N=r_3P_{11}$, then the appropriate controller gain can be solved
	\begin{equation*}
		K=B_{2,s}^{-1}P_{11}^{-1}\Xi_2,
	\end{equation*}
	and the external disturbance $d(t)$ gradually decays with a prescribed attenuation level $\gamma_s$.
\end{theorem}

Let $\rho^{2}=\gamma_s$, we consider the following optimized problem:
\begin{equation}\label{eq53}\begin{aligned}
	\underset{\textbf{u}}{\min}\ \rho\ \text{subject}\ \text{to}\ \text{LMIs}\ (\ref{e31}) (\ref{e41}) (\ref{e51})
\end{aligned}\end{equation}
where $\textbf{u}=\{\rho>0,\ P>0,\ U>0,\ Q,\ K,\ \varpi_i>0,\ M_1,\ M_2,\ M_3,\ N,\ \text{and}\ L>0,\ \beta_1>0\}$.
\begin{algorithm}[H]
	\caption{Optimization Process of attenuation level $\gamma_s$.}\label{alg:alg3}
	\begin{algorithmic}
		\STATE 
		\STATE {\textsc{Step 1}}
		\STATE \hspace{0.5cm} By solving the LMI (\ref{e31}) (\ref{e41}) (\ref{e51}) by Algorithm \ref{alg.1}, we can obtain an initially feasible solution, from which suitable $K$. Set $k = 0$.
		\STATE {\textsc{Step 2}}
		\STATE \hspace{0.5cm} Solve the subsequent LMI optimized problem utilizing the $K$ obtained in Step 1:
		\begin{equation*}\label{eq54}\begin{aligned}
				\underset{\textbf{u}}{\min}\ \rho\ \text{subject}\ \text{to}\ \text{LMIs}\ (\ref{e31}) (\ref{e41}) (\ref{e51})
		\end{aligned}\end{equation*}
		where $\textbf{u}=\{\rho>0,\ P>0,\ U>0,\ Q,\ \varpi_i>0,\ M_1,\ M_2,\ M_3,\ N,\ \text{and}\ L>0,\ \beta_1>0\}$.
		
		Let $\rho_{k}=\rho$. If $k>0$ and $\left| \rho_{k}-\rho_{k-1}\right| < \omega_{\rho}$, where $\omega_{\rho}$ is apredetermined tolerance, go to Step 4; Otherwise, continue.
		\STATE {\textsc{Step 3}}
		\STATE \hspace{0.5cm} Solve LMI optimized issue with $N$ gotten in Step 2:
	\begin{equation*}\label{eq55}\begin{aligned}
		\underset{\textbf{u}}{\min}\ \rho\ \text{subject}\ \text{to}\ \text{LMIs}\ (\ref{e31}) (\ref{e41}) (\ref{e51})
\end{aligned}\end{equation*}
where $\textbf{u}=\{\rho>0,\ P>0,\ U>0,\ Q,\ K,\ \varpi_i>0,\ M_1,\ M_2,\ M_3,\ \text{and}\ L>0,\ \beta_1>0\}$.
		
		Then, set $k=k+1$ and $\rho_k=\rho$. If $\left|\rho_{k}-\rho_{k-1}\right|<\omega_\rho$ , go to Step 4; Otherwise, return to Step 2.\\
		\STATE {\textsc{Step 4}}
		\STATE \hspace{0.5cm}  A sub-optimal solution of (\ref{eq53}) is derived and the optimized level is $\gamma_{opt}=\sqrt{\rho}$ and $\gamma_s>\gamma_{opt}$.\\
		
	\end{algorithmic}
	\label{alg3}
\end{algorithm}

Using the Algorithm \ref{alg3}, the sub-optimal solution of the above optimized problem is solved. Consecutively, the slow system (\ref{e8}) satisfies $H_\infty$ performance stable with a prescribed attenuation level $\gamma_s$.
\subsection{Stability and $H_\infty$ performance control of PDE systems}
In previous subsections, we analyzed the stability and $H_\infty$ performance of closed-loop systems (\ref{e17}) through Theorems 1-4. However, it is evident that these analyses do not directly generalize to the stability and $H_\infty$ performance of closed-loop systems of PDE system (1)-(3). Therefore, we will delve further into the analysis of PDE systems (1)-(3).

Incorporating the controller (\ref{e16}) into the infinite-dimensional ODE system (\ref{e6}),
we derive the following switching closed-loop system:
\begin{equation}\label{e53}
	\left\{\begin{aligned}
		\dot{\xi}_{s}(t)&=(A_{s}-B_{2,s}KC_s)\xi_{s}(t)+	\boldsymbol{W}^{*}\boldsymbol{\mu}(\boldsymbol{\varsigma}(t))+\Delta\boldsymbol{A}(\xi_{s}(t))\\
		&\ +\chi(t) B_{2,s}KC_s\int_{t-\tau(t)}^{t}\dot{\xi}_s(s) ds\\
		&\ -(1-\chi(t)) B_{2,s}Ke(t)+B_{1,s}d(t)+\Delta f_s,\\
		\dot{\xi}_{f}(t)&=A_{f}\xi_{f}(t)+f_{f}(\xi_s(t),\xi_f(t))-B_{2,f}KC_s \xi_s(t)\\
		&\ +\chi(t) B_{2,f}KC_s\int_{t-\tau(t)}^{t}\dot{\xi}_s(s) ds\\
		&\ -(1-\chi(t)) B_{2,f}Ke(t)+B_{1,f}d(t),\\
	\end{aligned}\right.
\end{equation}
where 
\begin{equation*}\begin{aligned}
		\Delta f_s&=f_{s}(\xi_s(t),\xi_f(t)))-f_{s}(\xi_s(t),0),\\
	\chi(t)&=\left\{\begin{aligned}
		&1,\ t\in[t_k,t_k+h],\\
		&0,\ t\in[t_k+h,t_{k+1}].\\
	\end{aligned}\right.
\end{aligned}\end{equation*} 
\begin{remark}
	In previous sections, we solely focused on the main dynamic characteristics of the PDE systems, neglecting the output $y_f$ of the fast subsystem. However, from this subsection onward, we will consider both $y_s(t)$	and $y_f(t)$ simultaneously, namely, the actual output $y(t)$ of the PDE system.
\end{remark}
\begin{theorem}
	Given scalars $\alpha$, $\epsilon$, $\beta_2$, matrix $\Lambda$, and the weights $\boldsymbol{V}^{*},\boldsymbol{W}^{*}$ obtained through MNN. For the closed-loop system (\ref{e17}), if there exist positive-definite matrices $P$, $U$, and matrices $Q_1,Q_2$, scalars $\varpi_i>0$, and appropriate matrices $M_1$, $M_2$, $M_3$, $N$, $L$, scalar $\beta_1$, with suitable switching controller gain matrix $K$, such that the inequalities (\ref{e31}) (\ref{e29}) are feasible. Then, there exists a positive real number $\varepsilon^\ast$ such that if $\varepsilon\in(0,\varepsilon^\ast)$, the controller (\ref{e16}) can ensure that PDE system (1)-(3) is SGUUB and satisfy the $H_\infty$ control performance.
\end{theorem}
\begin{IEEEproof}
	See Appendix \ref{B}.
\end{IEEEproof}
\begin{remark}
	Theorem 5 demonstrates from a different perspective that the slow subsystem separated by the Galerkin approach can relatively accurately capture the main dynamic characteristics of the PDE systems.
\end{remark}
\section{Simulation}\label{s5}
Two simulation examples are used in this section to confirm the efficacy of the previously mentioned theory.
\begin{example}\label{example1}
	The ammonia synthesis reaction takes place in a fixed-bed catalytic reactor, where nitrogen and hydrogen are converted into ammonia through the catalysis of catalytic rods. During this process, the concentrations of reactants and products are crucial for the progression of the reaction. However, achieving stable control over these concentrations is challenging due to the variability of the catalyst surface activity over time and under different conditions. Therefore, consider the following catalytic rod reaction model:
	\begin{equation}\label{e54}
		\left\{\begin{aligned}
			\frac{\partial\bar{\xi}(p,t)}{\partial t}&=\frac{\alpha\frac{\partial\bar{\xi}(p,t)}{\partial p}}{\partial p}+\bar{\beta}_{1}\bar{\xi}(p,t)+\bar{\beta}_{2}\bar{\xi}^2(p,t)\\
			&\ +b_{2}(p)u(t)+b_{1}(p)d(t),\\
			y(t)&=\int_{\Omega}\bar{c}(p)\bar{\xi}(p,t)dp,
		\end{aligned}\right.
	\end{equation}
	with $\bar{\xi}(0,t)=\bar{\xi}(\pi,t)=0$, $\bar{\xi}_{0}(p)=-0.4(\frac{2}{\pi})^{\frac{3}{2}}sin(p)-0.075\sqrt{\frac{2}{\pi}}cos(p)$, where $\bar{\xi}(p,t)$ represents the concentration at spatial position $p$ and time $t$, $\alpha$ is the diffusion coefficient, the latter function terms indicate the catalytic reaction function. $u(t)$ is the input signal, $d(t)$ is the external disturbance which includes temperature change caused by the input. And $\alpha=1$, $\bar{\beta}_{1}=1.65$, $\bar{\beta}_{2}=1.5$,  $\bar{c}(p)=\sqrt{\frac{2}{\pi}}sin(p)+\frac{3}{4}\sqrt{\frac{\pi}{2}}cos(p)$, $b_{1}(p)=\frac{4}{\pi}sin(p)+\frac{3}{8}\pi cos(p)$, $b_{2}(p)=[-\frac{9}{4}\sqrt{\frac{2}{\pi}}cos(p)-\frac{3}{\pi}\sqrt{\frac{2}{\pi}}sin(p),\ -\frac{5}{\pi}\sqrt{\frac{2}{\pi}}sin(p)]$.

	First, Galerkin technology is used to separate the PDE system (\ref{e54}). According to the definition of differential operator $\mathscr{A}$ in the Section \ref{s3}, the following eigenvalues and eigenvectors are capable of being resolved:
	
	$\lambda_{j}=-\alpha j^{2}$, $\phi_{j}(p)=\sqrt{\frac{\pi}{2}}sin(jp)$, $j=1,2,\ldots,\infty.$
	
	Select $m=2$, the slow subsystem can be obtained:
	\begin{equation}\label{e55}
		\left\{\begin{aligned}
			\dot{\xi}_{s}(t)&=
			\begin{bmatrix}
				-\alpha&0\\
				0&-4\alpha
			\end{bmatrix}
			\xi_{s}(t)+
			\begin{bmatrix}
				f_{1}(\xi_{s}(t),0)\\f_{2}(\xi_{s}(t),0)
			\end{bmatrix}\\
			&\ +B_{2,s}u(t)+B_{1,s}d(t),\\
			y(t)&=C_{s}\xi_{s}(t)+\varphi(t),
		\end{aligned}\right.
	\end{equation}
	where
	\begin{equation*}
		\xi_{s}(t)=[\xi_{1}(t),\ \xi_{2}(t)]^T,\ B_{2,s}=
		\begin{bmatrix}
			-1.5&-2.5\\-3&0
		\end{bmatrix},
	\end{equation*}
	\begin{equation*}
		B_{1,s}=[2\sqrt{\frac{2}{\pi}},\ 2\sqrt{\frac{2}{\pi}}]^{T}, \ C_{s}=[1,\ 1]^{T},
	\end{equation*}
	\begin{equation*}
		f_{i}(\xi_{s}(t),0)=\beta_{1}\xi_{i}(t)+\beta_{2}\bar{f}_{i}(\xi_{s}(t)),
	\end{equation*}
	\begin{equation*}
		\begin{aligned}
			\bar{f}_{i}(\xi_{s}(t))&=\int_{0}^{\pi}\phi_{i}(p)(\xi_{1}(t)\phi_{1}(p)+\xi_{2}(t)\phi_{2}(p))^2dp,\ i=1,2.
	\end{aligned}\end{equation*}
	
	Although we have presented a specific model here, accurately modeling unknown nonlinearities in practical situations remains challenging. Therefore, employing NNs to approximate these unknown nonlinearities represents a relatively accurate and practical approach. The specific approximation process is expressed as follows.
	\begin{figure}[!t]
		\centering
		\includegraphics[width=0.4\textwidth]{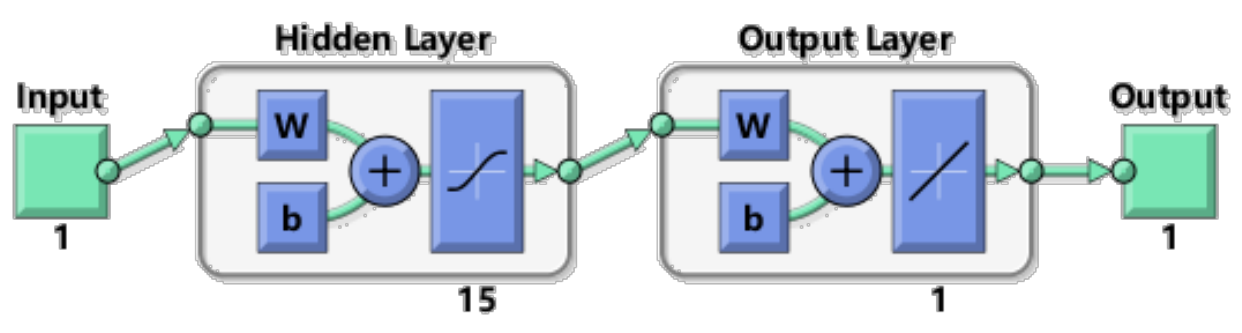}
		\caption{Three-layer neural network architecture (with a single hidden layer).}
		\label{Fig. 3}
	\end{figure}\\
	
	According to Algorithm \ref{alg1}, the MNN is set as Fig. \ref{Fig. 3}. We have established a three-layer NN with a hidden layer containing 15 neurons and an activation function of $\mu_i(\varsigma_i)$, select $q_i=1$, $r_i=1$. Based on the analysis, select $\Delta T_{s}=0.001$, and input $\xi_{s}(t)\in[0,2]$ at intervals of 0.001, the objective function can be approximately obtained:
	\begin{equation*}
		\begin{aligned}
			f_s(\xi_s(t),0)
			&\approx\frac{\xi_{s}(t+\Delta T_{s})-\xi_{s}(t)}{\Delta T_{s}}-A_{s}\xi_{s}(t),s=1,2.
		\end{aligned}
	\end{equation*}
	
	By inputting a set of corresponding data for training by two kinds of algorithms, the results are as Figs. \ref{Fig. 4}-\ref{Fig. 5}.
	\begin{figure}[!t]
		\centering
		\subfloat[]{\includegraphics[width=2.5in]{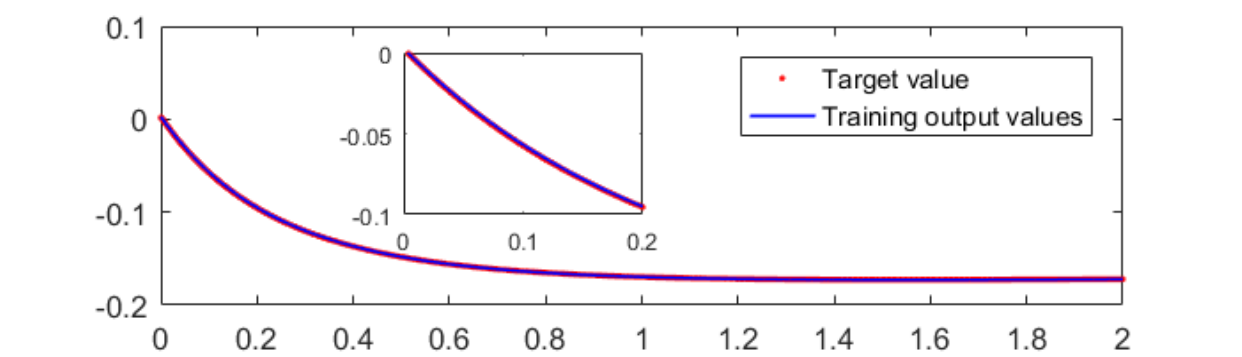}%
			\label{fig4_first_case}}
		\\
		\subfloat[]{\includegraphics[width=2.5in]{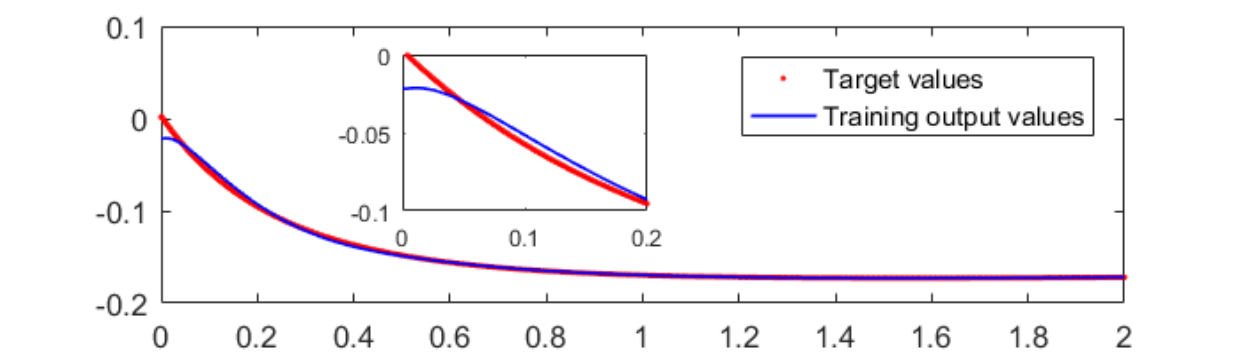}%
			\label{fig4_second_case}}
		\caption{Responses of $f_{nn,1}(\xi_s(t),\boldsymbol{V},\boldsymbol{W})$ and $f_1({\xi}_s(t))$ under different algorithms. (a) LM algorithm. (b) BP algorithm.}
		\label{Fig. 4}
	\end{figure}
	
	\begin{figure}[!t]
		\centering
		\subfloat[]{\includegraphics[width=2.5in]{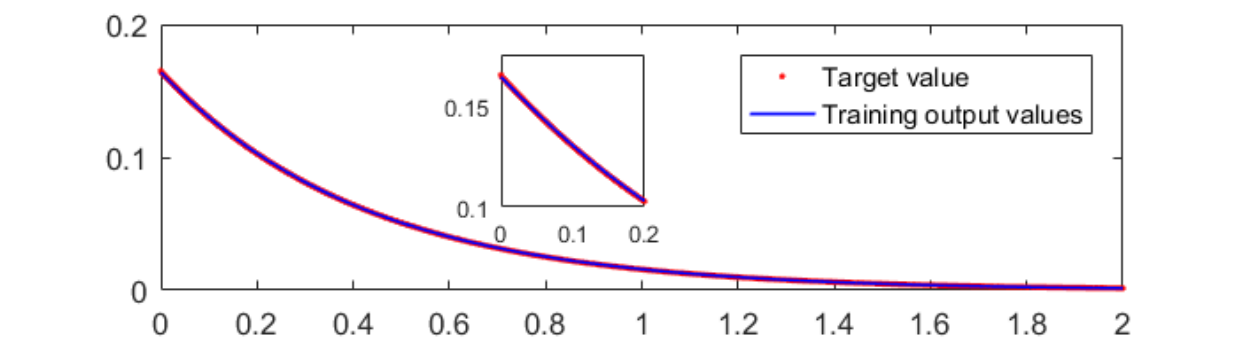}%
			\label{fig5_first_case}}
		\\
		\subfloat[]{\includegraphics[width=2.5in]{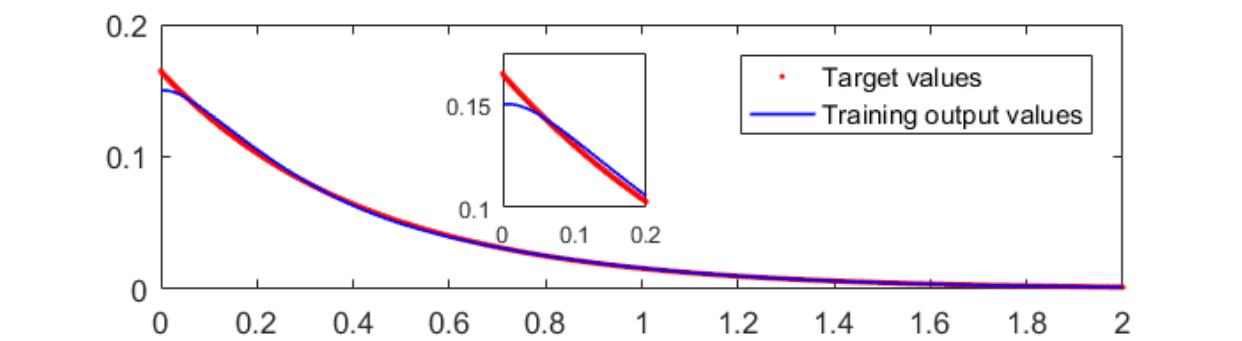}%
			\label{fig5_second_case}}
		\caption{Responses of $f_{nn,2}(\xi_s(t),\boldsymbol{V},\boldsymbol{W})$ and $f_2({\xi}_s(t))$ under different algorithms. (a) LM algorithm. (b) BP algorithm.}
		\label{Fig. 5}
	\end{figure}
	The approximation results obtained using the LM algorithm are depicted in Fig. \ref{fig4_first_case} and Fig. \ref{fig5_first_case}. It can be observed that after 7 and 6 epochs, respectively, the fitting effect of $f_{nn,i}(\xi_s(t),\boldsymbol{V},\boldsymbol{W})$ on $f_i({\xi}_s(t))(i=1,2)$ has reached a good level. Traditional gradient descent techniques, exemplified by the standard BP algorithm, often grapple with sluggish convergence rates. Fig. \ref{fig4_second_case} and Fig. \ref{fig5_second_case} show the outcomes of training with the BP algorithm at a fixed learning rate after 3666 and 2666 epochs, respectively. Conversely, it is clear that the LM algorithm used in this work exhibits much faster training dynamics and a higher convergence efficiency. Fig. \ref{Fig. 6} presents a comparative visualization of the target outputs (actual values) and the predicted outputs generated by the neural network under LM algorithm, superimposed on the same graph. It is strikingly apparent that all data points are densely concentrated around the $y=x$ line, signifying a high degree of congruence between the predictions of NN and the actual values. And the optimized Weights can be gotten as TABLE \ref{tab:1}.
	
	%
	\begin{table*}[!t]
		\caption{Optimized Wights $\boldsymbol{W}^{*}$, $\boldsymbol{V}^{*}$ (Only three decimal places are displayed here). \label{tab:1}}
		\centering
		\begin{tabular}{|c|}
			\hline
			$\boldsymbol{W}^{*}=\begin{pmatrix}
			0.905,	-0.019,	-0.126,	-0.039,	0.502,	0.094,	0.601,	-0.923,	0.715,	0.079,	-0.371,	-0.354,	-0.796,	-0.031	,0.406\\
			0.025,	-0.361,	0.680,	-0.850,	0.327,	-0.015,	0.056,	-0.041,0.545,	0.710,	0.011,	-0.649,	0.282,	-0.996,0.053
			\end{pmatrix}$\\
			\hline
			$\boldsymbol{V}^{*T}=\begin{pmatrix}
			0.998,-0.275,-0.325,0.409,-0.959,0.984,0.417,0.449,-0.274,0.920,0.058,0.211,0.271,0.507,0.832\\
			-0.058,0.085,0.543,-0.974,-0.413,-0.779,0.146,-0.855,-0.442,-0.055,0.206,0.869,-0.568,-0.571,0.547
			
			\end{pmatrix}$\\
			\hline
		\end{tabular}
	\end{table*}
	%
	
	\begin{figure}[!t]
		\centering
		\subfloat[]{\includegraphics[width=2.3in]{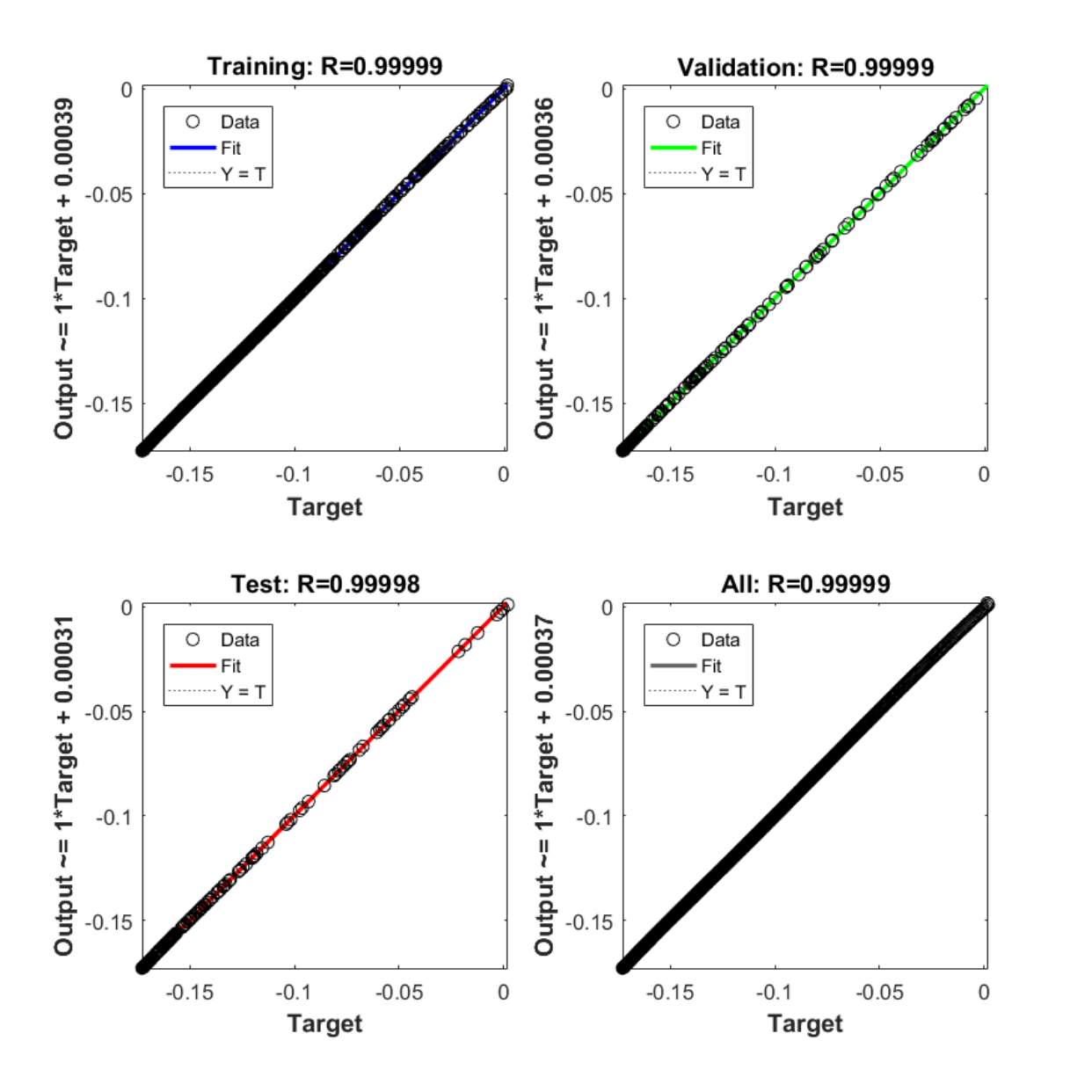}%
			\label{fig6_first_case}}
		\\
		\subfloat[]{\includegraphics[width=2.3in]{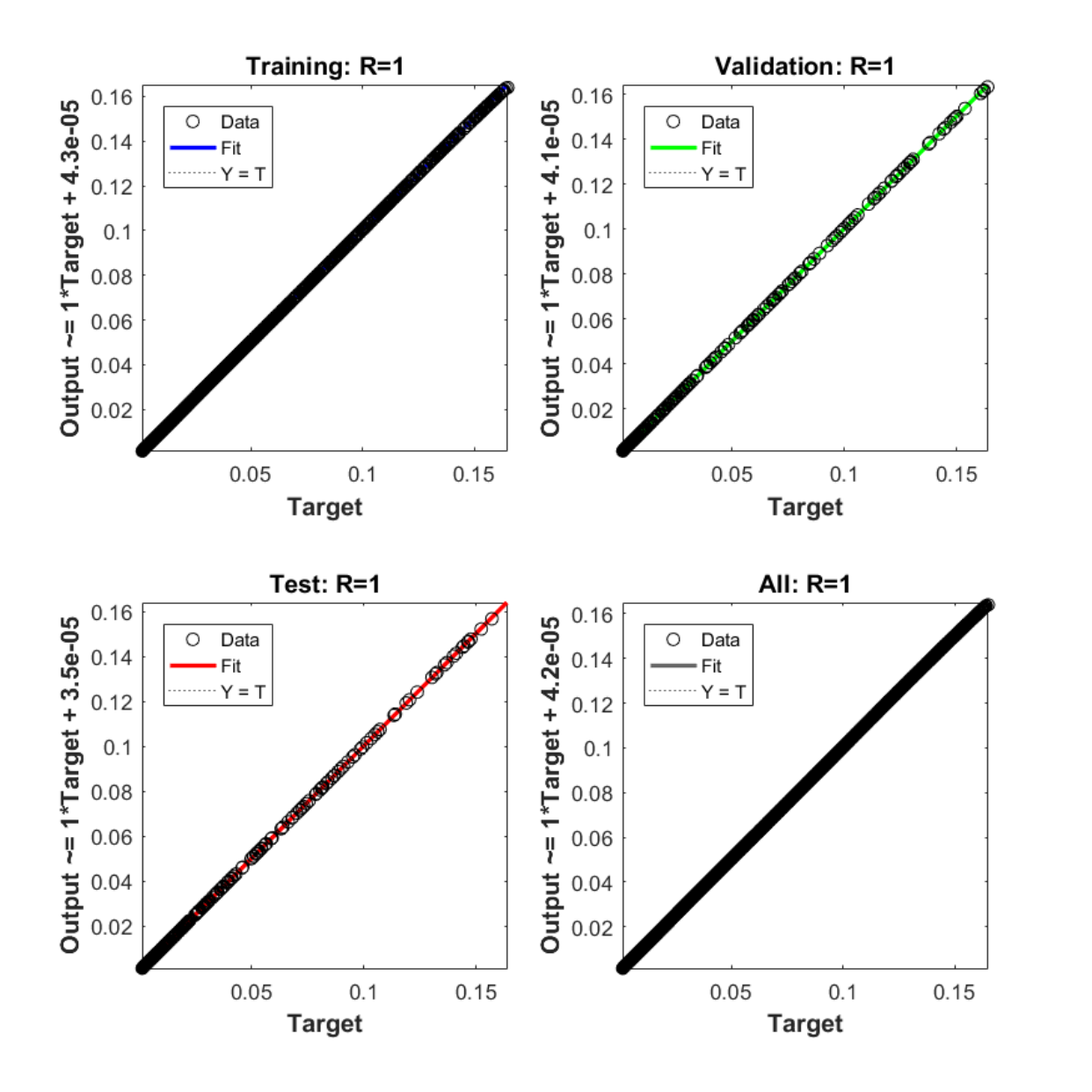}%
			\label{fig6_second_case}}
		\caption{Regression of $f_{nn,i}(\xi_s(t),\boldsymbol{V},\boldsymbol{W})$ and $f_i(\xi_s(t))$ under LM algorithm. (a) $i=1$. (b) $i=2$.}
		\label{Fig. 6}
	\end{figure}
	%
	
	According to (\ref{e12}) (\ref{e13}), we can obtain $\delta=0.0509$, $G_{\min}=0$, $G_{\max}=0.5$. Select $\xi_1(0)=-0.4$, $\xi_2(0)=0.1$, $\beta_2=1.110$, and $\epsilon=0.01$. Initially, we undertake a comprehensive examination of stability, where $h= 0.11$, $\beta_1=2.9061\times 10^3$, $\Lambda= 3.4901\times 10^3$. 
	
	
	Through solving LMIs (\ref{e31}) (\ref{e41}) by Algorithm \ref{alg.1}, the event-triggered output feedback controller gain $K$ can be obtained: 
	\begin{equation*}
		r_1 =\begin{bmatrix}
			8.5994 &  0.0856\\
			-0.6496 &   0.0579\\
		\end{bmatrix},\	K=\begin{bmatrix}
			-2.3002\\
			-0.5612\\
		\end{bmatrix}.
	\end{equation*}
	Fig. \ref{Fig. 7} shows that open and closed-loop evolutionary process of $\bar{\xi}_s(p,t)$. It can be observed that PDE system (\ref{e54}) become SGUUB under switching event triggered control. Fig. \ref{Fig. 8} shows the trajectory of state $\xi_s(t)$ of slow system. 
	\begin{figure}[!t]
		\centering
		\subfloat[]{\includegraphics[width=2in]{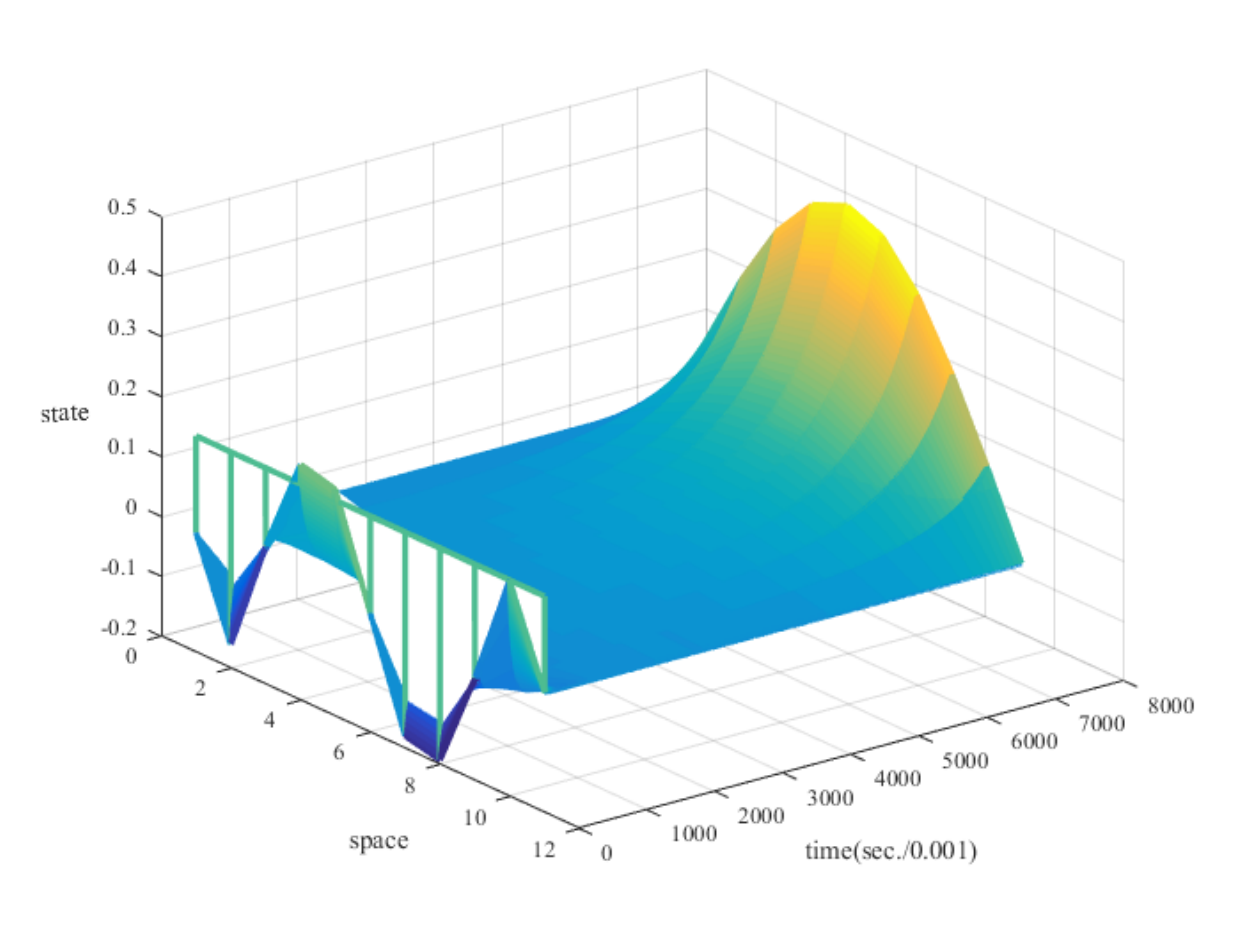}%
			\label{fig7_first_case}}
		\\
		\subfloat[]{\includegraphics[width=2in]{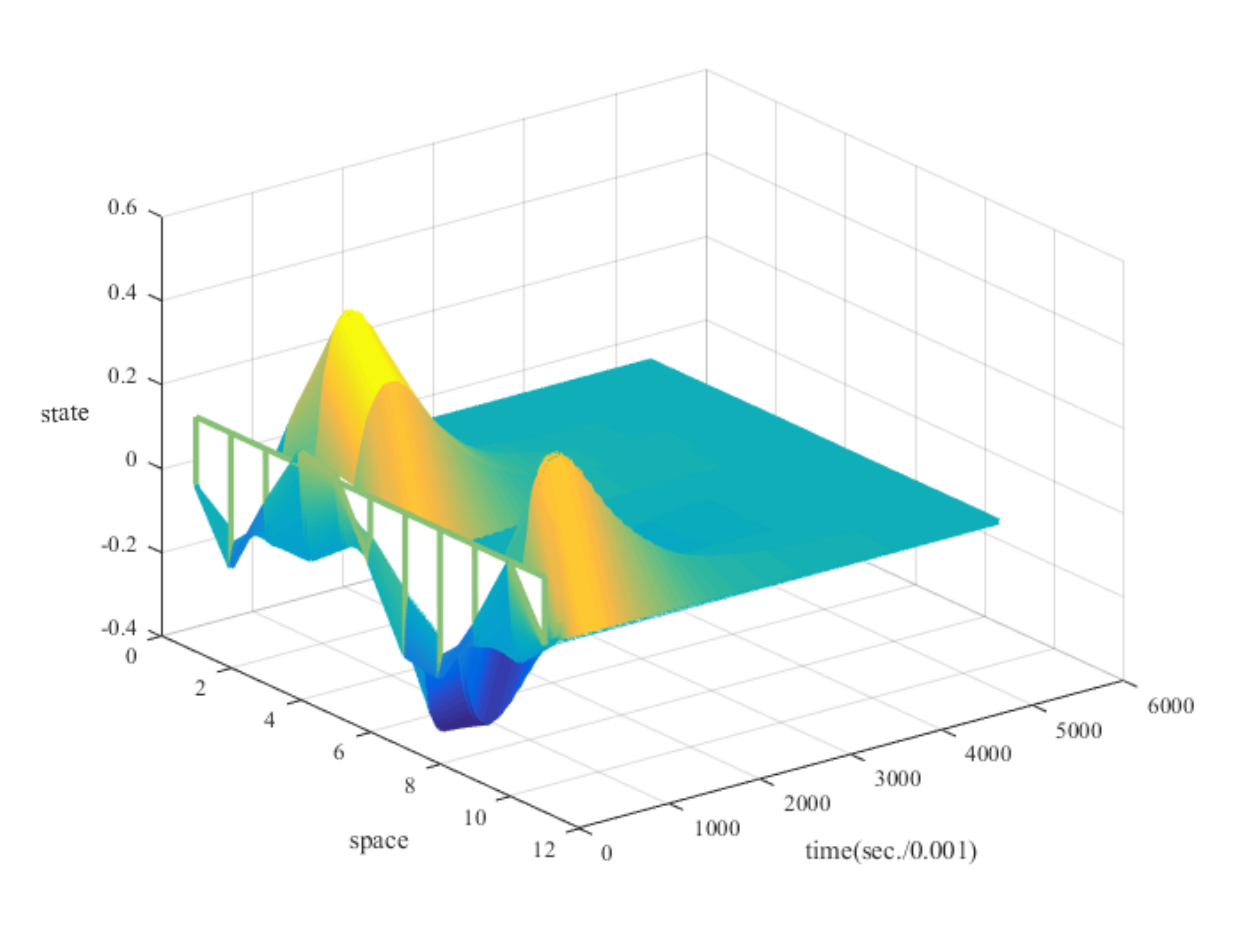}%
			\label{fig7_second_case}}
		\caption{Evolutionary process of $\bar{\xi}_s(p,t)$. (a) Open-loop. (b) Closed-loop.}
		\label{Fig. 7}
	\end{figure}
	
	\begin{figure}[!t]
		\centering
		\includegraphics[width=0.35\textwidth]{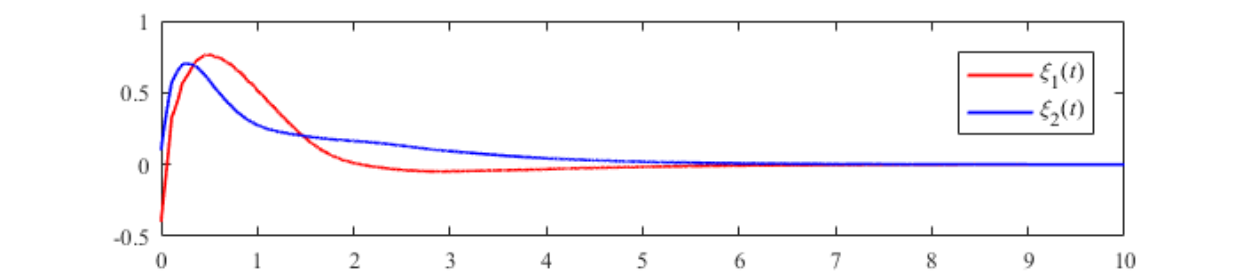}
		\caption{Trajectory of the state $\xi_s(t)$ of slow system.}
		\label{Fig. 8}
	\end{figure}
	
	\begin{figure}[!t]
		\centering
		\subfloat[]{\includegraphics[width=2.5in]{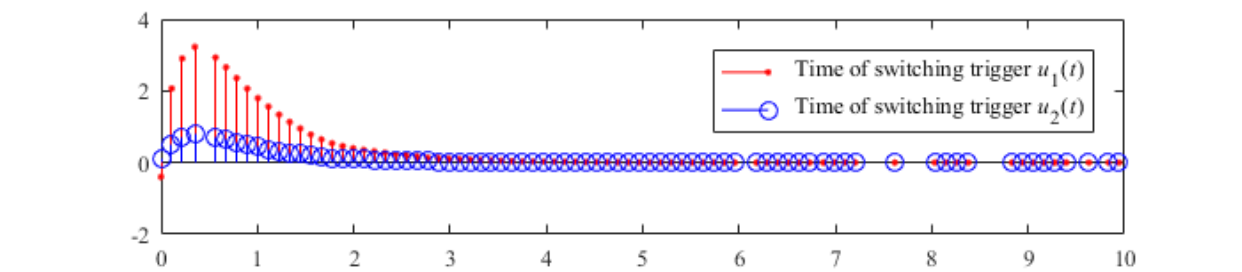}%
			\label{fig9_first_case}}
		\\
		\subfloat[]{\includegraphics[width=2.5in]{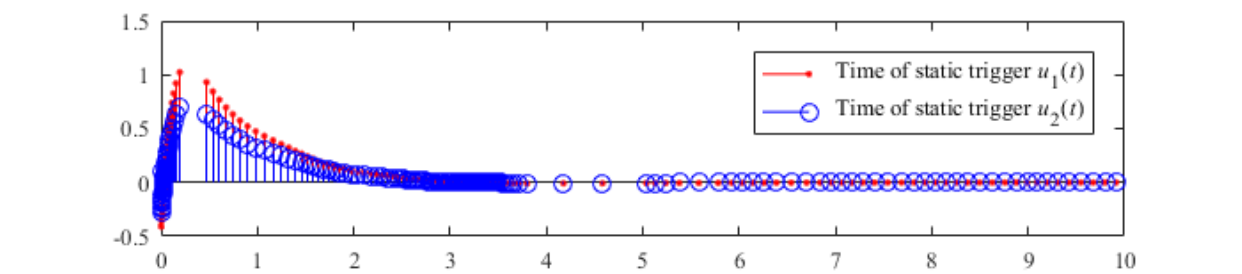}%
			\label{fig9_second_case}}
		\caption{Moment of event-triggered. (a) Switching event-triggered control. (b) Static event-triggered control.}
		\label{Fig. 9}
	\end{figure}
	
	In addition to presenting the system state trajectory variation, we have also facilitated a comparative analysis through trigger instant diagrams across diverse event-triggered control scenarios. The triggered moments for static and switching event-triggered control are shown in Fig. \ref{Fig. 9}, respectively. Upon comparison, it becomes evident that the switching event-triggered control exhibits a significantly lower number of triggered events compared to the static event-triggered control, showcasing a notable efficacy in conserving control resources. 
	\begin{figure}[!t]
		\centering
		\includegraphics[width=0.35\textwidth]{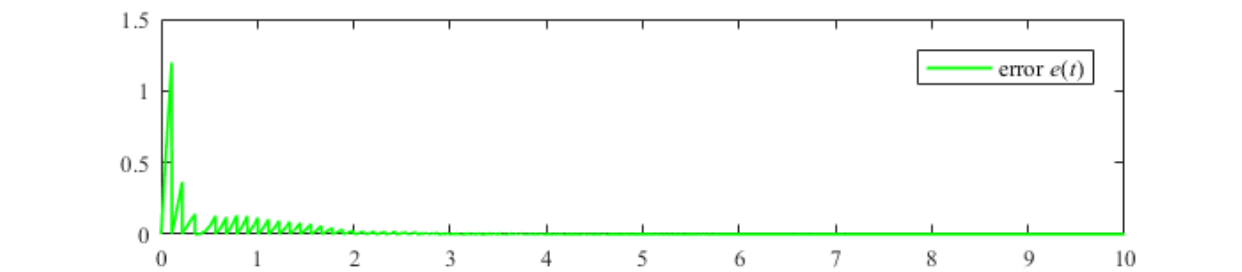}
		\caption{Error of switching triggered event control.}
		\label{Fig. 10}
	\end{figure}
	The error $e(t)$ is displayed in Fig. \ref{Fig. 10}. As the event-triggered control progresses, the error gradually converges to zero. We also show the $u(t)$ in Fig. \ref{Fig. 11}. Through a series of analyses, the favorable effectiveness of the switching event-triggered control has been demonstrated. And the output $y_s(t)$ of slow system can be seen in Fig. \ref{Fig. 12}. The number of event triggers for both static and switching event-triggered control at various sample intervals is then recorded in TABLE \ref{tab:table2}. It has been verified by thorough comparison analysis that there are substantially fewer event triggers in switching event-triggered control than in static event-triggered control.
	\begin{figure}[!t]
		\centering
		\includegraphics[width=0.35\textwidth]{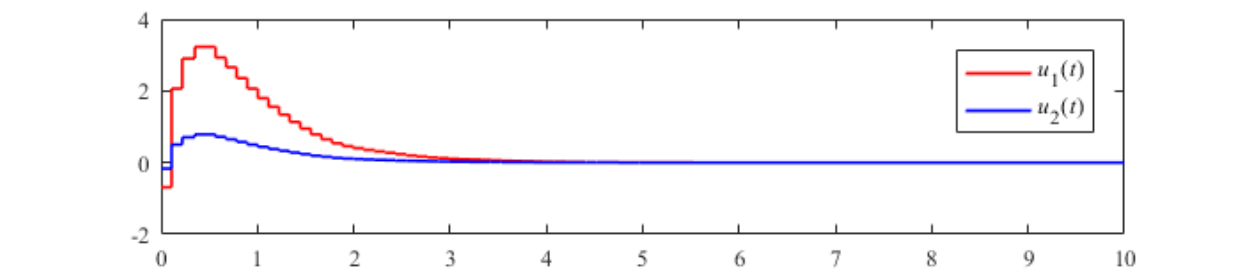}
		\caption{Trajectory of the switching event-triggered controller $u(t)$.}
		\label{Fig. 11}
	\end{figure}
	\begin{figure}[!t]
		\centering
		\includegraphics[width=0.35\textwidth]{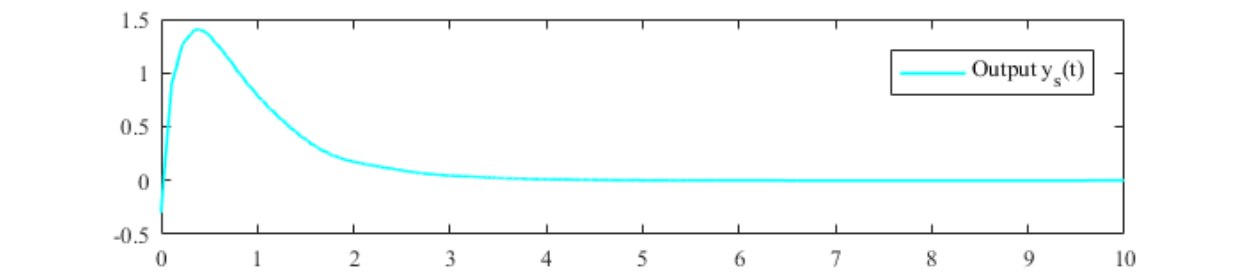}
		\caption{Trajectory of the output $y_s(t)$ of slow system.}
		\label{Fig. 12}
	\end{figure}

	\begin{table}[!t]
		\caption{Triggered Times of Event Triggered Control \label{tab:table2}}
		\centering
		\begin{tabular}{|c||c||c|}
			\hline
			Triggered times &  Switching event triggered &Static event triggered \\
			\hline
			h=0.055& 127&216\\
			\hline
			h=0.110 & 77&224\\
			\hline
		\end{tabular}
	\end{table}
	
		\begin{figure}[!t]
		\centering
		\subfloat[]{\includegraphics[width=2in]{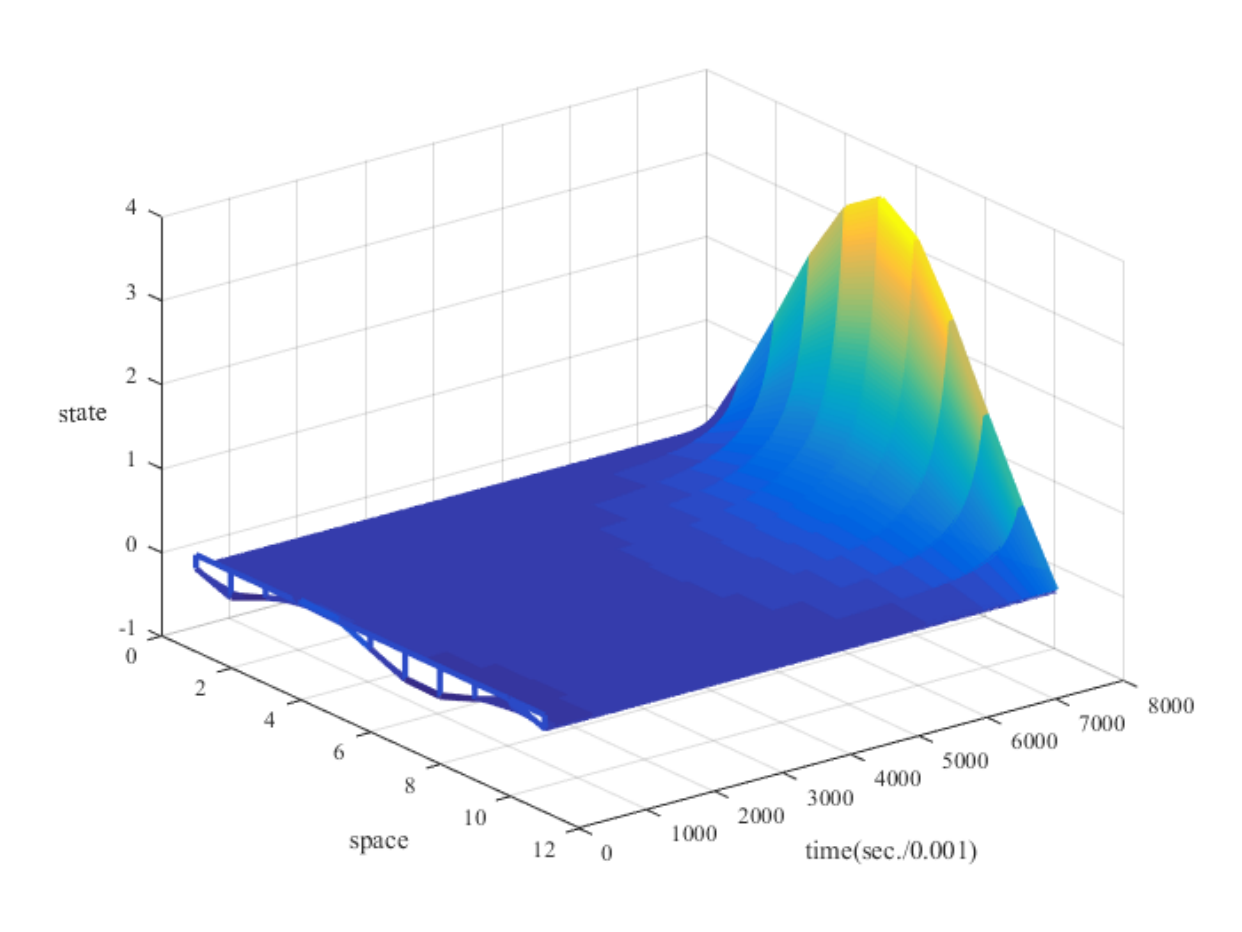}%
			\label{fig7_first_case}}
		\\
		\subfloat[]{\includegraphics[width=2in]{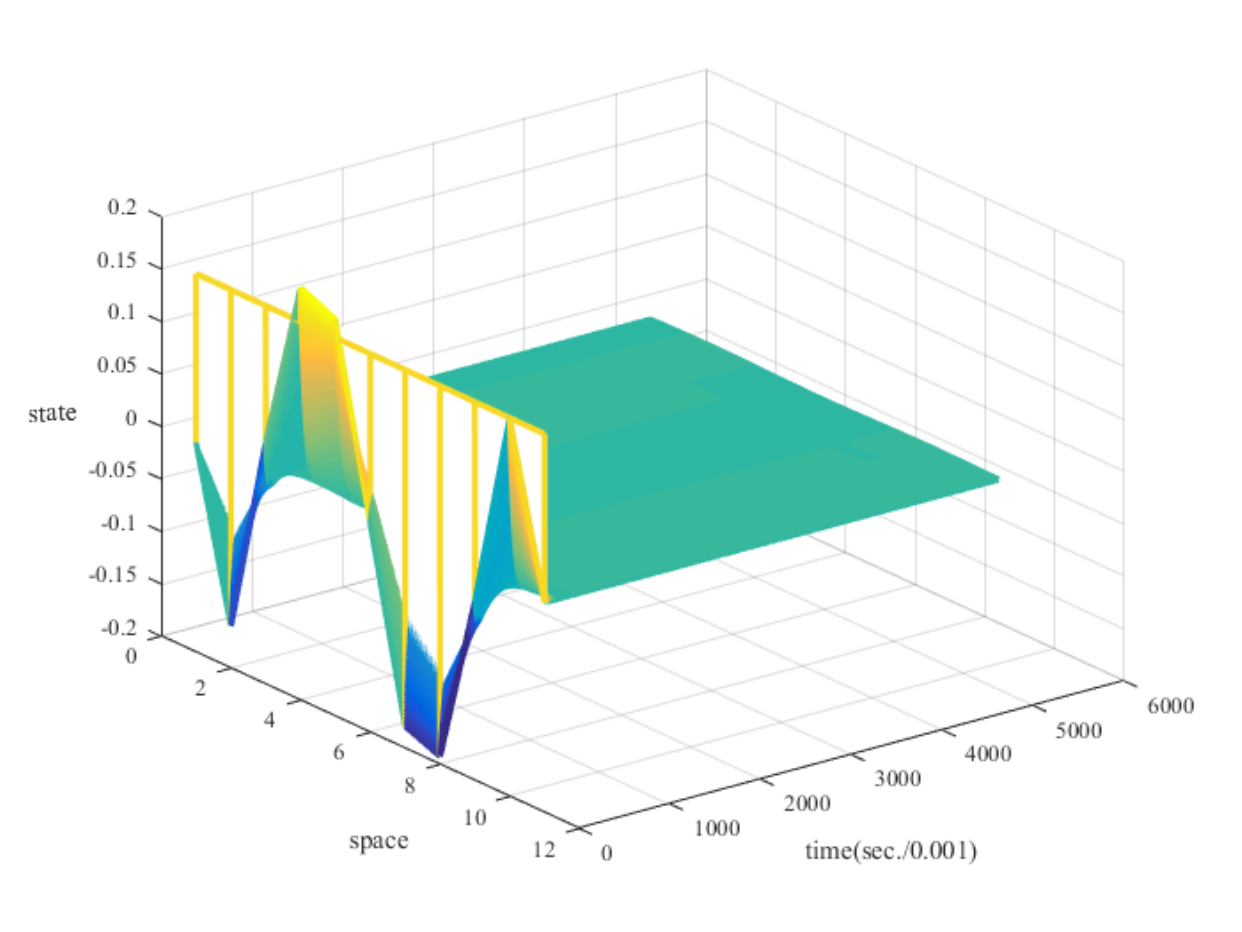}%
			\label{fig7_second_case}}
		\caption{Evolutionary process of $\bar{\xi}_s(p,t)$  when $d(t)\equiv0$. (a) Open-loop. (b) Closed-loop.}
		\label{Fig. 13}
	\end{figure}
	
Subsequently, the BMIs (\ref{e31}) (\ref{eq1}) can  be converted to LMIs similar to Theorem 2 and the controller gain can be solved: \begin{equation*}
		r_2=\begin{bmatrix}
		2.6594  & -0.0576\\
		1.2996  &  0.2289\\
		\end{bmatrix},\  K=\begin{bmatrix}
	\ \ 0.3514\\-0.8798
		\end{bmatrix}.
	\end{equation*}
	
	According to Corollary 1, the effectiveness of the event-triggered switching control is demonstrated in Fig. \ref{Fig. 13} when $d(t)\equiv0$. Under the anticipated monitoring, it is evident that the system finally achieves exponential stability. 
	
			In conclusion, the optimized problem (\ref{eq53}) can be solved by Algorithm \ref{alg3}. It can be gotten 
				\begin{equation*}
						r_3=  \begin{bmatrix}
							3.6594  & -0.0876\\
							-1.2996  &  0.8289\end{bmatrix},\ 
						K=\begin{bmatrix}
								-2.6006\\
								-0.6616\\
							\end{bmatrix}.
					\end{equation*}
	And $\rho=0.7290$. The controller makes ensure that the system satisfies $H_\infty$ performance sable with a sub-optimal attenuation level $\gamma_s=0.5315$. It is showed as Fig. 14, the energy decay degree of PDE system under $H_\infty$ performance control is obviously lower than the previous control method through comparing with trajectories of $\left\|\bar{\xi}(\cdot,t)\right\|_{2,\Omega}$. In summary, Example \ref{example1} effectively demonstrates the superior control performance of switching event-triggered control.
	
			\begin{figure}[!t]
		\centering
		\includegraphics[width=0.35\textwidth]{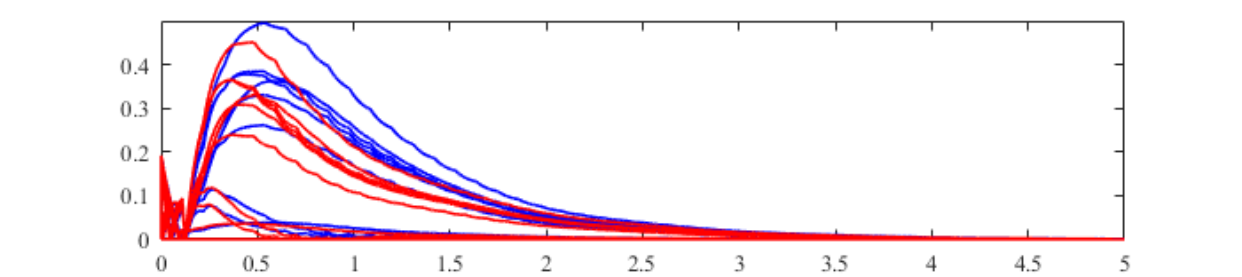}
		\caption{Trajectories of $\left\|\bar{\xi}(\cdot,t)\right\|_{2,\Omega}$ under switching event-triggered control (blue lines) and $H_\infty$ control (red lines) of Example 1 .}
		\label{Fig.}
	\end{figure}
\end{example}

\begin{example}\label{example2}
	Intelligent transportation systems have become a focus of transportation-related research. In addition to wasting people's time, traffic congestion increases energy use and pollutes the environment. According to the research presented in \cite{2009Bivariate,Fokker2010}, from a macroscopic standpoint, the Fokker-Planck equation can be employed to model the probability distribution of traffic flow. In order to determine the effectiveness of the approach presented in this work, we have modeled traffic flow on a specific one-way road stretch in Example \ref{example2} applying the following diffusion equation:
	\begin{equation}\label{e56}
		\left\{\begin{aligned}
			\frac{\partial\bar{\xi}(p,t)}{\partial t}&={\iota_1 \frac{\partial^2\bar{\xi}(p,t)}{\partial p^2}}+\iota_2\bar{\xi}(p,t)+\iota_3\bar{\xi}^2(p,t)\\
			&\ +b_{2}(p)u(t)+b_{1}(p)d(t),\\
			y(t)&=\int_{\Omega}\bar{c}\bar{\xi}(p,t)dp,
		\end{aligned}\right.
	\end{equation}
	with $\bar{\xi}(0,t)=\bar{\xi}(100,t)=0$, $\bar{\xi}_{0}(p)=0.1-0.1cos(p)$. And the specific data is as TABLE \ref{tab:table3} and Example \ref{example1}. In order to replicate more severe phenomena of traffic flow density diffusion, we choose a bigger diffusion coefficient directly and simplify the problem by ignoring the tiny interactions involved in diffusion.

	\begin{table*}[!t]
		\caption{Traffic flow data. \label{tab:table3}}
		\centering
		\begin{tabular}{|c||c||c|}
			\hline
			Parameter &  Interpretation &Value \\
			\hline
			$\bar{\xi}(p,t)$& Traffic density.&$\ast$ (vehicles/m)\\
			\hline
			$p$ & Road space.&$[0,100]$(m)\\
			\hline
			$t$ & Time.&$\ast$ (s)\\
			\hline
			$\iota_1$ &The diffusion coefficient of the unevenness of traffic density spreading on the road. &0.5\\
			\hline
			$\iota_2$&Coefficient of traffic density growth. &0.1\\
			\hline
			$\iota_3$& Coefficient of nonlinear interaction between vehicles (A certain mutual exclusion effect between vehicles).&-0.01 \\
			\hline
			$u(t)$ & Traffic lights or other command methods (here, the application refers to switching event-triggered control methods).&$\ast$\\
			\hline
			$d(t)$ & Interference factors such as road conditions, accidents, ambulances, weather, and other unexpected situations.&$\ast$\\
			\hline
		\end{tabular}
	\end{table*}			
	
	First, Galerkin technology is used to separate the PDE system (49). According to the definition of differential operator $\mathscr{A}$ in the Section \ref{s3}, the following eigenvalues and eigenvectors are capable of being resolved:
	
	$\lambda_{j}=-\iota_1 j^{2}$, $\phi_{j}(p)=\sqrt{\frac{\pi}{2}}sin(jp)$, $j=1,2,\ldots,\infty.$
	
	The slow subsystem can be obtained:
	\begin{equation}\label{e57}
		\left\{\begin{aligned}
			\dot{\xi}_{s}(t)&=
			\begin{bmatrix}
				-\iota_1&0\\
				0&-4\iota_1
			\end{bmatrix}
			\xi_{s}(t)+
			\begin{bmatrix}
				f_{1}(\xi_{s}(t),0)\\f_{2}(\xi_{s}(t),0)
			\end{bmatrix}\\
			&\ +B_{2,s}u(t)+B_{1,s}d(t),\\
			y(t)&=C_{s}\xi_{s}(t)+\varphi(t),
		\end{aligned}\right.
	\end{equation}
	where
	\begin{equation*}\begin{aligned}
			\xi_s(0)&=[0.2\sqrt{2/\pi},\ 2/15\sqrt{2/\pi}],\\
		f_{i}(\xi_{s}(t),0)&=\iota_2\xi_{i}(t)+\iota_3\bar{f}_{i}(\xi_{s}(t)).
		\end{aligned}
	\end{equation*}
	Firstly, the unknown nonlinear function can be modeled by MNN as TABLE \ref{tab:4}. 	
	\begin{table*}[!t]
		\caption{Optimized Wights $\boldsymbol{W}^{*}$, $\boldsymbol{V}^{*}$ (Only three decimal places are displayed here). \label{tab:4}}
		\centering
		\begin{tabular}{|c|}
			\hline
			$\boldsymbol{W}^{*}=\begin{pmatrix}
			-0.362,	-0.200,	0.755,	0.906,	-0.991,0.2030,	0.571,	0.277,	-0.766,	0.421,	-0.805,	0.170,	0.695,	-0.539,	-0.608\\
			-0.346,	0.397,	-0.369,	0.795,	0.221,	-0.589,	0.895,	0.151,	-0.717,	-0.346,	-0.295,	0.143,	0.388,	-0.083	,-0.409
			\end{pmatrix}$
			\\
			\hline
			$\boldsymbol{V}^{*T}=\begin{pmatrix}
				0.314,-0.713,-0.483,0.583,0.455,0.646,0.873,-0.618,0.054,-0.439,-0.183,-0.520,-0.285,-0.671,-0.995\\
				0.459,-0.335,-0.087,-0.049,-0.046,-0.245,0.999,0.299,-0.880,-0.639,-0.235,0.012,0.792,-0.431,0.158
			\end{pmatrix}$\\
			\hline
		\end{tabular}
	\end{table*} 
	Given $h=0.11$, $\delta=0.0114$, $\beta_1= 257.5950$, $\beta_2=1$, $\epsilon=0.05$, $\Lambda= 251.2520$. Then the control gain can be obtained by sloving LMIs (\ref{e31}) (\ref{e41}):
	\begin{equation*}
		r_1=\begin{bmatrix}
		0.7874  & -0.1019\\
		-0.2519  &  0.4408\\
		\end{bmatrix},\ K=\begin{bmatrix}
			  -1.1126\\
		\ \	0.4666
		\end{bmatrix}.
	\end{equation*}
	
	Without any control mechanisms in place, the traffic density shows a rising tendency, as shown in Fig. \ref{fig14_first_case}, which eventually results in traffic congestion. However, Fig. \ref{fig14_second_case} demonstrates an outstanding control effect, where the use of event-triggered control not only provides complete road clearance but also lowers traffic density, demonstrating a superb outcome of the control techniques.
	\begin{figure}[!t]
		\centering
		\subfloat[]{\includegraphics[width=2in]{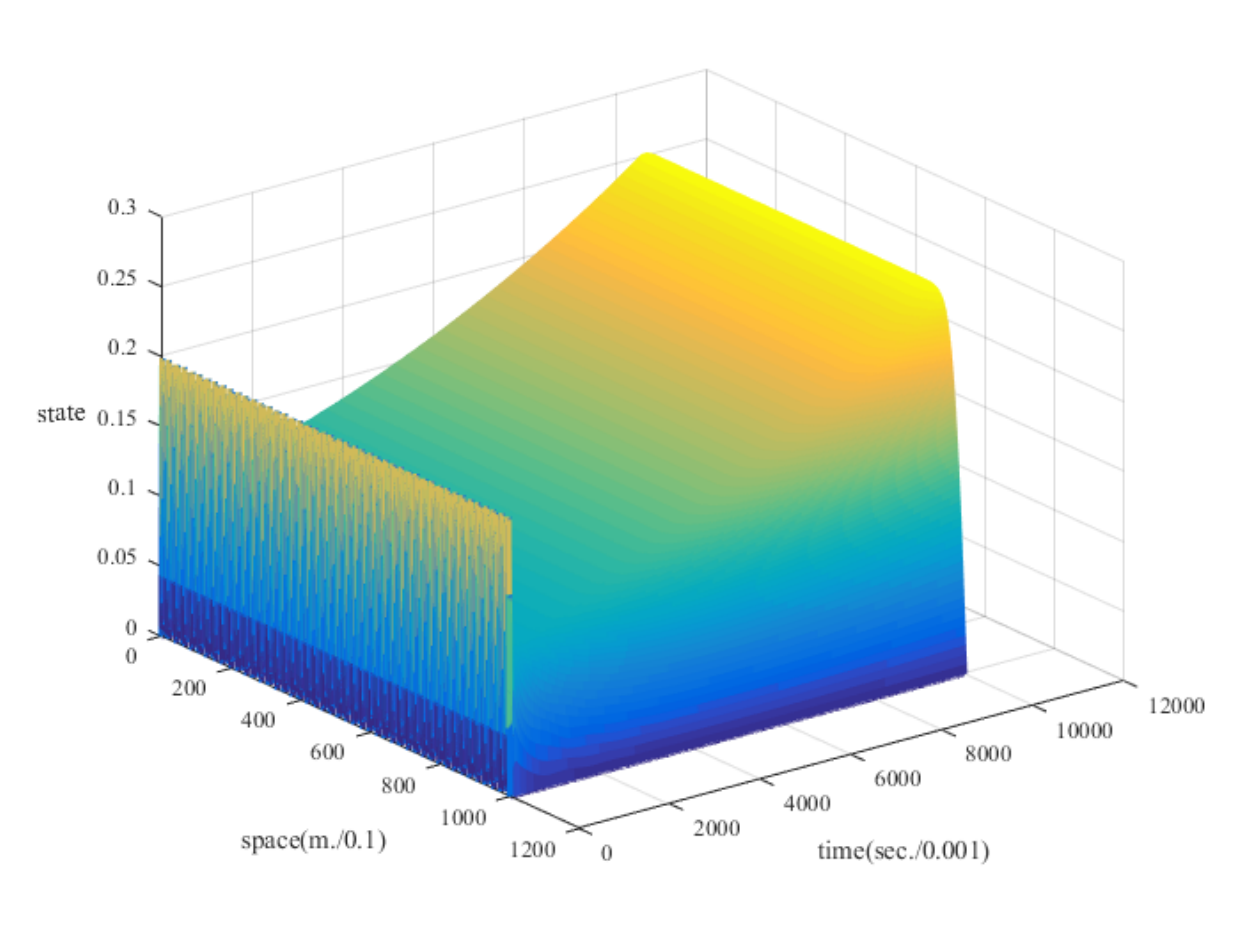}%
			\label{fig14_first_case}}
		\\
		\subfloat[]{\includegraphics[width=2in]{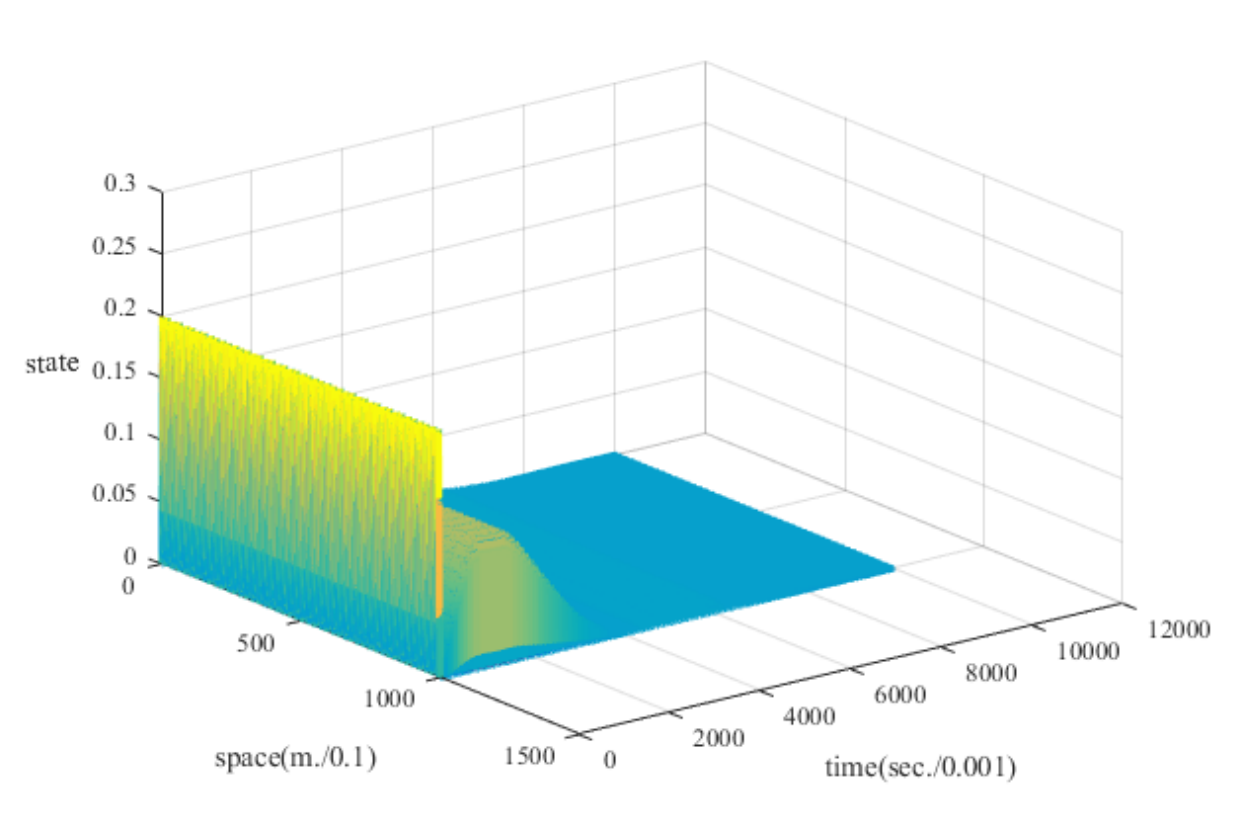}%
			\label{fig14_second_case}}
		\caption{evolutionary process of $\bar{\xi}_s(p,t)$. (a) Open-loop. (b) Closed-loop.}
		\label{Fig. 14}
	\end{figure} 
	The trajectory of state $\xi_s(t)$ of slow system is showed in Fig. \ref{Fig. 15}. Furthermore, Fig. \ref{Fig. 16} illustrates the triggered control moment. And Figs. \ref{Fig. 17}-\ref{Fig. 18} respectively illustrate the trajectory of the controller and error variation. The aforementioned analysis shows that switching event-triggered control in traffic flow is effective. Therefore, Example \ref{example2} demonstrates the control effect of the method proposed in this paper in a simple traffic flow modeling scenario.
	
	\begin{figure}[!t]
		\centering
		\includegraphics[width=0.35\textwidth]{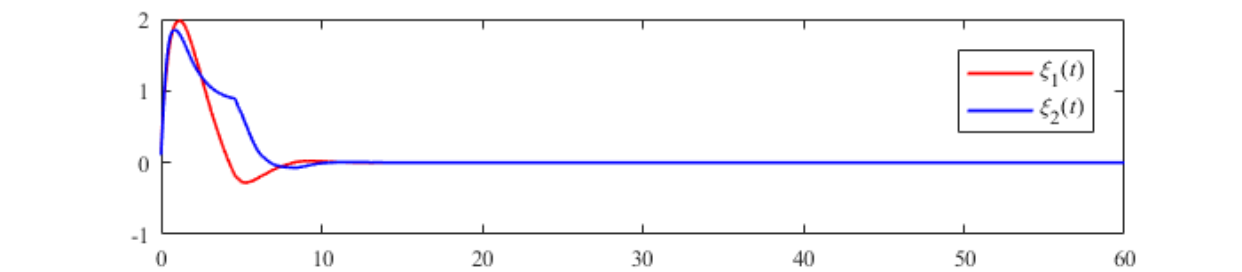}
		\caption{Trajectory of the state $\xi_s(t)$ of slow system.}
		\label{Fig. 15}
	\end{figure}
	\begin{figure}[!t]
		\centering
		\includegraphics[width=0.35\textwidth]{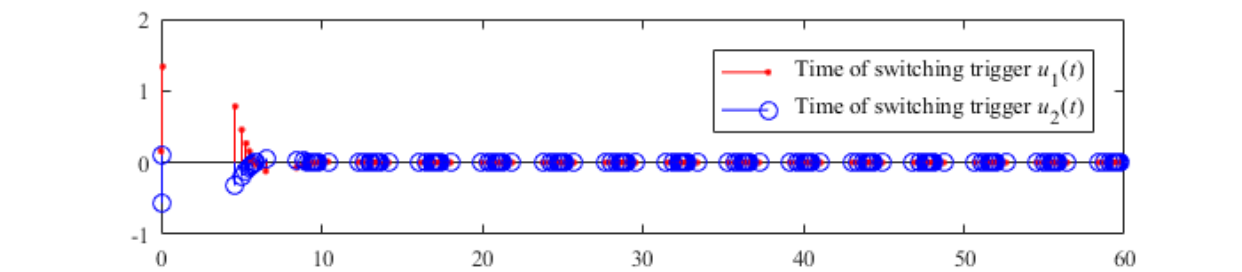}
		\caption{Triggered moment of the switching event-triggered controller $u(t)$ of slow system.}
		\label{Fig. 16}
	\end{figure}
		\begin{figure}[!t]
		\centering
		\includegraphics[width=0.35\textwidth]{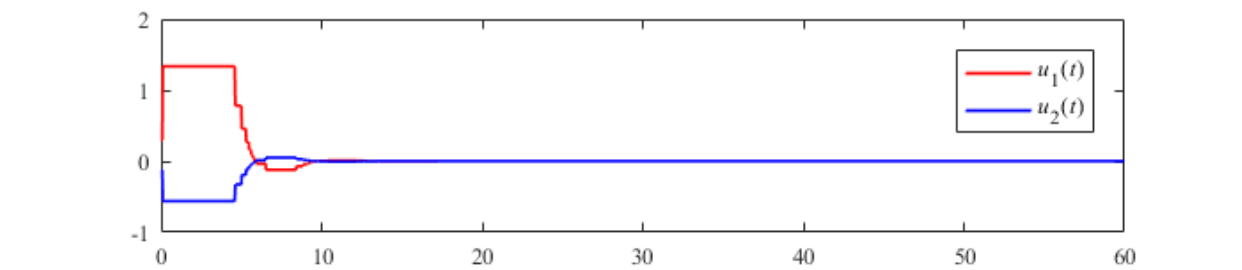}
		\caption{Trajectory of the switching event-triggered controller $u(t)$.}
		\label{Fig. 17}
	\end{figure}
	\begin{figure}[!t]
		\centering
		\includegraphics[width=0.35\textwidth]{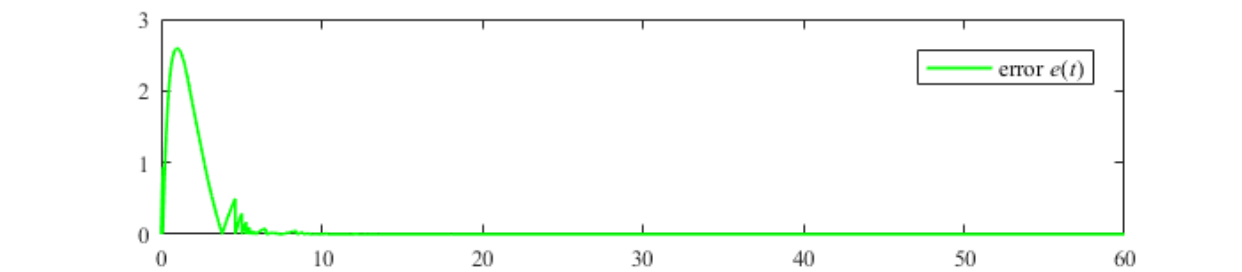}
		\caption{Triggered moment of the switching event-triggered controller $u(t)$ of slow system.}
		\label{Fig. 18}
	\end{figure}
\end{example}
\section{Conclusion}\label{s6}
In this article, the parabolic PDE systems with unknown nonlinearities is separated applying the Galerkin method. A finite-dimensional slow subsystem that faithfully captures the main characteristics of the original system can be constructed. Then, the MNN is used to model unknown nonlinear components. To ensure the stability of system, a switching event-triggered controller is designed and the corollary at absence of disturbances is presented. The $H_\infty$ performance control problem is then explored and the outcome of this issue is formulated as a BMI optimization problem. In conclusion, the simulation examples on the catalytic rod reaction model and traffic flow model show that the devised controller is more efficient at conserving control resources than static event-triggered controllers. The designed approach will be extended to the tracking issues of nonlinear PDE systems with unmatched disturbances in our future study.

\appendices

\section{Proof of Theorem 5}\label{B}
\begin{IEEEproof}
For the closed-system (49), the Lyapunov functional is as follows:	
\begin{equation*}
	\Gamma(t)=V(t)+V_f(t),
\end{equation*}
where $V(t)$ is from (20), $V_f(t)=q\xi_{f}^T(t)\xi_{f}(t)$, $q$ is a positive scalar. If (47) is holds, there exists a scalar $\kappa(\gamma_s)>0$ that it can be gotten:
\begin{equation*}
	\sum_{i=1}^{r}\sigma _{i}(\theta (t))\Phi_{5}<-\kappa(\gamma_s)I.			\end{equation*}
According to Lipschitz condition, there exist positive scalars $\kappa_1$, $\kappa_2$, $\kappa_3$, satisfy fact:
\begin{equation*}\left\{\begin{aligned}\|\Delta f_s\|&\leq\kappa_1\|\xi_f\|_{l^2},\\
		\|f_f(\xi_s,\xi_f)\|&\leq\kappa_2\|\xi_s(t)\|+\kappa_3\|\xi_f(t)\|_{l^2}.\\
		\end{aligned}\right.\end{equation*}
Otherwise, we have:
\begin{equation*}\begin{aligned}&\int_{\Omega}\bar{b}_{2}^{\mathrm{T}}(z)\bar{b}_{2}(z)dz= B_{2,s}^{\mathrm{T}}B_{2,s}+B_{2,f}^{\mathrm{T}}B_{2,f},\\&\int_{\Omega}\bar{b}_{1}^{\mathrm{T}}(z)\bar{b}_{1}(z)dz= B_{1,s}^{\mathrm{T}}B_{1,s}+B_{1,f}^{\mathrm{T}}B_{1,f},\\&\int_{\Omega}\bar{c}^{\mathrm{T}}(z)\bar{c}(z)dz= C_{s}^{\mathrm{T}}C_{s}+C_{f}^{\mathrm{T}}C_{f}.\end{aligned}\end{equation*}
Similar to (42)-(46), it can be derived:
\begin{equation*}\begin{aligned}
		&\ y_s^{T}(t)y_s(t)-\gamma_s^{2}{d_{1}(t)}^{T}d_{1}(t)+\dot{V}\\
		\leq&-\kappa(\gamma_s)\left\|\bar\xi_{2}(t)\right\|^2+2\vartheta_1 \kappa_1 \left\| \xi_s\right\|\left\| \xi_f\right\|,
\end{aligned}\end{equation*}	
where 
$\vartheta_1=\bar{\sigma}(N)=\sup_{0\neq \xi_s\in\mathbb{R}^m} \frac{\|N\xi_s\|}{\|\xi_s\|}\geq0.
$	

Define $\bar{\tau}=t/ \varepsilon$ and setting $ \varepsilon=0$, we can obtain the fast subsystem by (\ref{e7}):
$\frac{d\xi_f(\bar{\tau})}{d\bar{\tau}}=A_{f\varepsilon}\xi_f(\bar{\tau}),$
based on $\lambda_ {m+1}<0$ and the definition of $\varepsilon$, the fast subsystem has been found to be globally exponentially stable.

For any matrix $M_2$, $N_2$, $Y_1$, $Y_2$, the following inequalities will invariably be valid:
\begin{equation*}\begin{aligned}
		0&=2(\nu_1(t)N_2+\bar{\xi}_2(t)N_3)[\tau(t)\xi_s(t)-\tau(t)\xi_s(t-\tau(t))-\nu(t)],\\
		0&=2[\xi_s^T(t)Y_1+\xi_s^T(t-\tau(t))Y_2][y_s(t)-y_s(t_k)+e(t)]
		.\end{aligned}
\end{equation*}
Similarly, the following inequality can be deduced:
\begin{equation*}\begin{aligned}
		&\ y_f^{T}(t)y_f(t)-\gamma_f^{2}{d_{1}(t)}^{T}d_{1}(t)+\dot{V}_f\\
		&\leq2q\xi_f^T(t)A_f\xi_f(t)+2q\xi_f^T(t)\kappa_{2}\|\xi_s(t)\|\\
		&\ +2q\xi_f^T(t)\kappa_{3}\|\xi_f(t)\|+2q\xi_f^T(t)B_{1,f}d(t)\|\\
		&\ -2q\xi_f^T(t)B_{2,f}KC_s\xi_s(t)+2q\xi_f^T(t)\chi_{1}(t)B_{2,f}KC_s\tau(t)\nu_1(t)\\
		&\ -2(1-\chi_1(t))q\xi_f^T(t)B_{2,f}Ke(t)+\xi_f^T(t)C_f^TC_f\xi_f(t)\\
		&\ -\gamma_f^2\|d(t)\|^2+2(\nu_1^T(t)N_2+\bar{\xi}_2^T(t)N_3)\\
		&\ \times[\tau(t)\xi_s(t)-\tau(t)\xi_s(t-\tau(t))-\nu_1(t)]\\
		&\ +2[\xi_s^T(t)Y_1+\xi_s^T(t-\tau(t))Y_2][y_s(t)-y_s(t_k)+e(t)]\\
		&\ +\epsilon y_s^{T}(t)\Lambda  y_s(t)-(y_s(t)-y_s(t_{k}))^{T}\Lambda (y_s(t)-y_s(t_{k})).\\
\end{aligned}\end{equation*}
If $\chi(t)=1$, one has
\begin{equation*}\begin{aligned}	&\ y_f^{T}(t)y_f(t)-\gamma_f^{2}{d_{1}(t)}^{T}d_{1}(t)+\dot{V}_f\\
		&\leq2q\lambda_{m+1}\|\xi_f(t)\|^2+2q\kappa_2\|\xi_s(t)\|\|\xi_f(t)\|+2q\kappa_{3}\|\xi_f(t)\|^2\\
		&\ -2q\kappa_{4}\vartheta _2\vartheta _3\|\xi_s(t)\|\|\xi_f(t)\|+2q\kappa_{5}\|\xi_s(t)\|\|d(t)\|\\
		&\ +\kappa_{6}\|\xi_f(t)\|^2-\gamma_f^2\|d(t)\|^2+2q\vartheta _4\vartheta _5\kappa_{7}h\|\xi_f(t)\|\|\nu_1(t)\|\\
		&\ +2\vartheta _6 h\|\nu_1(t)\|\|\xi_s(t)\|-2\vartheta _7\|\nu_1(t)\|^2\\
		&\ +2\vartheta_8 h\|\bar{\xi}_2(t)\|\|\xi_s(t)\|-2\vartheta _9\|\xi_s(t)\|\|\nu_1(t)\|\\
		&\ +2\vartheta_5\vartheta_{10}\|\xi_s(t)\|^2-2\vartheta_3\vartheta _{11}\|\xi_s(t)\|\|\xi_s(t-\tau(t))\|\\
		&\ +2\vartheta _{10}\|\xi_s(t)\|\|e(t)\|+2\vartheta _5\vartheta _{12}\|\xi_s(t)\|\|\xi_s(t-\tau(t))\|\\
		&\ -2\vartheta_3\vartheta _{13}\|\xi_s(t-\tau(t))\|^2+2\vartheta _{12}\|\xi_s(t-\tau(t))\|\|e(t)\|\\
		&\ +\epsilon\Lambda\vartheta_{5}\|\xi_s(t)\|^2-\Lambda\|e(t)\|^2.
\end{aligned}\end{equation*}
where $\vartheta_{2}=(\underline{\sigma}(K^TK))^\frac{1}{2}=(\inf_{0\neq \xi_s\in\mathbb{R}^m} \frac{\|K^TK\xi_s\|}{\|\xi_s\|})^\frac{1}{2}\geq0$, $\vartheta_{3}=(\underline{\sigma}(C_s^T C_s))^\frac{1}{2}$, 
$\vartheta_{4}=(\bar{\sigma}(K^TK))^\frac{1}{2}$, $\vartheta_{5}=(\bar{\sigma}(C_s^T C_s))^\frac{1}{2}$, $\vartheta _6=\bar{\sigma}(N_2)$, $\vartheta _7=\underline{\sigma}(N_2)$, $\vartheta _8=\bar{\sigma}(N_3)$, $\vartheta _9=\underline{\sigma}(N_3)$,  $\vartheta _{10}=\bar{\sigma}(Y_1)$, $\vartheta _{11}=\underline{\sigma}(Y_1)$,  $\vartheta _{12}=\bar{\sigma}(Y_2)$, $\vartheta _{13}=\underline{\sigma}(Y_2)$,  $\kappa_{4}=(\underline{\sigma}(\int_{\Omega}\bar{b}_{2}^{\mathrm{T}}(z)\bar{b}_{2}(z)dz- B_{2,s}^{\mathrm{T}}B_{2,s}))^\frac{1}{2}, $ $\kappa_{5}=(\bar{\sigma}(\int_{\Omega}\bar{b}_{1}^{\mathrm{T}}(z)\bar{b}_{1}(z)dz- B_{1,s}^{\mathrm{T}}B_{1,s}))^\frac{1}{2}, $ $\kappa_{6}=(\bar{\sigma}(\int_{\Omega}\bar{c}^{\mathrm{T}}(z)\bar{c}(z)dz-C_{s}^{\mathrm{T}}C_{s}))^\frac{1}{2}$, $\kappa_{7}=(\bar{\sigma}(\int_{\Omega}\bar{b}_{2}^{\mathrm{T}}(z)\bar{b}_{2}(z)dz- B_{2,s}^{\mathrm{T}}B_{2,s}))^\frac{1}{2} $.

Let $\gamma=\gamma_s+\gamma_f$, then it can be obtained:
\begin{equation*}\begin{aligned}
		& \dot{\Gamma}(t)+y^{T}(t)y(t)-\gamma^{2}{d_{1}(t)}^{T}d_{1}(t)
		\leq \bar{\xi}_3^T(t)\tilde{\Lambda}_1\bar{\xi}_3(t),
\end{aligned}\end{equation*}
where $\bar{\xi}_3^T(t)=[\left\|\bar{\xi}_2(t)\right\|,\left\|\xi_s(t)\right\|,\left\|\xi_f(t)\right\|$, $\left\|e(t)\right\|,\left\|d(t)\right\|,\\\left\|\nu_1(t)\right\|,\left\|\xi_s(t-\tau(t))\right\|]$,
\begin{equation*}
	\tilde{\Lambda}_1=\begin{bmatrix}
		\Gamma_{11}&\Gamma_{12}&0&0&0&0&0\\
		\ast&\Gamma_{22}&\Gamma_{23}&\vartheta_{10}&\Gamma_{25}&\Gamma_{26}&\Gamma_{27}\\
		\ast&\ast&\Gamma_{33}&0&0&0&\Gamma_{37}\\
		\ast&\ast&\ast&-\Lambda&0&0&\vartheta_{12}\\
		\ast&\ast&\ast&\ast&\Gamma_{55}&0&0\\
		\ast&\ast&\ast&\ast&\ast&\Gamma_{66}&0\\
		\ast&\ast&\ast&\ast&\ast&\ast&\Gamma_{77}\\
	\end{bmatrix}
\end{equation*}
\begin{equation*}
	\begin{aligned}
		\Gamma_{11}&=-\kappa(\gamma_s),\ \Gamma_{12}=\vartheta_8h,\\
		\Gamma_{22}&=2\vartheta_5\vartheta_{10}+\epsilon\Lambda\vartheta_{5},\\ \Gamma_{23}&=\vartheta_5\kappa_{1}+q\kappa_{2}-q\kappa_{4}\vartheta_2\vartheta_3,\\
		\Gamma_{25}&=q\kappa_{5},\ \Gamma_{26}=-\vartheta_9,\\
		\Gamma_{27}&=\vartheta_6h-\vartheta_3\vartheta_{11}+\vartheta_5\vartheta_{12},\\
		\Gamma_{33}&=2q\lambda_{m+1}+2q\kappa_{3}+\kappa_{6},\\
		\Gamma_{37}&=q\vartheta_4\vartheta_5\kappa_{7}h,\ 
		\Gamma_{55}=-\gamma_f^2,\\
		\Gamma_{66}&=-2\vartheta_7,\ \Gamma_{77}=-2\vartheta_3\vartheta_{13}.\\
	\end{aligned}
\end{equation*}
Let
\begin{equation*}
	\begin{aligned}
		\Xi_3&=(\Lambda\Gamma_{77}+\vartheta_{12}^2)[(	-\Gamma_{11}\Gamma_{22}+\Gamma_{12}^2)\Gamma_{55}\Gamma_{66}\\
		&\ +	\Gamma_{11}\Gamma_{25}^2\Gamma_{66}+\Gamma_{26}^2\Gamma_{55}],\\
		\Xi_4&=-	\Gamma_{11}\Gamma_{55}\Gamma_{66}[\vartheta_{10}^2\Gamma_{77} -(\vartheta_{10}\vartheta_{12}+\Lambda\Gamma_{27})\Gamma_{27}],\\
		\Xi_5&=\Gamma_{55}[\Lambda\Gamma_{66}((-	\Gamma_{11}\Gamma_{22}+\Gamma_{12}^2)\Gamma_{37}^2-	\Gamma_{11}\Gamma_{23}^2\Gamma_{77}\\
		&\ +\Gamma_{23}\Gamma_{27}\Gamma_{37})+\Gamma_{26}^2\Gamma_{37}^2\Lambda].\\
	\end{aligned}
\end{equation*}
Based on the fact $|\lambda_L|=-\varepsilon\lambda_{m+1}$, and define
\begin{equation*}
	\varepsilon_1^*=\frac{2q(\Xi_3+\Xi_4)|\lambda_L|}{(2q\kappa_{3}+\kappa_{6})(\Xi_3+\Xi_4)-\Xi_5}.
\end{equation*}
If $\chi(t)=0$,, similar to the above analysis, one holds:
\begin{equation*}\begin{aligned}	&\ \ \ \  y_f^{T}(t)y_f(t)-\gamma_f^{2}{d_{1}(t)}^{T}d_{1}(t)+\dot{V}_f\\
		&\leq2q\lambda_{m+1}\|\xi_f(t)\|^2+2q\kappa_2\|\xi_s(t)\|\|\xi_f(t)\|+2q\kappa_{3}\|\xi_f(t)\|^2\\
		&\ -2q\kappa_{4}\vartheta _2\vartheta _3\|\xi_s(t)\|\|\xi_f(t)\|+2q\kappa_{5}\|\xi_s(t)\|\|d(t)\|\\
		&\ +\kappa_{6}\|\xi_f(t)\|^2-\gamma_f^2\|d(t)\|^2-2q\vartheta_2\kappa_{4}\|\xi_f(t)\|\|e(t)\|\\
		&\ +2\vartheta_6 h\|\nu_1(t)\|\|\xi_s(t)\|-2\vartheta _7\|\nu_1(t)\|^2\\
		&\ +2\vartheta_8 h\|\bar{\xi}_2(t)\|\|\xi_s(t)\|-2\vartheta _9\|\xi_s(t)\|\|\nu_1(t)\|\\
		&\ +2\vartheta _5\vartheta_{10}\|\xi_s(t)\|^2-2\vartheta_3\vartheta _{11}\|\xi_s(t)\|\|\xi_s(t-\tau(t))\|\\
		&\ +2\vartheta _{10}\|\xi_s(t)\|\|e(t)\|+2\vartheta _5\vartheta _{12}\|\xi_s(t)\|\|\xi_s(t-\tau(t))\|\\
		&\ -2\vartheta _3\vartheta _{13}\|\xi_s(t-\tau(t))\|^2+2\vartheta _{12}\|\xi_s(t-\tau(t))\|\|e(t)\|\\
		&\ +\epsilon\Lambda\vartheta_{5}\|\xi_s(t)\|^2-\Lambda\|e(t)\|^2.
\end{aligned}\end{equation*}
Then it can be derived:
\begin{equation*}\begin{aligned}
		& \dot{\Gamma}(t)+y^{T}(t)y(t)-\gamma^{2}{d_{1}(t)}^{T}d_{1}(t)
		\leq \bar{\xi}_3^T(t)\tilde{\Lambda}_2\bar{\xi}_3(t),
\end{aligned}\end{equation*}
where	\begin{equation*}
	\begin{aligned}
		\tilde{\Lambda}_1&=\begin{bmatrix}
			\Gamma_{11}&\Gamma_{12}&0&0&0&0&0\\
			\ast&\Gamma_{22}&\Gamma_{23}&\vartheta_{10}&\Gamma_{25}&\Gamma_{26}&\Gamma_{27}\\
			\ast&\ast&\Gamma_{33}&\Gamma_{34}&0&0&0\\
			\ast&\ast&\ast&-\Lambda&0&0&\vartheta_{12}\\
			\ast&\ast&\ast&\ast&\Gamma_{55}&0&0\\
			\ast&\ast&\ast&\ast&\ast&\Gamma_{66}&0\\
			\ast&\ast&\ast&\ast&\ast&\ast&\Gamma_{77}\\
		\end{bmatrix},\\
		\Gamma_{34}&=-q\vartheta_2\kappa_{4}.\\
	\end{aligned}
\end{equation*}
Let \begin{equation*}
	\begin{aligned}
		\Xi_6&=-\vartheta_{12}\Gamma_{11}\Gamma_{25}^2\Gamma_{34}\Gamma_{37}\Gamma_{66},\\
		\Xi_7&=(\vartheta_{12}\Gamma_{27}-\vartheta_{10}\Gamma_{77})\Gamma_{23}\Gamma_{34}-(\vartheta_{12}^2+\Gamma_{77}\Lambda)\Gamma_{23}^2,\\
		\Xi_8&=(-\Gamma_{11}\Gamma_{25}^2\Gamma_{66}+\Gamma_{26}^2\Gamma_{55}-\Gamma_{22}\Gamma_{55}\Gamma_{66}\\
		&\ +\Gamma_{12}^2\Gamma_{55}\Gamma_{66})(\Lambda\Gamma_{77}+\vartheta_{12}^2),\\
		\Xi_9&=\Gamma_{55}\Gamma_{66}(\vartheta_{10}\vartheta_{12}\Gamma_{27}+\Lambda\Gamma_{27}^2
		-\vartheta_{10}^2\Gamma_{77}+\vartheta_{10}\vartheta_{12}\Gamma_{27}).
	\end{aligned}
\end{equation*}
Similar to the above process, we define
\begin{equation*}
	\varepsilon_2^*=\frac{2q(\Xi_8+\Xi_9)|\lambda_L|}{(2q\kappa_{3}+\kappa_{6})(\Xi_8+\Xi_9)-\Xi_6-\Xi_7}.
\end{equation*}
We can obtain if $\varepsilon\in (0,\varepsilon^*)$:
\begin{equation*}\begin{aligned}
		& y^{T}(t)y(t)-\gamma^{2}{d_{1}(t)}^{T}d_{1}(t)
		\leq -\dot{\Gamma}(t),
\end{aligned}\end{equation*}
where $\varepsilon^*=\min\{\varepsilon_1^*,\ \varepsilon_2^*\}$.

Integrating both sides of the aforementioned inequality simultaneously yields:
\begin{equation*}
	\begin{aligned}
		& \int_{0}^{T}y^{T}(t)y(t)dt-\int_{0}^{T}\gamma^{2}{d_{1}(t)}^{T}d_{1}(t)dt
		\leq 0,
	\end{aligned}
\end{equation*}
then the closed-system of PDE system (1)-(3) satisfies the $H_\infty$ performance.

Furthermore, with reference to Theorem 1 and the aforementioned proof process, it is straightforward for us to analyze and deduce the existence of a positive number $\kappa_0$, such that:
\begin{equation*}
	\dot{\Gamma}(t)\leq -\kappa_0 \Gamma(t)+3D_{1}^{2}.
\end{equation*}
According to Lemma 1, it can be deduced:
\begin{equation*}
	0\leq \Gamma(t)\leq \Gamma(0)\mathrm{e}^{-\kappa_0t}+\frac{ 3D_{1}^{2}}{\kappa_0}(1-e^{-\kappa_0t}),
\end{equation*}
Thus, the closed-loop system of PDE systems (1)-(3) achieves SGUUB and satisfies the $H_\infty$ performance stable under the action of the controller (\ref{e16}). 

And when $d(t)\equiv 0$, it is satisfied that $\Gamma(t)\leq \Gamma(0)\mathrm{e}^{-\kappa_0t}$, then we can get:
\begin{equation*}
\left\|\xi_s(t)\right\|^2+\left\|\xi_f(t)\right\|^2\leq \alpha_9 (\left\|\xi_s(0)\right\|^2+\left\|\xi_f(0)\right\|^2)\mathrm{e}^{-\kappa_0t},
	\end{equation*}
	where $\alpha_9=\frac{\alpha_7}{\alpha_8}$, $\alpha_7=\max\{\alpha_2,\lambda_{\max}(q)\}$, $\alpha_8=\min\{\alpha_5,\lambda_{\min}(q)\}$.
	Thus it can be obtained 
	\begin{equation*}
	\left\|\bar{\xi}(t)\right\|_{2,\Omega}\leq \sqrt{\alpha_9} \left\|\bar{\xi}(0)\right\|_{2,\Omega}\mathrm{e}^{-0.5\kappa_0t},
	\end{equation*}
In conclusion, the closed-loop system of PDE systems (1)-(3) is exponentially stable when $d(t)\equiv 0$. The proof is complete.
\end{IEEEproof}
\bibliographystyle{IEEEtran}

\end{document}